\documentclass[11pt,english]{article}
\usepackage[T1]{fontenc}
\usepackage[latin1]{inputenc}
\usepackage{babel}
\usepackage{setspace}
\usepackage{indentfirst}

\makeatletter

\usepackage{verbatim}

\makeatother

\usepackage{amsmath}
\usepackage{amssymb}
\usepackage{stmaryrd}
\usepackage{wasysym}
\usepackage{pifont}
\usepackage{graphicx}
\usepackage[all]{xy}
\usepackage{amsthm}
\theoremstyle{plain}% default
\newtheorem{thm}{Theorem}[section]
\newtheorem{lem}[thm]{Lemma}
\newtheorem{prop}[thm]{Proposition}
\newtheorem{cor}[thm]{Corollary}
\newtheorem{defn}[thm]{Definition}
\theoremstyle{definition}
\newtheorem{example}[thm]{Example}% single example
\theoremstyle{remark}
\newtheorem{rmk}[thm]{Remark}% produces a remark with no number

\DeclareMathOperator{\ad}{ad}
\DeclareMathOperator{\Ad}{Ad}
\DeclareMathOperator{\AD}{\mathbf{Ad}}

\DeclareMathOperator{\iso}{iso}
\DeclareMathOperator{\out}{out}
\DeclareMathOperator{\loc}{loc}
\DeclareMathOperator{\Iso}{Iso}
\DeclareMathOperator{\Inn}{Inn}
\DeclareMathOperator{\Out}{Out}
\DeclareMathOperator{\Aut}{Aut}
\DeclareMathOperator{\ev}{ev}
\DeclareMathOperator{\idn}{id}
\DeclareMathOperator{\hol}{hol}
\DeclareMathOperator{\Ver}{Ver}
\DeclareMathOperator{\Hor}{Hor}
\DeclareMathOperator{\Lie}{Lie}
\DeclareMathOperator{\End}{End}
\DeclareMathOperator{\Hol}{Hol}
\DeclareMathOperator{\hor}{hor}
\DeclareMathOperator{\obs}{obs}
\DeclareMathOperator{\supp}{supp}
\DeclareMathOperator{\pr}{pr}

\renewcommand{\epsilon}{\varepsilon}
\newcommand{\lie}[2]{[#1,#2]}
\newcommand{\gendex}[2]{\left\{#1\right\}_{#2}}
\newcommand{\genrel}[2]{\left\{#1|#2\right\}}
\newcommand{\sect}[1]{\Gamma(#1)}
\newcommand{\derlie}[1]{\mathcal{L}_{#1}}
\newcommand{\ii}{\mathbin{\vrule width1.5ex height.4pt\vrule height1.5ex}}
\newcommand{\ip}{\ii}
\newcommand{\cty}{\mathop{C^\infty}}

\newcommand{\rond}{\smalcirc} \newcommand{\smalcirc}{\mbox{\tiny{$\circ$}}}
\newcommand{\inv}{^{-1}}
\newcommand{\half}{\tfrac{1}{2}}
\newcommand{\simplicial}{_{\scriptscriptstyle\bullet}}
\newcommand{\degree}{^{\bullet}}
\newcommand{\isomorphism}{\simeq}
\newcommand{\mfg}{\mathfrak{g}}
\newcommand{\toto}{\rightrightarrows}
\newcommand{\lto}{\longrightarrow}
\newcommand{\R}{\mathbb{R}}
\newcommand{\N}{\mathbb{N}}

\newcommand{\bt}{\mathbf{t}}
\newcommand{\bs}{\mathbf{s}}
\newcommand{\kernel}{\mathcal{K}}
\newcommand{\liekernel}{\mathfrak{K}}
\newcommand{\liecenter}{Z(\liekernel)}
\newcommand{\hkernel}{{H^{\kernel}}}

\newcommand{\hiso}{H^{\iso}}
\newcommand{\hout}{H^{\out}}
\newcommand{\alphakernel}{\alpha^{\kernel}}
\newcommand{\alphaiso}{\alpha^{\iso}}
\newcommand{\alphaout}{\alpha^{\out}}
\newcommand{\omegakernel}{\omega^{\kernel}}
\newcommand{\omegaiso}{\omega^{\iso}}
\newcommand{\omegaout}{\omega^{\out}}
\newcommand{\wtr}{\blacktriangleright}
\newcommand{\wtl}{\blacktriangleleft}
\newcommand{\btr}{\vartriangleright}
\newcommand{\btl}{\vartriangleleft}
\newcommand{\ptr}{\partial^{\btr}}
\newcommand{\ptl}{\partial^{\btl}}
\newcommand{\pt}{\partial}
\newcommand{\pcoho}[4]{H_{\ptr }^{#1,#2}(#3\simplicial\to #4\simplicial)}
\newcommand{\hcoho}[4]{H_{\hor}^{#1,#2}(#3\simplicial\to #4\simplicial)}
\newcommand{\pcohop}[5]{H_{\ptr }^{#1,#2}(#3\simplicial[#5]\to #4\simplicial[#5])}
\newcommand{\hcohop}[5]{H_{\hor}^{#1,#2}(#3\simplicial[#5]\to #4\simplicial[#5])}
\newcommand{\EEE}[2]{\mathcal{E}^{#1,#2}}
\newcommand{\FFF}[2]{\mathcal{F}^{#1,#2}}
\newcommand{\GGG}[2]{\mathcal{G}^{#1,#2}}
\newcommand{\rr}[2]{r_{#1}^{#2}}
\newcommand{\EE}{\mathcal{E}}
\newcommand{\FF}{\mathcal{F}}
\newcommand{\GG}{\mathcal{G}}

\newcommand{\sigt}{\tilde{\sigma}}
\newcommand{\etab}{\bar{\eta}}
\newcommand{\etau}{\underline{\eta}}
\newcommand{\ADb}{\overline{\AD}}
\newcommand{\ADu}{\underline{\AD}}
\newcommand{\gammab}{\bar{\gamma}}
\newcommand{\gammadot}{\overset{\scriptscriptstyle\bullet}{\gamma}}
\newcommand{\lambdab}{\underline{\lambda}}

\def\dar[#1]{\ar@<2pt>[#1]\ar@<-2pt>[#1]}

\newcommand{\MM}{\mathcal{M}}
\newcommand{\XX}{\mathfrak{X}}
\newcommand{\UU}{\mathfrak{U}}
\newcommand{\VV}{\mathfrak{V}}

\newcommand{\bbt}{\tau}
\newcommand{\covder}{ D }
\newcommand{\covderPB}{ d^{\nabla^{\bt}} }
\newcommand{\gm}{\Gamma}
  
\begin{document}

\title{Non Abelian Differentiable Gerbes}
\author{Camille \textsc{Laurent-Gengoux} \\ D\'epartement de math\'ematiques \\ Universit\'e 
de Poitiers \\ 86962 Futuroscope-Chasseneuil, France \\ 
\texttt{laurent@math.univ-poitiers.fr} \and
 Mathieu \textsc{Sti\'enon} \thanks{Francqui fellow of the Belgian American Educational 
Foundation}\\ Department of Mathematics \\ Pennsylvania State University \\ University Park, 
PA 16802, USA \\ \texttt{stienon@math.psu.edu} \and
 Ping \textsc{Xu}
\thanks{Research partially supported by NSF
       grants DMS03-06665, DMS-0605725  and NSA grant 03G-142. } \\ 
Department of Mathematics \\ Pennsylvania State University \\ University Park,
 PA 16802, USA \\  \texttt{ping@math.psu.edu} }
\date{
%%%
\bigskip \small{\textit{Dedicated to Jean-Luc Brylinski at the occasion of his 55th birthday}} 
%%%
}
%\date{\texttt{file: \jobname.tex}}
\maketitle

\begin{abstract}
We study non-abelian differentiable gerbes over
stacks using the theory of Lie groupoids.
More precisely, we  develop the theory 
of   connections on Lie groupoid $G$-extensions, which we call ``connections on gerbes'', 
and study the induced connections on various associated
bundles. We also  prove analogues of the Bianchi identities.
In particular, we develop a cohomology theory which
measures the existence of connections and curvings
for $G$-gerbes over stacks.
We also introduce $G$-central extensions of groupoids,
generalizing the standard groupoid $S^1$-central extensions.
 As an example, we apply our theory to study the differential geometry
of $G$-gerbes over a manifold.
\end{abstract}

%%%%%%%%%%%% COMMENT THIS WHEN PRODUCING FINAL VERSION %%%%%%%%%%%%%%%%%%%%%%
%\pagestyle{fancy}
%\fancyhead{}
%\renewcommand\headrulewidth{0pt}
%\fancyfoot{\texttt{\jobname.tex}\hfill\thepage}
%%%%%%%%%%%%%%%%%%%%%%%%%%%%%%%%%%%%%%%%%%%%%%%%%%%%%%%%%%%%%%%%%%%%%%%%%%%%%

\tableofcontents

\section{Introduction}
The general theory of gerbes over a site was developed by Giraud \cite{Giraud}.
 Its main motivation is the study of non-abelian
cohomology theory. Indeed, Giraud's theory shows, in particular, that 
$S^1$-bound gerbes over a manifold are in one-one
correspondence with the third integer-coefficients cohomology group of 
the manifold. Due to this fact, gerbes are often referred
to, in the literature, as geometric realizations of degree 3 integer cohomology
classes. Recently, string theory \cite{GN, Witten, Witten1} fostered the study of 
the differential geometry of gerbes. The differential geometry of $S^1$-bound gerbes over manifolds 
--- Murray named them bundle gerbes \cite{Murray} ---
was investigated by many authors among whom Brylinski \cite{brylinski1}, Murray \cite{Murray},
 Chatterjee \cite{Chatterjee} and
Hitchin \cite{Hitchin}.
A theory of connections, curvings and $3$-curvatures was developed,
which was used to define the characteristic classes, called
Dixmier Douday classes, of bundle gerbes \cite{Murray}.

However, the differential geometry of the more general non abelian gerbes
is a subtler question. For bundle gerbes, a connection
is defined to be a connection $1$-form for the usual
principal $S^1$-bundle satisfying some additional
property \cite{Murray, Hitchin}. However,
as non abelian $G$-gerbes can no longer be
considered principal $G$-bundles, 
the obvious generalization won't work,
and it is not clear what a connection on
a non abelian gerbe should be. 

Breen-Messing were the first to study the differential geometry of 
non abelian gerbes. Their important approach relies on the use of 
Kocks' synthetic geometry \cite{Kock} as a tool to study the differential geometry
of non abelian gerbes from an algebraic geometry perspective.
They introduced the concepts of connections, curvings and 3-curvatures
and obtained various identities between them.
Later on, Aschieri-Cantini-Jur\v{c}o investigated a non abelian analogue 
of bundle gerbes
 from a more physical point of view \cite{ACJ}. 
Yet another interesting approach pioneered by Baez-Schreiber 
\cite{Baez:2002, BS, BaezSchreiber:2004, Schreiber:2004e}
is to use higher gauge theory \cite{Bartels:2004, Gajer:1999, Pfeiffer:2003,GirelliPfeiffer:2004}. 
As was shown by Breen \cite{Breen, Br4}, a $G$-gerbe is equivalent to
a 2-group bundle in the sense of Dedecker \cite{Dedecker},
where the 2-group is the one corresponding to
the crossed module $G\to\Aut(G)$ (see \cite{GS:08} for
a geometrical construction).
Hence, a ``connection'' on a $G$-gerbe should correspond 
to a ``connection'' on the 2-group principal bundle.  
The latter is the viewpoint adopted by Baez-Schreiber \cite{BS, BaezSchreiber:2004}. See also
\cite{MackaayPicken:2000} for the case of bundle gerbes. 

The aim of this paper is to propose an approach based on Lie groupoid extensions 
to study geometry of non abelian differentiable gerbes.
It is based on the theory of differentiable stacks
developed in \cite{BX, BX1} (see also \cite{Metzler}).
Roughly speaking, differentiable stacks are Lie
groupoids \emph{up to Morita equivalence}. Any Lie groupoid $\gm$ defines a
differentiable stack: the differentiable stack $\XX_\gm$ of $\gm$-torsors. 
Two differentiable stacks $\XX_\gm$ and $\XX_\gm'$ are isomorphic if, 
and only if, the Lie groupoids $\gm$ and $\gm'$ are Morita equivalent. 
In a certain sense, Lie groupoids are like ``local charts'' on a differentiable stack. 
The relation between differentiable gerbes and Lie
groupoids is described by the following table, whose right hand side
lists the main topics investigated in this paper. 
\begin{center}
\begin{tabular}{|c|c|} 
\hline \textsc{Stack language} & \textsc{Groupoid Language} \\ 
\hline differentiable  stacks & Morita equivalent Lie groupoids \\ 
 gerbes over stacks & Morita equivalent groupoid extensions \\ 
 $G$-gerbes over stacks & Morita equivalent groupoid $G$-extensions \\
 $G$-bound gerbes over stacks & Morita equiv. groupoid $G$-extensions with trivial band \\
\hline\end{tabular}
\end{center}

Our approach has the advantage of avoiding abstract stacks 
by working directly with Lie groupoids, which are much more down-to-earth.
Moreover it allows us to use the usual techniques of differential geometry and formulate
our results in a global fashion without resorting to local charts.
Such a viewpoint has already been taken by a few authors. For example, by Murray \cite{Murray}
and Murray-Stevenson \cite{MS} in their study of bundle gerbes and
in the proof of their theorem concerning stable equivalence classes. 
For equivariant bundle gerbes, see \cite{GN, Meinrenken, Stienon}.
In \cite{BX, BX1} and  \cite{Bunke},
$S^1$-bound gerbes over differentiable stacks
were studied in terms of
$S^1$-central extensions of  Lie groupoids,  and
the characteristic classes of
gerbes, i.e. the Dixmier-Douady classes, were
introduced in terms of connection-like data.
This perspective is also central in a series of work of Moerdijk \cite{Moer:regular, Moerdijk}; 
in particular, \cite{Moer:regular} explores the deep connection 
between the classification of extensions of regular 
Lie groupoids, Giraud's non abelian cohomology \cite{Giraud} 
and subsequent work of Breen \cite{Breen2}.

In this paper, we develop a theory of connections on groupoid extensions,
which we call ``connections on gerbes''.
A connection on a groupoid extension $X_1\to Y_1 \toto M$ is an Ehresmann connection on the
fiber bundle $\phi: X_1\to Y_1$, which is compatible with the groupoid structure on $X_1\toto
M$. More precisely, a connection on a Lie groupoid extension
 $X_1 \xrightarrow{\phi} Y_1 \rightrightarrows M$ 
is a horizontal distribution $H$ on $X_1 \xrightarrow{\phi} Y_1$ which is 
also a subgroupoid of the tangent groupoid $TX_1\rightrightarrows TM$. 
For $S^1$-central extensions, such a connection is automatically 
a connection for the principal $S^1$-bundle $X_1 \xrightarrow{\phi} Y_1$, 
i.e. it is invariant under the action of $S^1$. 
This is in agreement with the definition of connections in 
\cite{BX, BX1, WeinsteinX:prequantization},
 which coincides with the one used
 by Brylinski (who uses the name ``connective structures'') \cite{brylinski1},
Murray \cite{Murray} and Hitchin \cite{Hitchin}.
However, in general, although $\phi: X_1\to Y_1$
is a bi-principal $\kernel$-$\kernel$-bundle, where $\kernel$ is the kernel of $\phi$, 
the horizontal distribution defining a groupoid extension connection is not invariant with respect to either
$\kernel$-action of the principal bundle.
We also introduce the notion of curvings and
$3$-curvatures for $G$-extensions.
Finally, a cohomology theory is developed
 that measures the obstructions to the existence of connections and curvings. 

We now describe the contents of this paper in more details.

Section 2 recalls some basic definitions and notions central to this paper. 
These include groupoid
extensions, Morita equivalence of Lie groupoid extensions, and generalized morphisms of Lie
groupoid extensions. Given a groupoid extension $X_1
\xrightarrow{\phi} Y_1\toto M$, let $\kernel=\ker \phi$ be its kernel. The outer action is a
canonical groupoid morphism from 
$Y_1\toto M$ to $\Out (\kernel, \kernel )\toto M$,
which  is essential to the definition of the band of a groupoid extension. 

Section 3 is devoted to the
study of $G$-gerbes, or more precisely, groupoid $G$-extensions. The notion of band is
introduced, which is a principal $\Out(G)$-bundle (or an $\Out(G)$-torsor) over the groupoid
$Y_1\toto M$. Groupoid $G$-extensions with trivial band --- they are called $G$-bound extensions in the paper --- are emphasized and related to the so called groupoid
$G$-central extensions, a natural generalization of groupoid $S^1$-central extensions
\cite{WeinsteinX:prequantization, BX, BX1}. A one-one correspondence between groupoid $G$-central extensions 
and $Z(G)$-central extensions is proved. 
As a  special case, we consider $G$-gerbes over a manifold, i.e. groupoid
 $G$-extensions of a \v{C}ech groupoid,
where
 the non-abelian \v{C}ech 2-cocycle conditions as in \cite{BreenMessing}
arise naturally. In particular, we
derive, by a direct argument, the following non-trivial theorem of Giraud:  the isomorphism
classes of $G$-bound gerbes over a manifold $N$ are in one-one correspondence with $H^2(N,
Z(G))$. 

Section 4 is devoted to the study of connections on groupoid extensions. As usual, we
give three equivalent definitions: horizontal distributions, parallel transportations, and
connection 1-forms. We then discuss the induced connections on the bundle of groups $\kernel\to M$
and its corresponding bundle of Lie algebras $\mathfrak{K}\to M$.
We also  prove analogues of the Bianchi identities.

In Section 5, we prove that a connection on a groupoid extension
induces a canonical connection on its band. 
We discuss the relations between the curvatures of various induced
connections, and in particular we derive an important formula which expresses the curvature
on the band in terms of the curvature on the group bundle $\kernel\to M$.

In the last Section, we develop a 
 cohomology theory, called horizontal cohomology in the paper, which captures
the obstruction to the existence of connections and curvings. In a certain sense, this is a
generalization of smooth Deligne cohomology \cite{brylinski1} to the non-abelian context. This 
cohomology is shown to be invariant under Morita equivalence, and hence depends only on the
underlying gerbe. We then compute explicitly the horizontal cohomology groups for $G$-gerbes
over a manifold. As a consequence we show that connections always exist on any $G$-gerbe over
a manifold. Finally we introduce flat $G$-gerbes and prove that in the case of central
extensions they are in one-one correspondence with flat $Z(G)$-gerbes. 

After the present paper was completed, we learned that 
non-abelian gerbes were also investigated by Stevenson \cite{Stevenson} from a different viewpoint.
The relation between our approach and Breen-Messing's
\cite{BreenMessing} is recently explored in \cite{BC}.
An explicit correspondence between $G$-gerbes and
$[G\to\Aut(G)]$-bundles is constructed geometrically in
\cite{GS:08}.
Non-abelian gerbes have also appeared in
many recent works in deformation quantization \cite{Tsygan, Tsygan1, Kashiwara, KS1, KS2, K2, AP, PS, Y1, Y2}.
We hope that our theory will be useful in understanding the geometry underlying these works.

\paragraph{Acknowledgments}
We would like to thank several institutions for their hospitality while work on this project
was being done, among which Penn State University (Laurent-Gengoux) and Universit\'e Pierre et Marie
Curie, Paris 6 (Xu). We also wish to thank many people for useful discussions
 and comments: Lawrence Breen,
Marius Crainic, Gr\'egory Ginot,  Andr\'e Haefliger,
Jim Stasheff, Jean-Louis Tu and Alan Weinstein to name a few, 
and especially Kai Behrend, who was participating in this project at its earlier stage.

\section{Differentiable  gerbes as groupoid extensions}
%groupoid multiplication : $X_2=X_1\times_{{\bt},M,{\bs}}X_1\to X_1$

\subsection{Lie groupoid extensions}
\label{sec:2.1}

The purpose of this section is 
to set up the basic notions of 
differentiable gerbes in terms of
Lie groupoid extensions. Some materials here
might be standard for experts \cite{Breen, Moer:regular}.

\begin{defn}
A Lie \textbf{groupoid extension} is a morphism of Lie groupoids
$$ \xymatrix{
X_1 \ar[r]^{\phi} \ar@<0.5ex>[d] \ar@<-0.5ex>[d] & Y_1 \ar@<0.5ex>[d] \ar@<-0.5ex>[d] \\
 M \ar@2{-}[r] & M }$$
where $\phi$ is  a fibration.
\end{defn}
For short, we denote a groupoid extension simply by $
X_{\simplicial}\xrightarrow{\phi}Y_{\simplicial}$ or $X_1\xrightarrow{\phi}Y_1\toto M$.

The \textbf{kernel}  of a
Lie  groupoid extension $X_{\simplicial}\xrightarrow{\phi}Y_{\simplicial}$ is 
the preimage $\kernel$ (or $\ker\phi$) in $X_1$ of the unit space 
of $Y\simplicial$. 
It is obviously a subgroupoid of $X_{\simplicial}$ on which the source and target maps coincide making it into a bundle of groups over $M$.
Alternatively, we can define a 
Lie groupoid extension as a short exact sequence of Lie groupoids 
\begin{equation}
\label{eq:ext}
 \xymatrix{ \kernel \ar[r] \ar@<-0.5ex>[d]_{\pi} \ar@<0.5ex>[d]^{\pi} & X_1 \ar[r]^{\phi} \ar@<0.5ex>[d]^t \ar@<-0.5ex>[d]_s & Y_1 \ar@<0.5ex>[d]^t \ar@<-0.5ex>[d]_s \\ M \ar@2{-}[r] & M \ar@2{-}[r] & M } 
\end{equation}
where $\kernel  \toto M$ is a
Lie groupoid whose target coincides with  the source, i.e.
a  bundle of Lie  groups.

There are a few important features that deserve to be pointed out.
First $\kernel  \toto M$ is a  subgroupoid of $X_1\toto M$. Therefore we can multiply an element
of $X_1$ by a composable element of $\kernel$ from both the left and the
right. These two actions commute and the quotient space by any of
these actions is $Y_1$. Thus $X_1\to Y_1$ is a
$\kernel$-$\kernel$ principal bibundle.

Second, the groupoid $X_1\rightrightarrows M$ acts on $\kernel \to M$ by
conjugation. More precisely, any $x\in X_1$ induces a group isomorphism
$$\AD_x: \kernel _{{\bt}(x)} \to \kernel _{ {\bs}(x)} ; g\mapsto  x\cdot g\cdot x\inv  ,$$
where the multiplication on the right hand side stands for
the product on $X_1$.

\subsection{Morita equivalence of Lie groupoid extensions}

Let $\Gamma_1\rightrightarrows \Gamma_0$ be a Lie groupoid,
and $J:P_0\to \Gamma_0$ a surjective submersion.
 Let $P_1$ denote the fibered product
$P_0\times_{J,\Gamma_0,{\bs}}\Gamma_1\times_{{\bt},\Gamma_0,J}P_0$. 
Then $P_1\rightrightarrows P_0$ has a natural structure of
 Lie groupoid with structure maps $ s :P_1\to P_0: (p,x,q)\mapsto p$,
$t:P_1\to P_0: (p,x,q)\mapsto q$, and
$m: P_2\to P_1: ((p,x,q),(q,y,r)) \mapsto (p,xy,r)$.
This  is called the \textbf{pullback groupoid} of
 $\Gamma_1\rightrightarrows \Gamma_0$ through $J$ and is denoted by
 $\Gamma_1[P_0] \rightrightarrows P_0$.
Recall that a morphism of Lie groupoids  $J$
from  $P_1\toto P_0$ to $ \Gamma_1\toto \Gamma_0$
is said to be a \textbf{Morita morphism} if $J: P_0\to \Gamma_0$ is a
surjective submersion and  $P_1\rightrightarrows P_0$ is
 isomorphic to the pullback groupoid of $\Gamma_1\rightrightarrows \Gamma_0$
 through  $J$.
Two Lie groupoids $\Gamma_1\rightrightarrows \Gamma_0$
 and $\Delta_1\rightrightarrows \Delta_0$ are said to be
 \textbf{Morita equivalent} if
 there exists a third Lie groupoid $P_1\rightrightarrows P_0$ together with
a Morita morphism from $P_{\simplicial}$
 to $\Gamma_{\simplicial}$ and a Morita morphism from $P_{\simplicial}$ to
 $\Delta_{\simplicial}$ \cite{Moer:regular}.
%$$ \xymatrix{
%\Gamma_1 \ar@<0.5ex>[d] \ar@<-0.5ex>[d] & P_1 \ar[l] \ar[r] \ar@<0.5ex>[d] \ar@<-0.5ex>[d] & \Delta_1 \ar@<0.5ex>[d] \ar@<-0.5ex>[d] \\
%\Gamma_0 & P_0 \ar[l] \ar[r] & \Delta_0 } $$
Equivalently, two Lie groupoids $\Gamma_1\rightrightarrows \Gamma_0$
and $\Delta_1\rightrightarrows \Delta_0$ are Morita equivalent if there exists
 a manifold $P_0$, two surjective submersions
 $P_0\xrightarrow{f} \Gamma_0$, $P_0\xrightarrow{g} \Delta_0$ and an isomorphism
of groupoids between $\Gamma_1[P_0]\toto P_0$ and
$\Gamma_1[P_0]\toto P_0$.

%$$ \xymatrix{
%\Gamma_1[P_0] \ar[r]^{\isomorphism} \ar@<0.5ex>[d] \ar@<-0.5ex>[d] & \Delta_1[P_0] \ar@<0.5ex>[d] \ar@<-0.5ex>[d] \\
%P_0 \ar@2{-}[r] & P_0 } .$$

%\begin{rmk}
%\begin{enumerate}
%\item If one drops the differential structures, Morita equivalent groupoids are nothing but weakly equivalent (small) categories.

Note that   there is a 1-1 correspondence between Morita equivalence classes
 of Lie groupoids and (equivalence classes of) differentiable stacks \cite{BX, BX1}.
%\end{enumerate}
%\end{rmk}

Now we are ready to introduce the definition of
Morita equivalence of groupoid  extensions.

\begin{defn}
A \textbf{Morita morphism} $f$ from a Lie groupoid extension
 $X'_1\xrightarrow{\phi} Y'_1\rightrightarrows M'$ to another extension
 $X_1\xrightarrow{\phi} Y_1\rightrightarrows M$ consists 
of Morita morphisms
$$ \xymatrix{
X'_1 \ar@<0.5ex>[d] \ar@<-0.5ex>[d] \ar[r]^{f} & X_1 \ar@<0.5ex>[d]
 \ar@<-0.5ex>[d] \\
M' \ar[r]_{f} & M} \quad \text{and} \quad 
\xymatrix{
Y'_1 \ar@<0.5ex>[d] \ar@<-0.5ex>[d] \ar[r]^{f} & Y_1 \ar@<0.5ex>[d] \ar@<-0.5ex>[d] \\
M' \ar[r]_{f} & M } $$
such that the diagram
\begin{equation}
\label{eq:morita}
 \xymatrix{
X'_1 \ar[d]^{\phi'} \ar[r]^{f} & X_1 \ar[d]^{\phi} \\
Y'_1 \ar@<0.5ex>[d] \ar@<-0.5ex>[d] \ar[r]^{f} & Y_1 
\ar@<0.5ex>[d] \ar@<-0.5ex>[d] \\
M' \ar[r]_{f} & M } 
\end{equation}
 commutes.
\end{defn}

It is simple to see that any Morita morphism of Lie groupoids
$$\xymatrix{
X'_1 \ar@<0.5ex>[d] \ar@<-0.5ex>[d] \ar[r]^{f} & X_1 \ar@<0.5ex>[d] \ar@<-0.5ex>[d] \\
M' \ar[r] & M }$$
such that
\begin{equation}
f(\ker\phi')=\ker\phi
\label{6}
\end{equation}
induces a Morita morphism between the groupoid extensions $
X'\simplicial\xrightarrow{\phi'}Y'\simplicial$ and $X\simplicial\xrightarrow{\phi}Y\simplicial$.
The converse is also true.

\begin{defn}
Two Lie groupoid extensions
$X_1\xrightarrow{\phi} Y_1\rightrightarrows M$ and $X'_1\xrightarrow{\phi'} Y'_1\rightrightarrows M'$
are said to be \textbf{Morita equivalent} if
 there exists a third extension
$X''_1\xrightarrow{\phi''} Y''_1\rightrightarrows M''$
 together with a  Morita morphism from $X''\simplicial\xrightarrow{\phi''}Y''\simplicial$ to $X\simplicial\xrightarrow{\phi} Y\simplicial$ and
a  Morita morphism from $X''\simplicial\xrightarrow{\phi''}Y''\simplicial$ 
to  $X'\simplicial\xrightarrow{\phi'}Y'\simplicial$.
\end{defn}

Equivalently, two Lie groupoid extensions $X_1\xrightarrow{\phi}
Y_1\rightrightarrows M$ and $X'_1\xrightarrow{\phi'}
Y'_1\rightrightarrows M'$  are Morita equivalent if there exists
 a manifold $P$, two surjective submersions
 $P\xrightarrow{\pi} M$, $P\xrightarrow{\pi'} M'$ and an isomorphism
of Lie  groupoid extensions between $X_1[P]\xrightarrow{\phi}
Y_1[P]\rightrightarrows P$ and $X'_1[P]\xrightarrow{\phi'}
Y'_1[P]\rightrightarrows P$.

One easily checks that this yields an equivalence
 relation on Lie  groupoid extensions.
Morita equivalence can also be defined  in terms of
 bitorsors \cite{hilsum-skandalis87}.

Recall that
a  Lie groupoid $\Gamma_1\rightrightarrows\Gamma_0$
 is said to \textbf{act on a manifold} $P$ (from the left)
 if there exists a map $J:P\to \Gamma_0$, called the momentum map, 
and an action map $\Gamma_1\times_{{\bt},\Gamma_0,J} P\to P:
 (x,p) \mapsto x\cdot p$ such that 
$x_1\cdot(x_2\cdot p)=(x_1 x_2)\cdot p$ and $J(p)\cdot p=p$,
 for all $(x_1,x_2,p)\in \Gamma_1\times_{{\bt},\Gamma_0,{\bs}}
\Gamma_1\times_{{\bt},\Gamma_0,J} P$.
Let $\Gamma_1\rightrightarrows\Gamma_0$ be a Lie groupoid and $M$ a manifold. A (left) \textbf{$\Gamma_{\simplicial}$-torsor over $M$} is a manifold $P$ together with a surjective submersion $\pi:P\to M$ called structure map and a left action of $\Gamma_{\simplicial}$ on $P$ such that the action is free and the quotient space $\Gamma_1\backslash P$ is diffeomorphic to $M$.  In other words, for any pair of elements $p_1,p_2$ in $P$ satisfying $\pi(p_1)=\pi(p_2)$, there exists a unique solution $x\in \Gamma_1$ to the equation $x\cdot p_1=p_2$.
A \textbf{$\Gamma_{\simplicial}$-$\Delta_{\simplicial}$-bitorsor} is a
 manifold $P$ together with two smooth maps $f:P\to \Gamma_0$ 
and $g:P\to \Delta_0$,  and commuting actions of $\Gamma_{\simplicial}$
 from the left and
of  $\Delta_{\simplicial}$ from the right on $P$, 
 with momentum maps $f$ and $g$ respectively,
 such that $P$ is at the same time a left-$\Gamma_{\simplicial}$-torsor over $\Delta_0$ (via $g$) and a right $\Delta_{\simplicial}$-torsor over $\Gamma_0$ (via $f$).

\begin{prop} \label{7}
Let $\Gamma_{\simplicial}$ and $\Delta_{\simplicial}$ be two Lie groupoids. 
There exists a $\Gamma_{\simplicial}$-$\Delta_{\simplicial}$-bitorsor if and only if
 $\Gamma_{\simplicial}$ and $\Delta_{\simplicial}$ are Morita equivalent.
\end{prop}

\begin{proof} 
This is a  standard result. We  will sketch a proof here, which
is needed later on.

\fbox{$\Rightarrow$}
Choose a bitorsor $P_0$. Let $P_1$ be the fibered product
$\Gamma_1\times_{{\bt}, \Gamma_0} P_0\times_{\Delta_0, {\bs}} \Delta_1$.
There is a natural groupoid structure on $P_1\rightrightarrows P_0$. The source and target maps are defined by
\begin{equation}
\label{eq:bist}
s(\gamma,p,\delta)=p \quad \text{and} \quad {\bt}(\gamma,p,\delta)=\gamma\cdot p\cdot \delta .
\end{equation}
Hence a pair $(\gamma_1,p_1,\delta_1),(\gamma_2,p_2,\delta_2)$ of elements of $P_1$ is composable if and only if $\gamma_1\cdot p_1\cdot \delta_1=p_2$.
The product is defined by
\begin{equation}
\label{eq:bim}
(\gamma_1,p_1,\delta_1)\cdot (\gamma_2,p_2,\delta_2)=(\gamma_2 \gamma_1,p_1,\delta_1 \delta_2) .
\end{equation}
One checks that the maps
\begin{gather*}
\Gamma_1 \gets P_1 \to \Delta_1 \\
\gamma\inv \mapsfrom (\gamma,p,\delta) \mapsto \delta
\end{gather*}
induce Morita morphisms 
$$ \xymatrix{
\Gamma_1 \ar@<0.5ex>[d] \ar@<-0.5ex>[d] & P_1 \ar[l] \ar[r] \ar@<0.5ex>[d] \ar@<-0.5ex>[d] & \Delta_1 \ar@<0.5ex>[d] \ar@<-0.5ex>[d] \\
\Gamma_0 & P_0 \ar[l] \ar[r] & \Delta_0 } .$$
Indeed, the maps $$P_1 \to \Delta_1[P_0]:(\gamma,p,\delta)\mapsto (p,\delta,\gamma p\delta)$$ and
$$P_1 \to \Gamma_1[P_0]:(\gamma,p,\delta)\mapsto (p,\gamma\inv,\gamma p\delta)$$ induce the required Lie groupoid isomorphisms between $P_1\rightrightarrows P_0$ and the pull-backs of $\Delta_{\simplicial}$ and $\Gamma_{\simplicial}$ over $P_0$.

\fbox{$\Leftarrow$}
This follows from the following two facts.
\begin{itemize}
\item If $\Gamma_{\simplicial}\to \Delta_{\simplicial}$ is a Morita morphism,
 then $P=\Gamma_0\times_{\Delta_0,{\bs}}\Delta_1$ is naturally
a $\Gamma_{\simplicial}$-$\Delta_{\simplicial}$-bitorsor.
\item If $P$ is a $\Gamma_{\simplicial}$-$\Delta_{\simplicial}$-bitorsor
 and $Q$ is a $\Delta_{\simplicial}$-$E_{\simplicial}$-bitorsor, then 
$(P\times_{\Delta_0}Q)/\Delta_1$ is a $\Gamma_{\simplicial}$-$E_{\simplicial}$-bitorsor.
\end{itemize}
\end{proof}

\begin{prop}
\label{pro:morita}
Let $X_1\xrightarrow{\phi} Y_1\rightrightarrows M$ and $X'_1\xrightarrow{\phi'} Y'_1\rightrightarrows M'$ be two Lie groupoid extensions with kernels $\kernel $ and $\kernel '$ respectively. Then $X\simplicial\xrightarrow{\phi}Y\simplicial$
 and $X'\simplicial\xrightarrow{\phi'}Y'\simplicial $ 
are \textbf{Morita equivalent} if and only if there exists an $X_{\simplicial}$-$X'_{\simplicial}$-bitorsor $B$ such that the orbits of the 
induced actions of $\kernel $ and $\kernel '$ on $B$ coincide and the corresponding orbit space is a manifold.
\end{prop}

A bitorsor as in the above proposition
 will be called an
 \textbf{extension $\phi_{\simplicial}$-$\phi'_{\simplicial}$-bitorsor}.
We refer the reader to \cite{BX1, TXL} for the discussion on Morita
equivalence of $S^1$-central extensions.

The proof of this proposition will be postponed to the next section.

\begin{rmk}
There is a 1-1 correspondence between Morita equivalence classes of Lie groupoid extensions and (equivalence classes of)
differentiable gerbes over  stacks.
\end{rmk}

\subsection{Generalized morphisms of Lie groupoid extensions}
\label{sec:gen}

The material here is parallel to Section 2.1 in 
 \cite{TXL} and follows the approach 
in \cite{hilsum-skandalis87}, to which we refer the reader for details.

\begin{defn}
A \textbf{strict morphism of groupoid extensions} is a
commutative diagram  as below,
 where the  horizontal arrows are groupoid homomorphisms:
\begin{equation}
\label{eq:strict}
 \xymatrix{
X'_1 \ar[d]^{\phi'} \ar[r]^{f} & X_1 \ar[d]^{\phi} \\
Y'_1 \ar@<0.5ex>[d] \ar@<-0.5ex>[d] \ar[r]^{f} & Y_1
\ar@<0.5ex>[d] \ar@<-0.5ex>[d] \\
M' \ar[r]_{f} & M }
\end{equation}
and, for any $m' \in M' $, the restriction
of  $f:X_1'\to X_1$ to a map $\ker \phi'|_{m'}
\to \ker \phi|_{f(m')}$ is an isomorphism.
\end{defn}

In particular, Morita morphisms of groupoid extensions
are strict homomorphisms of groupoid extensions.

\begin{defn}
Let $X_1\xrightarrow{\phi} Y_1\rightrightarrows M$ and
$X'_1\xrightarrow{\phi'} Y'_1\rightrightarrows M'$ be two Lie
groupoid extensions with kernels $\kernel $ and $\kernel '$
respectively. A \textbf{generalized morphism  of groupoid
extensions} from  $X_1'\xrightarrow{\phi'} Y_1'\rightrightarrows M'$
to $X_1\xrightarrow{\phi} Y_1\rightrightarrows M$ is a
$X_{\simplicial}'$-$X_{\simplicial}$ bimodule $M' \xleftarrow{f} B
\xrightarrow{g} M$ (i.e. $B$ is a $X_{\simplicial}'$-left space and
$X_{\simplicial}$-right space and the  $X_{\simplicial}'$ and $X_{\simplicial}$
actions commute) satisfying
\begin{enumerate}
\item   $B$ is a   $X_{\simplicial}$-torsor, and
\item   the induced actions of $\kernel $ and $\kernel '$
 on $B$ are free, and their orbits
coincide so that  the corresponding orbit space is a manifold.
\end{enumerate}
\end{defn}

It is easy to see that, in this case, 
 $\kernel \backslash B=B / \kernel '$
is a generalized morphism from 
$Y'_1\rightrightarrows M'$ to $Y_1\rightrightarrows M$.

\begin{lem}
\label{lem:gen}
Strict homomorphisms of groupoid extensions are 
generalized homomorphisms.
\end{lem}
\begin{proof}
Let  $f$ be   a strict homomorphism of groupoid extensions
from   $X_1'\xrightarrow{\phi'} Y_1'\rightrightarrows M'$ to
$X_1\xrightarrow{\phi} Y_1\rightrightarrows M$.
Set $B_f=  M' \times_{f, M, {\bs}} X_1$, 
 where $B_f \to M'$ is $(m' , x)\mapsto m'$, and $B_f \to M$
is $(m',x)\mapsto{\bt}(x)$.
The left $X'_{\simplicial}$-action is given by
$x'\cdot ({\bt}(x'),x)=( {\bs}(x'),f(x')x)$
 and the right $X_{\simplicial}$-action is
given by
$(m',x_1)\cdot x_2=(m',x_1 x_2)$.
It is simple to see that  the induced actions
 of $\ker\phi'$ and $\ker\phi$  on $B$ are free
and their orbits coincide.
The corresponding orbit space is the manifold
$M'\times_{f, M,{\bs}} Y_1$.
This concludes the proof.
\end{proof}

Just like strict homomorphisms, generalized homomorphisms of
 groupoid extensions can be composed.

\begin{prop}
\label{pro:com}
 Let $M'' \xleftarrow{f'} B' \xrightarrow{g'} M'$ be a
generalized morphism of groupoid extensions from
$X_1''\xrightarrow{\phi''} Y_1''\rightrightarrows M''$ to
$X_1'\xrightarrow{\phi'} Y_1'\rightrightarrows M'$,  and $M'
\xleftarrow{f} B \xrightarrow{g} M$   a generalized morphism of
groupoid extensions from $X_1'\xrightarrow{\phi}
Y_1'\rightrightarrows M'$ to $X_1\xrightarrow{\phi}
Y_1\rightrightarrows M$, their composition is the generalized
morphism of groupoid extensions from $X_1''\xrightarrow{\phi''}
Y_1''\rightrightarrows M''$ to $X_1\xrightarrow{\phi}
Y_1\rightrightarrows M$ given by the
$X_{\simplicial}''$-$X_{\simplicial}$-bimodule
 $$B'' =  {(B' \times_{M'} B)}/X'_1 ,$$
where $X'_1\toto M'$ acts on  ${B' \times_{M'} B}$ diagonally:
$(b',b)\cdot x'= (b' \cdot  x', {x'}\inv  b)$, for all
compatible $x'\in X'_1$ and $(b',b)\in {B' \times_{M'} B}$.
In particular, if both $M' \xleftarrow{f} B \xrightarrow{g} M$ and
 $M'' \xleftarrow{f'} B' \xrightarrow{g'} M'$ are bitorsors,
the resulting composition is  a $\phi\simplicial-\phi''\simplicial$-bitorsor.

Moreover, the composition of generalized morphisms
is associative.
\end{prop}

As a consequence, we obtain a category, where the
objects are groupoid extensions and the morphisms are
generalized homomorphisms of groupoid extensions.
Invertible morphisms exactly correspond to Morita equivalence of groupoid
extensions. As usual, we can decompose a generalized homomorphism
of groupoid extensions as the composition of
a Morita  equivalence with a  strict  homomorphism.

\begin{prop}\label{prop:decom}
Any generalized homomorphism of groupoid extensions
$M'\xleftarrow{f}B\xrightarrow{g}M$
from $X_1'\xrightarrow{\phi'}Y_1'\toto M'$ to
$X_1\xrightarrow{\phi}Y_1\toto M$
can be decomposed as the composition of the canonical
 Morita equivalence between $X'\simplicial \to Y'\simplicial$ and $X'\simplicial[B] 
\to Y'\simplicial[B]$,
with a strict  homomorphism  of groupoid extensions
from $X'\simplicial[B] \to Y'\simplicial[B]$ to $X\simplicial \to Y\simplicial$.
\end{prop}

\begin{proof}
Denote by $X'_1[B]\xrightarrow{} Y'_1[B]\rightrightarrows B$
the pull pack extension of $X_1'\xrightarrow{} Y_1'\rightrightarrows M'$
 via the surjective submersion
$B\xrightarrow{f}M'$. Then the projection
from    $X'\simplicial[B] \to Y'\simplicial[B]$ 
to $X'\simplicial \to Y'\simplicial$ is a Morita morphism,
 which induces a Morita
equivalence between these two groupoid extensions. 

As in the proof of Proposition \ref{7}, consider the
fiber product $X'_1\times_{{\bt}, M'}B\times_{M, {\bs}}X_1$.
Thus $X'_1\times_{{\bt}, M'}B\times_{M, {\bs}}X_1\toto B$ is
a Lie  groupoid, where the source, target, and multiplication
are given by Eqs.~\eqref{eq:bist} and \eqref{eq:bim}.
Introduce an equivalence relation in 
$X'_1\times_{{\bt}, M'}B\times_{M, {\bs}}X_1\toto B$
by 
$ (x' ,b, x)\sim (x' k' ,b, k\inv  x)$, iff $k'b=bk$,
where $k'\in \kernel'_{s(b)}$ and $k\in \kernel_{{\bt}(b)}$.
It follows from a direct verification, using Eqs.~\eqref{eq:bist} and \eqref{eq:bim},
 that the groupoid structure on $X'_1\times_{{\bt}, M'}B\times_{M, {\bs}}X_1\toto B$
descends to the quotient and
$X'_1\times_{{\bt}, M'}B\times_{M, {\bs}}X_1\to (X'_1\times_{{\bt}, M'}B\times_{M, {\bs}}X_1)
/\sim\toto B$ is
a groupoid extension. Moreover, one has the following commutative
diagram
\begin{multline}
\label{eq:XBX'}
\xymatrix{
X'_1[B] \ar[d] & X'_1\times_{{\bt}, M'}B\times_{M, {\bs}}X_1 \ar[l] \ar[d] \ar[r] & X_1 \ar[d] \\
Y'_1[B] \ar@<0.5ex>[d] \ar@<-0.5ex>[d] &( X'_1\times_{{\bt}, M'}B\times_{M, {\bs}}X_1)/\sim 
 \ar[l] \ar@<0.5ex>[d] \ar@<-0.5ex>[d] \ar[r]
 & Y_1 \ar@<0.5ex>[d] \ar@<-0.5ex>[d] \\
B & B \ar[l] \ar[r] & M } \\
\xymatrix{
(b, (x')\inv, x'bx) \ar[d] & (x',b,x) \ar[l] \ar[d] \ar[r] & x \ar[d] \\
(b, \phi'((x')\inv), x'bx) \ar@<0.5ex>[d] \ar@<-0.5ex>[d] & [x',b,x] \ar[l] \ar@<0.5ex>[d] \ar@
<-0.5ex>[d] \ar[r] & \phi(x) \ar@<0.5ex>[d] \ar@<-0.5ex>[d] \\
b & b \ar[l] \ar[r] & g(b) }
\end{multline}
where the left arrow is an isomorphism of groupoid
extensions, while the right arrow is a  strict homomorphism
 of groupoid extensions.

This concludes the proof of the proposition.
\end{proof}

Now we are ready to prove Proposition \ref{pro:morita}.

\begin{proof}
\fbox{$\Rightarrow$}
As in the proof of Lemma \ref{lem:gen}, 
$M'\xleftarrow{f}B_f\xrightarrow{g}M$ is indeed a $\phi'\simplicial-\phi\simplicial$-bitorsor.
Thus the conclusion follows from Proposition \ref{pro:com}.

\fbox{$\Leftarrow$}
Note that since $M'\xleftarrow{f}B\xrightarrow{g}M$ is a bitorsor,
the right arrow in diagram \eqref{eq:XBX'} is
indeed Morita morphism.
\end{proof}

\subsection{The outer action}

Let $G$ and $H$ be any Lie groups. We define $\Iso(G,H)$ as the set of group isomorphisms $G\xleftarrow{f} H$.
Let $\AD_h\in \Aut (H)$ denote the  conjugation by $h\in H$. The map
$$\Iso(G,H)\times H\to \Iso(G,H):(f,h)\mapsto f\rond\AD_h$$ defines an action of $H$ on $\Iso(G,H)$. The quotient of this action is 
 $\Out(G,H)$ consisting of outer isomorphisms from $H$ to $G$.
Since, for all $f\in\Aut(G)$ and $g\in G$, $f\rond\AD_g=\AD_{f(g)}\rond f$, 
the set $\Out(G,G)$ is a Lie group which will be denoted by $\Out(G)$.

 For any group  bundle  $\mathcal{G}\to M$,
let $\Iso(\mathcal{G},\mathcal{G})=\coprod_{m,n\in M}\Iso(\mathcal{G}_m,
\mathcal{G}_n)$,
 where $\mathcal{G}_m$ stands for the group fiber of $\mathcal{G}$
at  the point $m$. Similarly, 
  $\Out(\mathcal{G},\mathcal{G})=\coprod_{m,n\in M}\Out(\mathcal{G}_m,\mathcal{G}_n)$.

The following proposition is obvious.

\begin{prop}
\begin{enumerate}
\item The maps
$s:\Iso(\mathcal{G},\mathcal{G})\to 
M:(\mathcal{G}_m\xleftarrow{f}\mathcal{G}_n)
\mapsto m$, $t:\Iso(\mathcal{G},\mathcal{G})\to M: \ 
(\mathcal{G}_m\xleftarrow{f}\mathcal{G}_n)\mapsto n$
 and $m:\Iso(\mathcal{G},\mathcal{G})\times_{{\bt},M,{\bs}}
\Iso(\mathcal{G},\mathcal{G})\to\Iso(\mathcal{G},\mathcal{G}): 
(\mathcal{G}_m\xleftarrow{f}\mathcal{G}_n,\mathcal{G}_n\xleftarrow{g}\mathcal{G}_p)\mapsto
(\mathcal{G}|_m\xleftarrow{f\rond g}\mathcal{G}|_p)$
endow $\Iso(\mathcal{G},\mathcal{G})\overset{s}{\underset{t}{\rightrightarrows}}M$ with a groupoid structure.
\item These maps descend to the quotient $\Out(\mathcal{G},\mathcal{G})$, yielding a groupoid $\Out(\mathcal{G},\mathcal{G})\rightrightarrows M$.
\item The quotient map $\Iso(\mathcal{G},\mathcal{G}) \to \Out(\mathcal{G},\mathcal{G})$ is a groupoid morphism.
\end{enumerate}
\end{prop}

%\begin{proof}
%\begin{enumerate}
%\item trivial
%\item \label{point2} The multiplication descends to the quotient
 %since $f\rond\AD_h=\AD_{f(h)}\rond f$ for any group isomorphism $f$.
%\item This is an immediate consequence of \ref{point2}.
%\end{enumerate}
%\end{proof}

\begin{defn}
Let $\Gamma_1\rightrightarrows M$ be a Lie
 groupoid and  $\mathcal{G}\to M$  a bundle of  Lie groups over the same base.
\begin{enumerate}
\item An action by isomorphisms 
of $\Gamma_1\rightrightarrows M$ on $\mathcal{G}\to M$ is a
Lie  groupoid morphism
$$ \xymatrix{
\Gamma_1 \ar[r] \ar@<0.5ex>[d] \ar@<-0.5ex>[d] & \Iso(\mathcal{G},\mathcal{G}) \ar@<0.5ex>[d] \ar@<-0.5ex>[d] \\
M \ar@2{-}[r] & M }.$$
\item An  action by outer isomorphisms of $\Gamma_1\rightrightarrows M$ on $\mathcal{G}\to M$ is a Lie groupoid morphism
$$ \xymatrix{
\Gamma_1 \ar[r] \ar@<0.5ex>[d] \ar@<-0.5ex>[d] & \Out(\mathcal{G},\mathcal{G}) \ar@<0.5ex>[d] \ar@<-0.5ex>[d] \\
M \ar@2{-}[r] & M }.$$
\end{enumerate}
\end{defn}

\begin{prop}
\label{pro:outer}
Let $X_1\xrightarrow{\phi}Y_1\toto M$ be a Lie groupoid extension
with kernel $\kernel$. Then,
\begin{enumerate}
\item the groupoid $X_1\rightrightarrows M$ acts on $\kernel \to M$ by
conjugation. I.e.
\begin{equation}
\label{eq:conj}
 \xymatrix{
X_1 \ar[r]^{\AD} \ar@<0.5ex>[d] \ar@<-0.5ex>[d] & \Iso(\kernel ,\kernel ) \ar@<0.5ex>[d] \ar@<-0.5ex>[d] \\
M \ar@2{-}[r] & M }
\end{equation}
given by $\AD: x\to \AD_x, \forall x\in X_1$ is a groupoid morphism;
\item the composition of the groupoid morphism (\ref{eq:conj})
with the quotient morphism
$$ \xymatrix{\Iso(\kernel ,\kernel ) \ar[r]^{q} \ar@<0.5ex>[d] \ar@<-0.5ex>[d] & \Out(\kernel ,\kernel ) \ar@<0.5ex>[d] \ar@<-0.5ex>[d] \\
M \ar@2{-}[r] & M }$$
factorizes through $\phi$ :
$$ \xymatrix{ X_1 \ar[rd]^{q\rond\AD } \ar[d] \ar[d]_{\phi} \\
Y_1 \ar[r]_{\ADb} & \Out(\kernel ,\kernel ) ,}$$
yielding an action by outer isomorphisms $\ADb$
 of $Y_1 \rightrightarrows M$ on $\kernel \to M$, which is called the \textbf{outer action}.
\end{enumerate}
\end{prop}

\begin{rmk}
Let $Y_1\rightrightarrows M$ be a Lie groupoid
 and $X_1\xrightarrow{\phi}Y_1$ a surjective submersion. 
Assume $m_1$ and $m_2$ are two groupoid multiplications on $X_1$ making $X_1\xrightarrow{\phi}Y_1\rightrightarrows M$ into 
an extension of the groupoid $Y_1\rightrightarrows M$ in two different ways. 
Then the outer actions $\ADb_1$ and $\ADb_2$ induced by the two
 groupoid structures on $X_1$ are equal.
\end{rmk}

\section{Differentiable $G$-gerbes as groupoid $G$-extensions}

\subsection{Groupoid $G$-extensions}

Let $G$ be a fixed Lie group.

\begin{defn}
A Lie groupoid extension $X_1\xrightarrow{\phi}Y_1\rightrightarrows M$ is
called a \textbf{$G$-extension} if its kernel
$\kernel\to M$ is a locally trivial bundle of groups with
fibers isomorphic to $G$.
\end{defn}

It is simple to see that this notion is  invariant under Morita equivalence.
Namely, any groupoid extension Morita equivalent to a groupoid $G$-extension
must be  a  groupoid $G$-extension  itself.

The following proposition shows that any $G$-extension admits a Morita equivalent
$G$-extension whose kernel is a trivial bundle of groups.

\begin{prop}
Let $X_1\xrightarrow{\phi}Y_1\rightrightarrows M$ be  a groupoid
extension. Then $\phi$ is a
 $G$-extension if and only if there is an open covering
$\gendex{U_i}{i\in I}$ of $M$ such that the kernel of the pullback extension $X_1
[\coprod U_i]\xrightarrow{\phi'}Y_1[\coprod U_i]\rightrightarrows \coprod U_i$ via the
covering map $\coprod_{i\in I} U_i \to M$ is isomorphic to the trivial
bundle of groups $\coprod_{i\in I} U_i\times G$.
\end{prop}

\begin{proof}
Let $\kernel$ denote the kernel of $\phi$. Since $\phi$ is a $G$-extension, there exists an open cover $\gendex{U_i}{i\in I}$ of $M$ such that the pullback $\kernel'$ of $\kernel\to M$ to $\coprod_{i\in I}{U_i}$ is isomorphic to the trivial bundle $\coprod_{i\in I}{U_i}\times G \to \coprod_{i\in I}{U_i}$. To conclude, it suffices to notice that $\kernel'$ is the kernel of the pullback extension $X_1[\coprod{U_i}]\xrightarrow{\phi'}Y_1[\coprod{U_i}]\toto\coprod{U_i}$.
\end{proof}

Finally we say that an extension is a {\bf trivial $G$-extension}
if it is isomorphic to the extension $Y_1\times G\to Y_1\toto M$,
where $Y_1\times G$ is equipped with the product
groupoid structure.

\subsection{Band of a groupoid $G$-extension}

Now we introduce an important notion for
a groupoid $G$-extension, namely, the band. It  corresponds
to the band of  a $G$-gerbe in terms of stack language.
First of all, let us recall the definition of $G$-torsors
(also known as $G$-principal bundles) over a
groupoid \cite{LTX}.

\begin{defn}
Let $\Gamma_1\rightrightarrows\Gamma_0$ be a Lie groupoid and $G$ a Lie group.
 A \textbf{$G$-torsor}, or \textbf{$G$-principal bundle 
over $\Gamma_{\simplicial}$} consists of a right principal 
$G$-bundle $P\xrightarrow{J} \Gamma_0$ endowed with a left 
action of $\Gamma_1\rightrightarrows\Gamma_0$ on $P$ with momentum
 map $J$ such that the actions of $G$ and $\Gamma_{\simplicial}$ on $P$ commute.
\end{defn}

Let  $X_1\xrightarrow{\phi}Y_1\rightrightarrows M$ be a Lie 
groupoid {$G$-extension} with kernel $\kernel$. Let $\Iso(\kernel,G)$ (resp.
$\Out(\kernel,G)$) be the group bundle $\coprod_{m\in M}\Iso(\kernel_m,G)$
 (resp. $\coprod_{m\in M}\Out(\kernel_m,G)$) and $J$ (resp. $\bar{J}$) the canonical projection of 
$\Iso(\kernel,G)$ (resp. $\Out(\kernel,G)$) onto $M$. 
First note that $\Iso(\kernel,G)\xrightarrow{J}M$ (resp.
$\Out(\kernel,G)\xrightarrow{\bar{J}}M$) is a right principal
$\Aut(G)$-bundle (resp. $\Out(G)$-bundle).

On the other hand, $\Iso(\kernel,G) \xrightarrow{J} M$ (resp. $\Out(\kernel,G)
\xrightarrow{\bar{J}} M$) admits a natural left action of the groupoid
$\Iso(\kernel,\kernel)\toto M$ (resp. $\Out(\kernel,\kernel)\toto M$). Moreover, these two actions commute. Hence $\Iso(\kernel,G)\xrightarrow{J}M$
(resp. $\Out(\kernel,G)\xrightarrow{\bar{J}}M$) is an
$\Aut(G)$-torsor (resp. $\Out(G)$-torsor) over $\Iso(\kernel,\kernel)\toto M$ (resp. $\Out(\kernel,\kernel)\toto M$). According to
 Proposition \ref{pro:outer}, there is a groupoid
morphism $\AD: X\simplicial\to
\Iso(\kernel, \kernel)\simplicial$ (resp. $\ADb: Y\simplicial\to
\Out(\kernel, \kernel)\simplicial$). 
Therefore, one may pull-back the
$\Out(G)$-torsor $\Out(\kernel,G) \xrightarrow{\bar{J}} M$
over $\Out(\kernel,\kernel)\toto M$
to an $\Out(G)$-torsor over $Y\simplicial$ via
$\ADb: Y\simplicial\to \Out(\kernel, \kernel)\simplicial$. More precisely, one
defines the left $Y\simplicial$-action on $\Out(\kernel ,G)
\xrightarrow{J} M$ by \begin{equation} \label{eq:yaction} y\cdot
f=\ADb (y) \rond f, \quad \forall f\in
\Out(\kernel ,G)|_{{\bt}(y)} .\end{equation}

\begin{defn} \label{def:band}
Let  $X_1\xrightarrow{\phi}Y_1\rightrightarrows M$ be a Lie  groupoid
{$G$-extension} with kernel $\kernel$.
Then $\Out(\kernel,G)\xrightarrow{\bar{J}} M$, considered as
an $\Out(G)$-torsor over $Y\simplicial$, is called
the \textbf{band} of the $G$-extension $\phi$.
\end{defn}

It is well-known  that for a given Lie group $G$ there is
a bijection between 
 $G$-torsors over a pair of Morita equivalent groupoids.
They are called  Morita equivalent $G$-torsors. 
The following proposition shows that bands are preserved under Morita morphisms, thus also by Morita equivalences in the above sense.

\begin{prop} \label{pro:3.5}
Let $f$ be  a  Morita morphism of Lie groupoid extensions
from  $X'_1\xrightarrow{\phi'}Y'_1\rightrightarrows M'$ 
to $X_1\xrightarrow{\phi}Y_1\rightrightarrows M$. 
%$$ \xymatrix{
%X'_1 \ar[r]^f \ar[d]_{\phi'} & X_1 \ar[d]^{\phi} \\
%Y'_1 \ar[r]^f \ar@<0.5ex>[d] \ar@<-0.5ex>[d] & Y_1 \ar@<0.5ex>[d] \ar@<-0.5ex>[d] \\
%M' \ar[r]^f & M }$$
%be a Morita morphism of groupoid extensions.
Let $\kernel '$ and $\kernel $ denote the kernels of $\phi'$ and $\phi$ respectively.
If  $\phi$ and $\phi'$ are groupoid   $G$-extensions, the band of $\phi'$:
$$ \xymatrix{ Y_1' \ar@<0.5ex>[dr] \ar@<-0.5ex>[dr] & \Out(\kernel',G) \ar[d] \\  & M' } $$
is isomorphic to the pullback of the band of $\phi$:
$$ \xymatrix{ Y_1 \ar@<0.5ex>[dr] \ar@<-0.5ex>[dr]  & \Out(\kernel,G) \ar[d] \\ & M } $$
via $M'\xrightarrow{f}M$.
\end{prop}

\begin{proof}
First, we observe that $\kernel'$ is isomorphic to the pullback bundle $M'\times_{M}\kernel$ and $f$ restricts to a bundle map $\kernel'\xrightarrow{f}\kernel$, which is a fiberwise  group isomorphism.
Therefore, the torsor 
$$ \xymatrix{ \Out(\kernel',\kernel') \ar@<0.5ex>[dr] \ar@<-0.5ex>[dr] & \Out(\kernel',G) \ar[d] \\ & M' } $$
is isomorphic to the pullback of the torsor 
$$ \xymatrix{ \Out(\kernel ,\kernel) \ar@<0.5ex>[dr] \ar@<-0.5ex>[dr] & \Out(\kernel ,G) \ar[d] \\ & M } $$
via $M'\xrightarrow{f}M$. Let
$$ \xymatrix{
Y'_1 \ar[r]^{\ADb'} \ar@<0.5ex>[d] \ar@<-0.5ex>[d] & \Out(\kernel ',\kernel ') \ar@<0.5ex>[d] \ar@<-0.5ex>[d] \\ M' \ar@2{-}[r] & M' } \qquad \text{and} \qquad \xymatrix{
Y_1 \ar[r]^{\ADb} \ar@<0.5ex>[d] \ar@<-0.5ex>[d] & \Out(\kernel ,\kernel ) \ar@<0.5ex>[d] \ar@<-0.5ex>[d] \\ M \ar@2{-}[r] & M }$$
be the outer actions  of $\phi'$ and $\phi$ respectively.
The morphism $f$ induces a Morita morphism
$$ \xymatrix{ \Out(\kernel ',\kernel ') \ar[r]^{\mathbf{f}} \ar@<0.5ex>[d] \ar@<-0.5ex>[d] & \Out(\kernel ,\kernel ) \ar@<0.5ex>[d] \ar@<-0.5ex>[d] \\ M' \ar[r]_{f} & M } $$
such that the diagram
$$ \xymatrix{ Y'_1 \ar[r]^{\ADb'} \ar[d]_{f} & \Out(\kernel',\kernel') \ar[d]^{\mathbf{f}} \\ Y_1 \ar[r]_{\ADb} & \Out(\kernel ,\kernel ) } $$
 commutes.
This completes the proof
\end{proof}

\begin{rmk}
In terms of stack and gerbe language, the band of a $G$-gerbe
 over a stack is an  $\Out(G)$-torsor over this stack \cite{Giraud}.
\end{rmk}

\subsection{$G$-bound gerbes and groupoid central extensions}

This  subsection is devoted to the study of an important class of $G$-gerbes,
 namely,  $G$-bound gerbes. These correspond to groupoid central extensions.

Recall that when $G$ is abelian,
 a groupoid $G$-central extension is a groupoid
 extension $X_1\xrightarrow{\phi}Y_1\toto M$ such that $M\times G\isomorphism \ker \phi , \ (m, g)\mapsto
 g_m\in \ker \phi|_m$ and
$$  x \cdot g_{{\bt}(x)}=g_{ {\bs}(x)} \cdot x.$$

\begin{example}
Consider the particular case of $S^1$-extensions.
%Since $S^1$ is abelian, $Z(S^1)=S^1$ and $\Out(S^1)=\Aut(S^1)$.
In this case, $X_1\xrightarrow{\phi}Y_1\rightrightarrows M$ is 
an $S^1$-central extension if and only if
% with kernel $\kernel$. Here, the outer action of an element $y\in Y^1$ is simply an isomorphism between the fibers $\kernel_{{\bt}(y)}$ and $\kernel_{ {\bs}(y)}$. And the extension is central if, and only if, 
there exists a trivialization of the kernel $$M\times S^1\to \ker \phi 
:(m,g )\mapsto g_m$$ such that 
\begin{equation} x\cdot g_{{\bt}(x)}=g_{ {\bs}(x)}\cdot x, 
\quad \forall x\in X_1,\;\forall g\in G . \label{gloup} \end{equation}
In other words, $X_1\xrightarrow{\phi}Y_1 $ is 
an $S^1$-principal bundle such that
$$ (g_1 x_1) \cdot (g_2 x_2)=(g_1g_2)  \cdot (x_1x_2), \ \ \forall g_1, g_2\in S^1, 
(x_1, x_2)\in X_2.$$
See \cite{WeinsteinX:prequantization} for details.
Gerbes represented by $S^1$-central extensions are called $S^1$-bound
gerbes\footnote{In \cite{BX, BX1} they are called $S^1$-gerbes for
simplicity. When $Y\simplicial$ is Morita equivalent to a manifold,
they are called bundle gerbes by Hitchin \cite{Hitchin} and Murray \cite{Murray}.}. 
\end{example}

However, when $G$ is non-abelian, the situation is more subtle.
First of all, we need to introduce the following

\begin{defn}
Let $X_1\xrightarrow{\phi}Y_1\rightrightarrows M$ be a groupoid $G$-extension with kernel $\kernel\to M$. Its band $\Out(\kernel,G)\xrightarrow{J}M$ is said to be trivial (as a $Y\simplicial$-torsor) if there exists a section $\etab$ of $\Out(\kernel,G)\to M$ invariant under the $Y\simplicial$-action: 
\begin{equation}
\label{eq:trivial}
y\cdot\etab({\bt}(y))
=\etab( {\bs}(y)), \quad \forall y\in Y_1 .
\end{equation}
Such a section is called a \textbf{trivialization of the band}.
\end{defn}

\begin{prop} \label{ads}
 If $X_1\xrightarrow{\phi}Y_1\rightrightarrows M$ and $X_1'\xrightarrow{\phi'}Y_1'\rightrightarrows M'$ are Morita equivalent groupoid extensions, then the band of the first one is trivial if, and only if, the band of 
the second one is trivial.
\end{prop}

\begin{proof}
 As in Proposition \ref{pro:3.5}, we can restrict ourselves to the case
that there is  a Morita morphism  of groupoid extensions 
from $X_1'\xrightarrow{\phi'}Y_1'\rightrightarrows M'$ to $X_1\xrightarrow{\phi}Y_1\rightrightarrows M$.
 Since $\kernel'\isomorphism M'\times_M \kernel$, there is a natural 1-1 correspondence between the $Y\simplicial$-invariant sections of $\Out(\kernel,G)$ and the $Y'\simplicial$-invariant sections of $\Out(\kernel',G)$. 
Hence $X_1'\xrightarrow{\phi'}Y_1'\rightrightarrows M'$ has trivial band if, and only if, 
$X_1\xrightarrow{\phi}Y_1\rightrightarrows M$  has trivial band.
\end{proof}

Any section $\eta$ of $\Iso(\kernel,G)\to M$ is called a 
\textbf{trivialization of the kernel}. 
For it defines a map $M\times G\to\kernel:(m,g)\mapsto g_m:
=\eta_m (g)$,
 which is an isomorphism of bundles of groups over $M$.
i.e.  $\eta_m (gh)=\eta_m (g) \eta_m (h)$.

The following proposition indicates that when the band is trivial,
a trivialization of the kernel always exists when passing
to a  Morita equivalent extension. 

\begin{prop}
 Let $X_1\xrightarrow{\phi}Y_1\rightrightarrows M$ be a groupoid $G$-extension with kernel $\kernel\to M$ whose band is trivial. Then there exists a Morita equivalent groupoid $G$-extension
$X'_1\xrightarrow{\phi'}Y'_1\rightrightarrows M'$ with kernel $\kernel'\to M'$ together with a section $\eta'$ of $\Iso(\kernel',G)\to M'$ such that its induced section $\etab'$ of $\Out(\kernel',G)\to M'$  is invariant under the $Y'\simplicial$-action.
\end{prop}

\begin{proof}
Choose a trivialization $\etab$ of the band $\Out(\kernel,G)\to M$ of $\phi$. Take a good open covering $\gendex{U_i}{i\in I}$ of $M$ and consider the pullback extension $X_1[\coprod U_i]\xrightarrow{\phi'}Y_1[\coprod U_i]\toto\coprod U_i$ of $\phi$ by the projection $\coprod_{i\in I}U_i\to M$. 
Let us call $\etab'$, the $Y\simplicial[\coprod U_i]$-invariant section of $\Out(\kernel',G)\to\coprod U_i$ associated to $\etab$ as in Proposition \ref{ads} above. Since the $U_i$'s are contractible, $\etab'$ can 
be lifted to a section $\eta'$ of $\Iso(\kernel',G)\to\coprod U_i$.
\end{proof}

The following proposition can be verified directly.

\begin{prop}
Let $H$ be  a subgroup of $Z(G)$.
Assume that  $\tilde{X}_1\xrightarrow{\phi}Y_1\rightrightarrows M$
is  an $H$-central  extension.  Let 
$$X_1=\frac{\tilde{X}_1 \times G}{H}, $$
where $H$ acts on $\tilde{X}_1 \times G$ diagonally:
$ (\tilde{x}, g) \cdot h=(\tilde{x}h_{{\bt}(x)}, h\inv g)$, $\forall
h\in H$. Then ${X}_1\xrightarrow{\phi}Y_1\rightrightarrows M$
is a $G$-extension endowed with a trivialization of its kernel,
called the induced $G$-extension.
\end{prop}

We are now ready to state the main result of this subsection.

\begin{thm} \label{thm:3.9}
Let $X_1\xrightarrow{\phi}Y_1\rightrightarrows M$ be a groupoid $G$-extension with kernel $\kernel\to M$, and let 
$\eta$ be a trivialization of the kernel. The following assertions are equivalent.
\begin{enumerate}
\item The section $\etab$ of $\Out(\kernel,G)\to M$ induced by the trivialization of the kernel $\eta:M\to\Iso(\kernel, G)$ is a trivialization of the band.
\item There exists a groupoid morphism $$\xymatrix{ X_1 \ar@<0.5ex>[d] \ar@<-0.5ex>[d] \ar[r]^{r} & G/Z(G) \ar@<0.5ex>[d] \ar@<-0.5ex>[d] \\ M \ar[r] & {*} }$$ with kernel $Z(\kernel)$ such that
 $$x\cdot g_{{\bt}(x)}=(\ADu_{r(x)}g)_{ {\bs}(x)}\cdot x, \forall x\in X_1,\;\forall g\in G,$$
where $\ADu$ denotes the canonical isomorphism from $G/Z(G)$ to $\Inn(G)$:
$$\ADu_{[g']}(g)=g' g (g')\inv , \ \ \forall g, g'\in  G. $$

\item There exists a global section $\sigma$ of
the induced morphism 
 $X_1/Z(\kernel)\to Y_1$ such that, for any local lift $\sigt:Y_1\to X_1$ of $\sigma$, $$\sigt(y)\cdot g_{{\bt}(y)}=g_{ {\bs}(y)}\cdot\sigt(y), \quad \forall y\in Y_1,\;\forall g\in G.$$ 
Moreover,  $\sigma$ is  a groupoid morphism.
\item The extension   $X_1\xrightarrow{\phi}Y_1\rightrightarrows M$ 
is isomorphic to the induced $G$-extension
$\frac{\tilde{X}_1 \times G}{Z(G)}\to Y_1 \rightrightarrows M$ of a $Z(G)$-central 
extension $\tilde{X}_1\xrightarrow{\phi}Y_1\rightrightarrows M$.
\end{enumerate}
\end{thm}

When  any of the conditions  above 
is   satisfied, the groupoid $G$-extension is called \textbf{central}.
differentiable gerbes represented by $G$-central extensions are
called \textbf{$G$-bound gerbes}.
\begin{proof}
\fbox{1$\Rightarrow$2} It is well known that the group morphism $\AD:G\to\Inn(G):g\mapsto\AD_g$ 
factorizes through $G/Z(G)$:
$$\xymatrix{ G \ar[dr]^{\AD} \ar[d] & \\ G/Z(G) \ar[r]_{\ADu} & \Inn(G), }$$
and $\ADu:G/Z(G)\to\Inn(G)$ is an isomorphism.

For any $x\in X_1$, it is clear that
$\eta_{ {\bs}(x)}\inv \smalcirc \AD_x \smalcirc \eta_{{\bt}(x)}$ is
an automorphism of $G$. According to Eq.~\eqref{eq:trivial},
it must be an inner automorphism, and therefore corresponds
to an element in $G/Z(G)$, which is defined
to be $r(x)$. Here we used the fact that 
$\ADu: \Inn(G)\to G/Z(G)$ is an isomorphism.
Thus we have
\begin{equation} \AD_x\rond\eta_{{\bt}(x)}=\eta_{ {\bs}(x)}\rond\ADu_{r(x)} \label{stara} \end{equation}
or 

\begin{equation}
 x\cdot g_{{\bt}(x)}=(\ADu_{r(x)}g)_{ {\bs}(x)} \cdot x,\quad \forall x\in X_1,\;\forall g\in G.
\end{equation}
From its definition, it is simple to see that
 $$\xymatrix{X_1\ar@<0.5ex>[d] \ar@<-0.5ex>[d] \ar[r]^r & G/Z(G) \ar@<0.5ex>[d] \ar@<-0.5ex>[d] \\ M \ar[r] & {*} }$$ 
is indeed a Lie  groupoid morphism.

\fbox{2$\Rightarrow$3} Since the groupoid morphism
 $\AD$ has kernel $Z(\kernel)$, it factorizes through $X_1/Z(\kernel)$: 
$$\xymatrix{X_1 \ar[dr]^{\AD} \ar[d] & \\ X_1/Z(\kernel) \ar[r]_{\ADu} & \Iso(\kernel,\kernel). }$$
Also the trivialization of the kernel $\eta:
M\to \Iso (\kernel, G)$ 
 induces a unique section $\etau$ of 
$\Iso\big(\kernel/Z(\kernel),G/Z(G)\big)\to M$ making
the diagram
 $$\xymatrix{ \kernel_m \ar[d] & G \ar[l]_{\eta_m} \ar[d] \\ 
{\big(\kernel/Z(\kernel)\big)}_m & G/Z(G) \ar[l]_{\etau_m} }$$ commutes. 
It thus follows that
$\eta_m\rond\AD_g = \AD_{\eta_m(g)}\rond \eta_m 
= \ADu_{[\eta_m(g)]}\rond\eta_m$.
Hence
\begin{equation}
\eta_m\rond\ADu_{[g]} = \ADu_{\etau_m[g]}\rond \eta_m, 
\label{squareplus} \end{equation}
where $[q]$ denotes the class of an element $q\in X_1$ (resp. $G$) in $X_1/Z(\kernel)$ (resp. $G/Z(G)$).
Using Eq.~\eqref{squareplus}, the formula \eqref{stara} becomes
$\AD_x\rond\eta_{{\bt}(x)}=\eta_{ {\bs}(x)}\rond\ADu_{r(x)}$. Hence
$\ADu_{[x]}\rond\eta_{{\bt}(x)}=\ADu_{\etau_{ {\bs}(x)}(r(x))}
\rond\eta_{ {\bs}(x)}$. Therefore 

\begin{equation}
\ADu_{\etau_{ {\bs}(x)}((r(x))\inv\cdot[x])}\rond\eta_{{\bt}(x)}
=\eta_{ {\bs}(x)}.
\label{hart} 
\end{equation}

For any $y\in Y_1$, take $x\in X_1$ any element such that
$\phi (x)=y$. Set
$$\sigma(y)=\etau_{ {\bs}(x)}(r(x))\inv\cdot[x] .$$
It is easy to see that  $\sigma$ defines a section of
 $X_1/Z(\kernel)\to Y_1$. 
Then Eq.~\eqref{hart} becomes
 \begin{equation}
 \ADu_{\sigma(\phi(x))}\rond\eta_{{\bt}(x)}=\eta_{ {\bs}(x)} , \ \ \forall x\in X_1. 
\label{cross} 
\end{equation}
 Therefore, for any
 local lift $\sigt:Y_1\to X_1$ of $\sigma$, we have 
\begin{equation}
 \AD_{\sigt(\phi(x))}\rond\eta_{{\bt}(x)}=\eta_{ {\bs}(x)} , \ \ \forall x\in X_1,
\label{sun} \end{equation} or,
 equivalently,
 \begin{equation} 
\sigt(\phi(x))\cdot g_{{\bt}(x)}=g_{ {\bs}(x)}\cdot\sigt(\phi(x)), \ \ \forall x\in X_1,
g\in G .  \end{equation}
I.e.
\begin{equation} 
\sigt(y)\cdot g_{{\bt}(y)}=g_{ {\bs}(y)}\cdot\sigt(y), \ \ \forall y\in Y_1, g\in G . 
\label{varstar}
 \end{equation}
Moreover, Eq.~\eqref{cross} implies that $\forall (x, y)\in X_2$,
$\eta_{{\bt}(y)}\inv\rond\ADu_{\big(\sigma((\phi(xy))\inv)\sigma(\phi(x))\sigma(\phi(y))\big)}\rond\eta_{{\bt}(y)}=\idn$.
Hence,
 $$\sigma((\phi(xy))\inv)\sigma(\phi(x))\sigma(\phi(y))=1 \quad \text{in } X_1/Z(\kernel),$$ 
since $\ADu$ is injective. It thus follows that  $\sigma$ is a
 groupoid morphism.

\fbox{3$\Rightarrow$4} Take $\tilde{X}_1=\pi \inv\big(\sigma(Y_1)\big)$,
 where $\pi $ denotes the projection $X_1\to X_1/Z(\kernel)$.
Since $\sigma$ is a Lie  groupoid morphism, it is simple to see
that  $\tilde{X}_1\toto M$ is
a Lie subgroupoid of ${X}_1\toto M$ and $\tilde{X}_1\to Y_1\toto M$ is 
a groupoid  $Z(G)$-extension.
The   trivialization of the kernel $\eta_m :  G\to \kernel_m$, when
being restricted to the center $Z(G)$, induces a
trivialization of the kernel $\tilde{\eta}_m :  Z(G)\to Z(\kernel)_m$
of the $Z(G)$-extension $\tilde{X}_1\to Y_1\toto M$.
Moreover, Eq.~\eqref{sun} implies that
 \begin{equation}
 \AD_{\tilde{x}}\rond\eta_{{\bt}(\tilde{x})}=\eta_{{\bs}(\tilde{x})},\quad\forall\tilde{x}\in\tilde{X}_1 . \label{moon} \end{equation}
It thus follows that  $\tilde{X}_1\to Y_1\toto M$ is
a $Z(G)$-central extension.
Consider the map
 $$\tau: \tilde{X}_1\times G\to X_1:(\tilde{x},g)\mapsto \tilde{x}\cdot 
g_{{\bt}(\tilde{x})}.$$ 
From Eq.~\eqref{varstar}, it follows that 
$\tau$ is a groupoid morphism.
According to Eq.~\eqref{cross}, we have 
 $$[x]=\etau_{ {\bs}(x)}(r(x))\cdot \sigma(\phi(x))=\sigma(\phi(x))\cdot 
\etau_{{\bt}(x)}(r(x)),\quad\forall x\in X_1.$$
It thus follows that
$x=\sigma(\phi(x))\cdot 
\etau_{{\bt}(x)}(r(x)) \cdot k$ for some $k\in Z(\kernel )$.
Hence $\tau$ is surjective.  And its kernel
 $\genrel{(z_m,z\inv)}{z\in Z(G)}$ is isomorphic to $Z(G)$. 
Therefore   $X\simplicial\xrightarrow{\phi}Y\simplicial$
is isomorphic to the induced $G$-extension  $\frac{\tilde{X}_1 \times G}{Z(G)} \to Y_1 \rightrightarrows M $.

\fbox{4$\Rightarrow$1} 
For any $x=[(\tilde{x}, g)]\in \frac{\tilde{X}_1 \times G}{Z(G)}$,
it follows from a simple computation that
$$ \eta_{ {\bs}(x)}\inv \smalcirc \AD_x \smalcirc \eta_{{\bt}(x)}=\AD_g.$$
The conclusion thus follows.
\end{proof}

As an immediate consequence, Morita equivalent classes of
$G$-extensions with trivial band are in one-one correspondence
with Morita equivalent classes of
$Z(G)$-central extensions. The latter is classified by
$H^2(Y\simplicial, Z(G))$. Thus we have recovered the following
result of Giraud \cite{Giraud} in the context of differential
geometry. 

\begin{thm}
 Morita equivalent classes of
$G$-extensions with trivial band (i.e. $G$-bound gerbes) over
the groupoid $Y_1\toto M$ are in one-one correspondence
with $H^2(Y\simplicial , Z(G))$.
\end{thm}

\subsection{$G$-gerbes over manifolds}

In this subsection, we study $G$-gerbes over a manifold.
This corresponds to a $G$-extension over  a
groupoid which is Morita equivalent to a  manifold, i.e.
a groupoid of the form $M\times_N M\toto M$ for a surjective
submersion $M\to N$ (see also \cite{Moerdijk}).
When $\gendex{U_i}{i \in I}$ is an open covering of
$N$ and $M$ is the disjoint union $\coprod_{i} U_i$ with the covering
map $\coprod_{i} U_i\to N$, the resulting groupoid
$M\times_N M\toto M$, which is easily seen to be
isomorphic to  $\coprod_{ij} U_{ij}\toto \coprod_{i} U_i$,
is called the \v{C}ech groupoid
 associated to an open covering $\gendex{U_i}{i \in I}$ of the manifold
$N$.

Let  $\gendex{U_i}{i \in I}$ be a {\bf good open covering} of $N$,
namely all $U_i$ and their finite intersections
are contractible.
 We form its \v{C}ech groupoid $Y\simplicial: \coprod_{ij} U_{ij}
 \rightrightarrows \coprod_{i} U_i$,
and consider a
 $G$-extension $\kernel \simplicial\to X\simplicial \xrightarrow{\phi} Y\simplicial$.
A point $x$ in $N$ will be denoted by
 $x_i$ when considered as a point
 in $U_i$ and by $x_{ij}$ when considered as a point in $U_{ij}$, etc.
 The source, target and
 multiplication maps of $Y\simplicial$ are given, respectively,  by
$$ {\bs}(x_{ij})=x_i, \ \ \ {\bt}(x_{ij})=x_j, \ \ \ x_{ij}\cdot x_{jk}=x_{ik}. $$

Since  $U_i$ is contractible for all $i\in I$, we can identify $\kernel $ with the trivial bundle of groups $\kernel \isomorphism \coprod_{i} U_i \times G$.
 Since $U_{ij}$ is contractible for all $i,j \in I$,
there exists a global section $\rho$ (which does not need to be
a groupoid morphism) of
$X_1\xrightarrow{\phi}Y_1$ such that 
$$\rho(x_{ii}) = \epsilon(x_i),\;\forall i\in I, \ \ \text{ and } \ 
\rho(x_{ji})={\big(\rho(x_{ij})\big)}\inv,\;\forall i,j\in I $$
where $\epsilon: M \to X_1$ is the unit map.

We identify any element $\tilde{x}\in X_1$ with a pair
$x_{ij}=\phi(\tilde{x})\in Y_1$ and $g\in\kernel _{x_j}$ such that $\rho(x_{ij})\cdot g=\tilde{x}$. In other words, we identify $X_1$
 with $\coprod_{ij}U_{ij}\times G$ in such a way that the multiplication of elements in $X_1$ by elements of $\kernel \isomorphism \coprod_{i} U_i \times G$ from the right coincides with the multiplication from the right on the group $G$,
i.e.
\begin{equation}
\label{eq:gh}
(x_{ij},g)\cdot(x_{j},h)=(x_{ij},gh) , \ \ \forall g, h \in G.
\end{equation}

First, we define $C^\infty (G, G)$-valued functions $\lambda_{ij}$ on the $2$-intersections $U_{ij}$ by comparing the right and the left actions of the group bundle $\kernel  \to \coprod_{i} U_i$. Indeed for any $i,j \in I$, $x_{ij}\in U_{ij}$
 and $g\in G$, there is a unique element $\lambda_{ij}(x_{ij})
\big(g\big)\in G$ such that
\begin{equation} \label{eq:lambda} (x_{ij},1)\cdot(x_{j},g)=(x_{i},\lambda_{ij}(x_{ij})(g))\cdot(x_{ij},1) .\end{equation}
Note that $g \to \lambda_{ij}(x_{ij})(g)$ is a smooth diffeomorphism 
by construction.

Second, we define $G$-valued functions $g_{ijk}$ on the 3-intersections $U_{ijk}$ by
$$(x_{ij},1)\cdot(x_{jk},1)=(x_{ik},g_{ijk}(x_{ijk})) .$$
These functions measure the default of $\rho$ from being a groupoid homomorphism.
In the sequel, the reference to the point $x$ where $\lambda_{ij}$ and $g_{ijk}$ are evaluated will be omitted.

The data $(\lambda_{ij},g_{ijk})$ determine completely the multiplication on the groupoid $X_1\rightrightarrows\coprod_{i}U_i$. More precisely, one has
\begin{equation} (x_{ij},g)(x_{jk},h)=(x_{ik},g_{ijk}\lambda_{jk}\inv (g)h),\quad\forall g,h\in G. \label{eq:defprod} \end{equation}

The associativity of the groupoid multiplication defined on $X_1$ imposes the following relations on the functions $\lambda_{ij}$ and $g_{ijk}$:
\begin{gather}
\lambda_{ij}(gh)=\lambda_{ij}(g)\lambda_{ij}(h),\quad\forall g,h\in G; \label{eq:cocycle0} \\
\lambda_{ij}\rond\lambda_{jk}=\AD_{g_{ijk}}\rond\lambda_{ik}; \label{eq:cocycle1} \\
g_{ijl}g_{jkl}=g_{ikl}\lambda_{kl}\inv (g_{ijk}); \label{eq:cocycle2}
\end{gather}
which are immediate consequences of  the following identities,
reflecting the associativity of the groupoid product:
\begin{gather*}
\big((x_{ij}, 1)(x_j, g)\big)(x_j, h)=(x_{ij}, 1) \big((x_j, g)(x_j, h) \big)\\
\big((x_i,g)(x_{ij},1)\big)(x_{jk},1)=(x_i,g)\big((x_{ij},1)(x_{jk},1)\big) \\
\big((x_{ij},1)(x_{jk},1)\big)(x_{kl},1)=(x_{ij},1)\big((x_{jk},1)(x_{kl},1)\big) .
\end{gather*}

Eq.~\eqref{eq:cocycle0} means that $\lambda_{ij}$ is
an $\Aut(G)$-valued function on $U_{ij}$. Alternatively,
one may write Eq.~\eqref{eq:lambda} as
$$\AD_{(x_{ij}, 1)}(x_j, g)=(x_i,  \lambda_{ij}(x_{ij})(g)).$$
From the fact that $\AD_{(x_{ij}, 1)}$ is a group isomorphism
from $\kernel_{x_j}$ to $\kernel_{x_i}$, it follows immediately
that $\lambda_{ij}$ is an $\Aut(G)$-valued function.

Conversely, if we are given some $\lambda_{ij}: U_{ij}\to \Aut(G)$,
$g_{ijk}:U_{ijk} \to G$ satisfying Eqs.~\eqref{eq:cocycle1} and \eqref{eq:cocycle2},
then the product defined in Eq.~\eqref{eq:defprod} defines a groupoid structure
on $X_1\toto\coprod_{ij}U_{ij}$, which makes $X\simplicial\to Y\simplicial$
into a groupoid $G$-extension.

The above discussion can be summarized  by  the following

\begin{prop}
Assume that 
 $\gendex{U_i}{i \in I}$ is a  good open covering of a manifold  $N$.
Then there is a one-one correspondence between 
$G$-extensions of  the  \v{C}ech groupoid
$  \coprod_{ij} U_{ij}
 \rightrightarrows \coprod_{i} U_i  $  and
the data $(\lambda_{ij}, g_{ijk})$, where
$\lambda_{ij}: U_{ij}\to \Aut(G)$ and
$g_{ijk}:U_{ijk} \to G$ satisfy
Eqs.~\eqref{eq:cocycle1}--\eqref{eq:cocycle2}.
\end{prop}

The data $(\lambda_{ij}, g_{ijk})$
is called a  non abelian 2-cocycle,
which coincides with the one  in \cite{Breen, Dedecker, Moer:regular}.

\begin{rmk}
Non abelian 2-cocycles have recently appeared in
many places in deformation quantization theory.
See \cite{Tsygan, Tsygan1, Kashiwara, KS1, KS2, K2, AP, PS, Y1, Y2}.
\end{rmk}

Now let us consider the band of this $G$-extension.
According to Definition \ref{def:band}, the band is an $\Out (G)$-torsor over the \v{C}ech groupoid $Y\simplicial :\coprod_{ij} U_{ij} \rightrightarrows \coprod_{i} U_i$, which is Morita equivalent to the manifold $N$. Thus the band is a principal $\Out(G)$-bundle over
the manifold $N$. To describe it more explicitly,
by $\lambdab_{ij}$ we denote the composition
of $\lambda_{ij}:U_{ij}\to\Aut(G)$ with the projection
$\Aut(G)\to\Out(G)$. Then Eq.~\eqref{eq:cocycle1} implies that
$$\lambdab_{ij}\rond\lambdab_{jk}\rond\lambdab_{ki}=1 .$$
In other words, $\lambdab_{ij}: U_{ij}\to \Out (G)$
is  a \v{C}ech 2-cocycle, which defines an $\Out(G)$-principal bundle 
over $N$. This is exactly the band of the $G$-extension.

\begin{prop} 
\label{lem:trivial}
Let $\gendex{U_i}{i\in I}$ be a good open covering of a manifold $N$. 
\begin{enumerate} 
\item A $G$-extension $X_1\xrightarrow{\phi}\coprod_{ij}U_{ij}\toto\coprod_{i}U_i$ 
of the \v{C}ech groupoid 
$\coprod_{ij}U_{ij}\toto\coprod_{i}U_i$ has a  trivial band 
if and only if there exists a trivialization of the kernel 
$\kernel \isomorphism \coprod_{i}U_i\times G$
and a section $\rho$ of  $X_1\xrightarrow{\phi}\coprod_{ij}U_{ij}$
 such that the associated data $(\lambda_{ij},g_{ijk})$
 satisfies $\lambda_{ij}=\idn$ and $(g_{ijk})\in\check{Z}^2\big(N,Z(G)\big)$.
\item Moreover, if $(g_{ijk})$ is a coboundary, then the section
 $\rho$ can be modified so that $\lambda_{ij}=\idn$ and $g_{ijk}=1$.
That is, the $G$-extension is isomorphic to the trivial $G$-extension.
\end{enumerate}
\end{prop}
\begin{proof}
\begin{enumerate}
\item Let $\bar{\eta}\in\sect{\Out(\kernel,G)\to\coprod U_i}$ be a trivialization of the band. Since the $U_i$'s are contractible, $\etab$ can be lifted to a section $\eta$ of $\Iso(\kernel,G)\to\coprod U_i$. This gives the desired trivialization of the kernel. 
According to Theorem \ref{thm:3.9}, since the band of the $G$-extension
 is trivial, and  all $U_{ij}$'s are contractible, there exists a section
 $\rho$ of $X_1\xrightarrow{\phi} 
\coprod U_{ij}$ such that \begin{gather} \rho(x_{ij})\cdot g_{x_j}=g_{x_i}\cdot \rho(x_{ij}) \label{spike}
\\ \rho(x_{ii})=1 \nonumber \\ \rho(x_{ji})=\rho(x_{ij})\inv \nonumber \end{gather} 
for all $x_{ij}\in U_{ij}$ and $g\in G$.
Now, we identify $X_1$ with
 $\coprod U_{ij}\times G$ through $$\coprod U_{ij}\times G\to X_1:(x_{ij},g)\mapsto \rho(x_{ij})\cdot g_{x_j}. $$ Then Eq.~\eqref{spike} becomes $(x_{ij},1)(x_{j},g)=(x_{i},g)(x_{ij},1)$. Hence $\lambda_{ij}=\idn$.
From $(x_{ij},1)(x_{jk},1)=(x_{ik},g_{ijk})$, it follows that $\rho(x_{ki})\rho(x_{ij})\rho(x_{jk})=g_{ijk}$.
By  Eq.~\eqref{spike}, we have $g_{ijk}\in Z(G)$.
\item Assume  that $(g_{ijk})$ is a coboundary: $\rho(x_{ki})\rho(x_{ij})\rho(x_{jk})=
g_{ijk}=h_{jk}h_{ik}\inv h_{ij}$, where  $h_{ij}:U_{ij}\to Z(G)$. Using
Eq.~\eqref{spike}, this can be rewritten as 
$\rho(x_{ki})h_{ki}\inv\cdot\rho(x_{ij})h_{ij}\inv\cdot\rho(x_{jk})
h_{jk}\inv=1$. 
Define a new section 
$\rho':\coprod U_{ij}\to X_1:x_{ij}\mapsto \rho(x_{ij})h_{ij}\inv$. It is easy to check that in the associated identification of $X_1$ with $\coprod U_{ij}\times G$, one has $\lambda'_{ij}=\lambda_{ij}\rond\AD_{h_{ij}\inv}=1$
 and $g_{ijk}'=\rho'(x_{ki})\rho'(x_{ij})\rho'(x_{jk})=1$.
\end{enumerate}
\end{proof}

As an immediate consequence, we  see that a $G$-extension with
trivial band over the \v{C}ech groupoid of a good cover
is completely determined by a \v{C}ech 2-cocycle in
$\check{Z}^2\big(N,Z(G)\big)$. Therefore
we have derived  the following result of Giraud \cite{Giraud}
by a direct argument.

\begin{cor} 
Isomorphism classes of $G$-bound gerbes over a manifold $N$
are in one-one correspondence with $H^2(N,  Z(G))$.
\end{cor}

\begin{prop} \label{prop:trivial}
Any groupoid $G$-extension of the \v{C}ech  groupoid associated to a good
 open covering of a contractible manifold is isomorphic to the
 trivial $G$-extension.
\end{prop}
\begin{proof}
Let $N$ be a contractible manifold, $(U_i )$  
a good covering of $N$,  and $X_1\to\coprod U_{ij}\toto\coprod U_i$
 a groupoid $G$-extension of the associated \v{C}ech groupoid. 
Its band must be trivial, since
 any principal bundle over a contractible manifold is trivial.
 According to Proposition \ref{lem:trivial},
 the multiplication on $X_1$ is
 entirely determined by the 2-cocycle
$(g_{ijk})$ in ${\check{Z}}^2\big(N,Z(G)\big)$.
 Since the manifold  $N$ is
 contractible and $(U_i )$ is a good covering,
 $(g_{ijk})$ must be  a 2-coboundary,
which implies that $X_1\to\coprod U_{ij}\toto\coprod U_i$ 
is isomorphic to the  trivial $G$-extension.
\end{proof}

By a \textbf{refinement} of
a $G$-extension  $X_1 \to \coprod_{i,j} U_{ij} \toto \coprod_j U_j$
 of a \v{C}ech groupoid, we mean the pull-back of this $G$-extension
 through a refinement of the covering $(U_j)_{j \in J}$ of $N$.

As an immediate consequence of Proposition \ref{prop:trivial},
we  have the following

\begin{cor} \label{lem:refinement}
Any groupoid $G$-extension of a \v{C}ech groupoid
over a contractible
manifold $N$ has a refinement which is isomorphic to a trivial $G$-extension.
\end{cor}

\section{Connections on groupoid extensions}

\subsection{Connections as horizontal distributions}

Recall that a \textbf{horizontal distribution} on a fiber bundle
 $X\xrightarrow{\phi} Y$ is an assignment to each point $x\in X$ of a subspace $H_x$ of $T_x X$ transversal  to the fiber of $\phi$ containing $x$.

Recall also that, to any Lie groupoid $\Gamma_1\rightrightarrows\Gamma_0$, one can associate its \textbf{tangent Lie groupoid} $T\Gamma_1\rightrightarrows T\Gamma_0$ whose structure maps are the differentials of those of $\Gamma_{\simplicial}$ \cite{MX}.
More precisely, if ${\bs}:\Gamma_1\to\Gamma_0$, $\bt:\Gamma_1\to\Gamma_0$ and $m:\Gamma_2\to\Gamma_1$ denote the source, target and multiplication maps of $\Gamma\simplicial$, then 
${\bs}_*:T\Gamma_1\to T\Gamma_0$ and ${\bt}_*:T\Gamma_1\to T\Gamma_0$ are respectively the source and target maps of $(T\Gamma)\simplicial$, and the multiplication map $T\Gamma_1\times_{{\bt}_*,\Gamma_0,{\bs}_*}T\Gamma_1\to T\Gamma_1$ is the composition of the canonical isomorphism 
$T\Gamma_1\times_{{\bt}_*,\Gamma_0,{\bs}_*}T\Gamma_1\isomorphism
 T(\Gamma_1\times_{{\bt},\Gamma_0,{\bs}}\Gamma_1)$
 with the differential of $m: \Gamma_1\times_{{\bt},\Gamma_0,{\bs}}\Gamma_1\xrightarrow{}\Gamma_1$.

\begin{defn}
A \textbf{connection} on a Lie groupoid extension
$X_1 \xrightarrow{\phi} Y_1 \rightrightarrows M$  is a
horizontal distribution $H$ on  $X_1 \xrightarrow{\phi} Y_1$ which is
also a Lie subgroupoid  of the tangent groupoid
 $TX_1\rightrightarrows TM$.
I.e.   we have
\begin{gather} H_{x_1}\cdot H_{x_2} \subset H_{x_1\cdot x_2} \ \ \ \  \text{for all $(x_1,x_2)\in X_2$}\label{3} \\
H_{x}\inv \subset H_{x\inv }   \ \ \ \  \text{for all $x \in X_1$}  \label{3b} .\end{gather}
\end{defn}

\begin{lem}
\begin{enumerate}
\item The distribution $H$ contains the entire unit space $TM$ of $X_1$.
\item $H_{x_1}\cdot H_{x_2} = H_{x_1\cdot x_2}$
\item $(H_x)\inv = H_{x\inv }$
\end{enumerate}
\end{lem}

\begin{proof}
\begin{enumerate}
\item  Since $H$ is a Lie  subgroupoid of $TX_1\toto TM$, it contains
its unit space, which is ${\bs}_* H$.
On the other hand, we have ${\bs}_* H={\bs}_* \phi_* H={\bs}_* TY_1=TM$.
The conclusion follows.
\item Since $ T_{{\bt}(x)}M \subset H_{x\inv } \cdot H_{x}$
and $H_{x_1}\inv \subset H_{x_1\inv }$, we
have
 $$H_{x_1 \cdot x_2} = H_{x_1 \cdot x_2} \cdot T_{{\bt}(x_2)} M \subset H_{x_1 \cdot x_2} \cdot H_{{x_2}\inv } \cdot H_{x_2} \subset H_{x_1} \cdot H_{x_2} .$$
\item Substituting $x\inv$ for $x$ in Eq.~\eqref{3b}, one gets $H\inv_{x\inv}\subset H_x$. Thus, $H_{x\inv}\subset H\inv_x$.
\end{enumerate}
\end{proof}

As before, let $\kernel$ be the kernel of $X\simplicial\xrightarrow{\phi}Y\simplicial$.
 Then $\kernel\to M$ is a group bundle,   which is also a Lie
 subgroupoid of
 $X\simplicial$. Set
\begin{equation}
\label{horker}
  \hkernel=T\kernel\cap H .
\end{equation}
Since $\ker(\phi_*) \subset T_k \kernel $ for any $k \in \kernel $,
$\hkernel$ defines a connection for the 
 bundle  $\kernel\xrightarrow{\phi} M$, which is also
a groupoid extension connection when $\kernel\to M$ is considered as a groupoid extension $\kernel\to M\toto M$.
Since this groupoid extension is simply a bundle of groups,
$\hkernel$ is a group bundle connection that we call
the \textbf{induced connection on the bundle of groups} $\kernel\to M$. 

Let $\liekernel\to M$ denote the Lie algebra
bundle  associated to the group bundle $\kernel\to M$.
Let $\exp:\liekernel\to\kernel$ denote the pointwise   exponential map.
The connection $\hkernel$ on the group  bundle  $\kernel\xrightarrow{\phi}M$
induces a connection on the Lie algebra bundle $\liekernel \to M$ that we call the
 \textbf{induced horizontal distribution $H^{\liekernel}$ on the bundle of Lie algebras}
$\liekernel\to M$
 defined as follows: $v\in T\liekernel$ is horizontal if $\exp_* v\in H^{\kernel}$.
 This connection is compatible with the Lie algebra bundle
 structure,
i.e.  it is given by a covariant derivative
on the vector bundle $\liekernel\to M$
and  satisfies Proposition \ref{prop:covderLAB}(4).

\subsection{Horizontal paths for groupoid extension connections}
\label{sec:4.2}

It is natural to ask what is the
geometrical meaning of  a connection on a groupoid extension.
This subsection is aimed  to answer this
question. 

Recall that, given an horizontal distribution on a fiber bundle,
a path  is said to be \textbf{horizontal} if it is tangent to the
horizontal distribution. Also recall that, given a groupoid, two paths
$\tau\mapsto\gamma_1(\tau),\;\tau\mapsto\gamma_2(\tau)$ on the space of 
arrows are said to be
\textbf{compatible} if, and only if, $\gamma_1(\tau),\gamma_2(\tau)$ are
composable  for all values of $\tau$. The \textbf{product} of these paths is the path $\tau\to\gamma_1(\tau)\cdot\gamma_2(\tau)$.

The following proposition gives an alternative definition of 
Lie groupoid extension connection.

\begin{prop} \label{pro:para}
Let $X_1 \xrightarrow{\phi} Y_1 \rightrightarrows M$ be a 
Lie groupoid extension.  An horizontal distribution on $X_1\xrightarrow{\phi}Y_1$
is a groupoid extension connection if and only if the product of any pair of horizontal composable paths in $X_1$ is still a horizontal
path, and the inverse of any horizontal path is still a horizontal path.
\end{prop}

%\begin{proof}
%This follows directly from the fact that for any pair of composable paths
% $(\gamma_1,\gamma_2)$ in $X_1$,
%$$\tfrac{d}{d\tau}(\gamma_1\cdot\gamma_2)=\tfrac{d}{d\tau}(\gamma_1)\cdot\tfrac{d}{d\tau}(\gamma_2) ,$$ and
%$$\tfrac{d}{d\tau}\big(\gamma\inv \big) =\big(\tfrac{d}{d\tau}\gamma\big)\inv .$$
%Here the right hand sides refer to the multiplication and the inverse in $TX_1\toto TM$.
%\end{proof}

Here is yet another alternative description of groupoid extension connections.

For any manifold $N$, denote $N^I$ the set of smooth paths from $[0,1]$ to $N$. If $\Gamma_1\toto\Gamma_0$ is a groupoid, then $\Gamma_1^I\toto \Gamma_0^I$ inherits a groupoid structure with 
the source map $\gamma\mapsto \bs\rond\gamma$,
the  target map $\gamma\mapsto \bt\rond\gamma$,
 and the  product $\tau\mapsto \gamma_1(\tau)\cdot\gamma_2(\tau)$,
 for any composable  $\gamma_1,\gamma_2\in\Gamma_1^I$. We call this groupoid the \textbf{path groupoid}.
Proposition \ref{pro:para} can be reinterpreted as follows: a horizontal distribution is a groupoid extension connection if, and only if, the horizontal paths form a subgroupoid of the path groupoid $X_1^I\toto M^I$.

Recall that an
 \textbf{Ehresmann connection} on a fiber bundle $X\xrightarrow{\phi} Y$ is a horizontal distribution satisfying the following
technical assumption: any path on $Y$ starting from $y\in Y$ has a unique horizontal lift on $X$ starting from a given point $x$ in the fiber of $X$ over $y$. This additional assumption is required to avoid horizontal lifts going to infinity in a finite time.

Let $\gamma$ be a path in $Y_1$. The parallel transport $\tau^{\gamma}_t$ along $\gamma$ is a transformation satisfying the following properties: 
\begin{itemize}
\item $\tau^{\gamma}_0(x)=x$;
\item for any $x\in X_1$ with $\phi(x)=\gamma(s)$, the relation  $\phi(\tau^{\gamma}_{t}(x))=\gamma(s+t)$ holds;
\item the path $t\mapsto \tau^{\gamma}_{t}(x)$ is horizontal.
\end{itemize}

\begin{prop}
Groupoid extension connections are Ehresmann connections.
\end{prop}

\begin{proof}
Let $X_1\xrightarrow{\phi}Y_1\toto M$ be a groupoid extension endowed with a groupoid extension connection.
Let us fix a path $\gamma:t\mapsto\gamma_t$ in $Y_1$ defined on an open interval $I$ of $\R$ containing $[0,1]$. For each $x\in X_1$ in the fiber over $\gamma_0$, we will denote by $\gammab^x:t \to \gammab^x_t$ the path
(if any) in $X_1$, which passes through $x$ and lifts $\gamma$ horizontally. 

\fbox{Step 1} 
Given a point $b$ in $\phi\inv(\gamma_0)$, we will show that there exists a positive real number $\epsilon$ such that, for all points $x$ in the same connected component as $b$ in $\phi\inv(\gamma_0)$, the path $\gammab^x$ is defined for $t\in (-\epsilon,\epsilon)$.

Indeed, since the solutions of ODE's depend smoothly on the parameters, there exists a connected neighborhood $U$ of $b$ in its fiber $\phi\inv(\gamma_0)$ and a positive real number $\epsilon$ such that, for all $u\in U$, the 
horizontal lift $\gammab^u$ is defined for all $t\in 
(-\epsilon,\epsilon)$.

For each $u\in U$, define the path $\delta^{b\inv\cdot u}\in
C^{\infty}((-\epsilon,\epsilon),  \kernel )$ by 
$\delta^{b\inv\cdot u}_t=(\gammab^b_t)\inv\cdot(\gammab^u_t)$. Since the horizontal distribution on $X_1$ is compatible with the groupoid multiplication, the paths $\delta^{b\inv\cdot u}$ with $u\in U$ as well as all products of 
such paths in the groupoid $C^{\infty}((-\epsilon,\epsilon), X_1 )$ 
are horizontal. 

Observe that $\big\{\delta^{b\inv\cdot u}_0|u\in U\big\}$ is the neighborhood $L_{b\inv} U$ of the identity in the Lie group $\kernel_{{\bt}(x)}$. 
Now, taking any element $x$ in the same connected component of $\phi\inv(\gamma_0)$ as $b$, the product $b\inv\cdot x$ lies in the connected component of the Lie group $\kernel_{{\bt}(x)}$ containing its unit element.
Hence $b\inv\cdot x$ can be written as a finite product $k_1\cdot k_2\cdot \cdots \cdot k_n$ of elements in
the neighborhood $L_{b\inv} U$.

Therefore, the horizontal lift $\gammab^x$ is defined for all $t\in (-\epsilon,\epsilon)$ because 
$$\gammab^x_t=\gammab^b_t\cdot\delta^{k_1}_t\cdot\delta^{k_2}_t\cdot\cdots\delta^{k_n}_t.$$ 

\fbox{Step 2} 
We  show that for any point $x$ in the connected component of $b$ in the fiber $\phi\inv(\gamma_0)$, the path $\gammab^x$ is defined for all $t\in [0,1]$. 

By the first step, it is clear that there exists an open covering
$$\bigcup_{s\in I ,\; U \in \mathcal{C}_s}\genrel{(s+t,\tau^{\gamma}_{t}(u))}{t\in(-\epsilon_U,\epsilon_U),\; u\in U}$$ 
of $I\times_{\gamma,Y_1,\phi}X_1$,
where $\mathcal{C}_s$ denotes the set of connected components of the fiber $\phi\inv(\gamma_s)$, $\epsilon_U$ depends on $U$
and $\tau$ denotes the parallel transport. 
Any connected component of $I\times_{\gamma,Y_1,\phi}X_1$ is thus a union of such open sets. Any horizontal lift of $\gamma$ does
 entirely lie in such a connected component. Consider the connected component $B$ of $I\times_{\gamma,Y_1,\phi}X_1$ containing $b$. 
Since the parallel transport maps fibers homeomorphically, it does preserve the connected components of the fibers. Hence the intersection of $B$ with any fiber must be a connected component of this fiber. It is thus clear that 
$$B=\bigcup_{s\in I}\genrel{(s+t,\tau^{\gamma}_{t}(u))}{t\in(-\epsilon_s,\epsilon_s),\; u\in B\cap\phi\inv(\gamma(s))},$$ 
where $\epsilon_s$ does now only depend on $s$.

Since $\bigcup_{s\in I}(s-\epsilon_s,s+\epsilon_s)$ is an open covering of the compact interval $[0,1]$, there exists a finite subcovering 
$\bigcup_{i=0,\cdots,n}(s_i-\epsilon_{s_i},s_i+\epsilon_{s_i})$ of $[0,1]$ with $0=s_0<s_1<s_2<\cdots <s_n=1$ and 
$(s_{i-1}-\epsilon_{s_{i-1}},s_{i-1}+\epsilon_{s_{i-1}})\cap(s_i-\epsilon_{s_i},s_i+\epsilon_{s_i})\ne\emptyset$.
Let $t_0=0$, $t_{n+1}=1$ and choose $n$ real numbers $t_1,\cdots,t_n$ such that $t_i\in 
(s_{i-1}-\epsilon_{s_{i-1}},s_{i-1}+\epsilon_{s_{i-1}})\cap(s_i-\epsilon_{s_i},s_i+\epsilon_{s_i})$. 
Then, for any $x$ in the same connected component of $\phi\inv(\gamma_0)$ as $b$, the path recursively defined by $\gammab^x_t=
\tau^{\gammab^x_{t_i}}_{t-t_i} \quad \text{for}\; t_i\le t\le t_{i+1}$
 is the horizontal lift of $\gamma$ through $x$.
Modifying the choice of the point $b$ and the path $\gamma$ in
 the two steps above,  the conclusion follows. 
\end{proof}

The following result is an immediate consequence of Proposition \ref{pro:para} applied to $\hkernel$.

\begin{cor}
Let $\hkernel$ be the induced connection on the group bundle $\kernel\to M$. Then the parallel transport in $\kernel\to M$ preserves the group structure on the fibers.
\end{cor}

Since  the horizontal distribution on $\kernel$ is the image of the
 horizontal distribution on $\liekernel$ under the differential of 
the exponential map, we have the following lemma.

\begin{lem}\label{lem:expopreserveshoriz}
The following diagram 
\begin{equation} \xymatrix{ \liekernel_{\gamma(0)} \ar[d]_{\exp} \ar[r]^{\tau^{\gamma}_t} & \liekernel_{\gamma(t)} \ar[d]_{\exp} \\ 
\kernel_{\gamma(0)} \ar[r]^{\tau^{\gamma}_t} & \kernel_{\gamma(t)} } \label{exptau} \end{equation}
commutes, where $\tau^{\gamma}_t$ stands, with an abuse of notation, for both the parallel transportation  in $\kernel $ and $ \liekernel$
over some path $\gamma$ in $M$.
\end{lem}

The following proposition shows that 
any groupoid extension of a  connected Lie groupoid that
 admits a  connection must be a  groupoid $G$-extension 
for some fixed Lie group $G$.

\begin{prop}
Let $X_1\xrightarrow{\phi}Y_1\rightrightarrows M$ be a Lie groupoid extension with kernel $\kernel $. Assume that the orbit space $M/Y_1$ is path connected. If there exists a groupoid extension connection, then
\begin{enumerate}
\item the groups $\kernel _m$, $m\in M$, are all isomorphic;
\item $\phi$ is a Lie groupoid $G$-extension for a
fixed  Lie group $G$.
\end{enumerate}
\end{prop}

\begin{proof}
\begin{enumerate}
\item Since $M/X_1$ is path connected, any two points $m, n\in M$ can be connected by a family of points $\gendex{m_i}{i=1,\dots,l}$ such that any two consecutive elements in the family are either on the same connected component of $M$, or are the source and target of some element in $X_1$.

For any pair of points $m$ and $n$ in the same connected component of $M$, an isomorphism $\kernel _m\isomorphism\kernel _n$ can be constructed by taking the parallel transport over a path joining $m$ and $n$. And if $x\in X_1$, the conjugation by $x$ induces an automorphism $\kernel _{ {\bs}(x)}\isomorphism\kernel _{{\bt}(x)}$.
\item Choose a Riemannian metric on $M$ and take an open covering of $M$ by contractible open normal neighborhoods $\gendex{U_i}{i\in I}$. 
Then $h_i:U_i\times[0,1]\to U_i:(\exp_{u_i} \xi,t)\mapsto \exp_{u_i} t\xi$,
where $\xi$ is in a small neighborhood of zero in $T_{u_i} U_i$,
 is a smooth deformation retraction of $U_i$ onto a fixed point $u_i$. Consider the pullback
$X_1[\coprod U_i]\to Y_1[\coprod U_i]\toto \coprod U_i$
of $X_1\to Y_1\toto M$.
%$$\xymatrix{ X_1[\coprod U_i] \ar[d] \ar[r] & X_1 \ar[d] \\
%  Y_1[\coprod U_i] \dar[d] \ar[r] & Y_1 \dar[d] \\ \coprod U_i
%  \ar[r] & M .}$$
Its kernel is $\coprod_{i\in I} U_i\times_M \kernel $. As each $U_i\times_M\kernel $ can be locally identified to $U_i\times\kernel _{u_i}$ using parallel transport along the paths $t\mapsto h_i(m,t)$, $m\in U_i$, as above, it follows that $\coprod_{i\in I} U_i\times_M \kernel $ can be identified with $\coprod_{i\in I} U_i\times G$ for a fixed
 Lie group $G$. Hence $X_1[\coprod U_i]\to Y_1[\coprod U_i]\rightrightarrows \coprod U_i$ is a Lie groupoid $G$-extension.
\end{enumerate}
\end{proof}

From now on, if a Lie 
 groupoid extension admits a groupoid extension connection, it is
 always  assumed to be a $G$-extension.
Take an arbitrary path $\gamma$ in $X_1^I$. The horizontal distribution being an Ehresmann connection, there exists an unique horizontal path $\gammab$ starting at $\gamma(0)$ and satisfying $\phi\rond\gamma=\phi\rond\gammab$.
There is therefore a unique map $g:[0,1]\to\kernel$ such that $\gammab(\tau)=\gamma(\tau)\cdot g_\gamma(\tau)$; note that, by construction, $g_\gamma(\tau)\in\kernel_{{\bt}\rond\gamma(\tau)}$ for all $\tau\in [0,1]$.

We call \textbf{right holonomy} of a path $\gamma$, denoted by
$\hol(\gamma)$,  the element $g_{\gamma}(1)\in\kernel_{{\bt}\rond\gamma(1)}$. 
The left holonomy can be defined similarly using the left action of $\kernel$.

\begin{prop}
\label{pro:holo}
For any pair of paths $\gamma_1,\gamma_2$ composable in $X_1$, the following relation holds,
$$\hol(\gamma_1\cdot\gamma_2)= \big(\AD_{(\gamma_2(1))\inv}\hol(\gamma_1)\big)\cdot\hol(\gamma_2) .$$
\end{prop}
\begin{proof}The paths $\gamma_1,\gamma_2$ being composable, so are the associated horizontal paths $\gammab_1,\gammab_2$. The product of horizontal paths being horizontal, their product $\gammab_1\cdot\gammab_2$ is an horizontal path. Moreover, this path starts at the point 
$\gammab_1(0)\gammab_2(0)=(\gamma_1\cdot\gamma_2)(0)$. Therefore $\overline{\gamma_1\gamma_2}=
\gammab_1\cdot\gammab_2$. It is now a simple matter to check that:
\begin{align*}(\gamma_1\gamma_2)(1)\cdot\hol_{\gamma_1\gamma_2}&=
\overline{\gamma_1\gamma_2}(1)\\ &
%=\gammab_1(1)\cdot\gammab_2(1) \\ &
=\gamma_1(1)\cdot\hol(\gamma_1)\cdot\gamma_2(1)\cdot\hol(\gamma_2) \\
&=\gamma_1(1)\cdot\gamma_2(1)\cdot\big(\AD_{(\gamma_2(1))\inv}\hol(\gamma_1)\big)\cdot\hol(\gamma_2) .
\end{align*}
\end{proof}

\subsection{Connections as 1-forms}
\label{sec:4.3}
Recall that any Lie groupoid $X_1\rightrightarrows M$ gives rise to a simplicial manifold
\begin{equation} \label{sim.ma}
\xymatrix{
\cdots \ar@<-1.5ex>[r]\ar@<-.5ex>[r]\ar@<.5ex>[r]\ar@<1.5ex>[r] & X_{2}
\ar@<1ex>[r]\ar[r]\ar@<-1ex>[r] & X_{1}\ar@<-.5ex>[r]\ar@<.5ex>[r] & X_{0}\,,}
\end{equation}
where $$X_n=\{(x_1,\dots,x_n)|{\bt}(x_i)= {\bs}(x_{i+1})\; i=1,\dots,n-1\}$$
is the set of composable $n$-tuples of elements of $X_1$, and $X_0=M$ and the face maps are defined as follows \cite{Segal}.
The maps $\epsilon^n_i:X_n\to X_{n-1}$ are given by, for $n>1$,
\begin{align*} & \epsilon_0^n(x_1,x_2,\dots,x_n)=(x_2,\dots,x_n) \\
& \epsilon_n^n(x_1,x_2,\dots,x_n)=(x_1,\dots,x_{n-1}) \\
& \epsilon_i^n(x_1,x_2,\dots,x_n)=(x_1,\dots,x_i x_{i+1},\dots,x_n), \; \; 1\le i \le n-1 ,
\end{align*} and, for $n=1$, $\epsilon_0^1(x)= {\bs}(x)$, $\epsilon_1^1(x)={\bt}(x)$. They satisfy the simplicial relations
$$\epsilon_i^{n-1}\rond\epsilon_j^n=\epsilon_{j-1}^{n-1}\rond\epsilon_i^n \quad \forall i<j .$$
We also define the maps $\bs:X_{n}\to M:(x_1,\dots,x_n)\mapsto  {\bs}(x_1)$ and $\bt:X_{n}\to M:(x_1,\dots,x_n)\mapsto {\bt}(x_n)$.
When $n=1$, we recover the source and target of the groupoid, justifying the notation.

Consider a Lie  groupoid extension \eqref{eq:ext}.
%$$ \xymatrix{
%\kernel \ar[r] \ar@<-0.5ex>[d]_{\pi} \ar@<0.5ex>[d]^{\pi} & X_1 \ar[r]^{\phi} \ar@<0.5ex>[d]^{\bt} \ar@<-0.5ex>[d]_{\bs} & Y_1 \ar@<0.5ex>[d]^{\bt} \ar@<-0.5ex>[d]_{\bs} \\ M \ar@2{-}[r] & M \ar@2{-}[r] & M } .$$
Since $X_1\toto M$ acts on $\kernel \to M$ by conjugation, it acts on $\liekernel \to M$ by adjoint action. Therefore one 
obtains a left representation of
the groupoid $X_1\toto M$ on $\liekernel \to M$. For any $x\in X_1$, the adjoint representation is denoted by  $\Ad_x: \liekernel _{{\bt}(x)}\to \liekernel _{ {\bs}(x)}$, which, by definition, is the 
derivative at the identity of the 
conjugation $\AD_x: \kernel _{{\bt}(x)} \to \kernel _{ {\bs}(x)}$.
Therefore one can talk about Lie groupoid cohomology with values in
$\liekernel $ \cite{Mackenzie}. Here a $\liekernel $-valued cochain is
a smooth map which associates an element in $\liekernel _{
  {\bs}(x_1)}$ to a composable $n$-tuple $(x_1,\dots,x_n)$. Thus the
space of $n$-cochains can be identified with
$C^{\infty}(X_n,\bs^*\liekernel )=\Gamma(\bs^*\liekernel \to X_n)$,
the space of smooth sections of the vector bundle 
$\bs^*\liekernel \to X_{n}$, i.e. the
pull back bundle of $\liekernel\to M$ via
$\bs: X_n\to M$.
 The differential $\ptl :C^{\infty}(X_{n-1},
\bs^*\liekernel )
\to C^{\infty}(X_{n},\bs^*\liekernel ) $ is 
given by
\begin{equation} \label{eq:left}
\ptl |_{(x_1,\cdots,x_{n})}=\Ad_{x_1} \rond {(\epsilon_0^n)}^* + \sum_{i=1}^n {(-1)}^i {(\epsilon_i^n)}^*.
\end{equation}

%re precisely, 
%\[\big(\ptl \rho\big)(x_1,\dots,x_{n}) = \Ad_{x_1} \big((\epsilon_o^n)^*\rho\big)(x_1,\cdots,x_n)
%+\sum_{i=1}^n (-1)^i \big((\epsilon_i^n)^*\rho\big)(x_1,\cdots,x_n),\]
%where $\rho\in \Gamma(\bs^*\liekernel \to X_n)$.

In general, as a cochain complex, 
one can consider $\Omega^l(X_{n},\bs^*\liekernel )$, the  space of differential forms on $X_{n}$ with values in the vector bundle
 $\bs^*\liekernel \to X_{n}$, i.e. 
the space of smooth sections of the vector bundle $\wedge^l T^*X_n  \otimes \bs^*\liekernel \to X_{n}$ and the operator $\ptl $ given by exactly the same formula \eqref{eq:left}
as the differential. Thus $(\ptl )^2=0$.
Similarly, by taking the inverse of the adjoint action of $X_1\toto M$ on $\liekernel \to M$, one obtains a right representation of $X_1\toto M$ on $\liekernel \to M$, and thus can consider the cochain complex formed by $\Omega^l(X_{n},\bt^*\liekernel )$, the space of differential forms on $X_{n}$ with values in the vector bundle $\bt^*\liekernel \to X_{n}$, 
and the differential $\ptr :\Omega^l(X_{n-1},\bt^*\liekernel)\to\Omega^l(X_{n},\bt^*\liekernel )$ given by
\begin{equation} \label{eq:right} \ptr |_{(x_1,\dots,x_n)}=\sum_{i=0}^{n-1} {(-1)}^i {(\epsilon_i^n)}^* + (-1)^n \Ad_{{(x_n)}\inv } \rond {(\epsilon_n^n)}^*. \end{equation}
We thus have $(\ptr )^2=0$.

In the sequel, if $\xi$ is an element of $\liekernel $, the right and left fundamental vector
 fields generated by $\xi$ will be denoted by $\xi^{\btr}$ and $\xi^{\btl}$ respectively.
Thus,
\begin{equation} \xi^{\btr}_x = 
\left.\tfrac{d}{d\tau}x\exp (\tau \xi)\right|_{\tau =0}, \qquad \qquad 
\xi^{\btl}_y = \left.\tfrac{d}{d\tau}\exp (\tau \xi)y\right|_{\tau =0},
 \label{eq:selection} \end{equation} 
where $x,y\in X_1$ with $\bt(x)=m=\bs(y)$ and $\xi\in\liekernel_m$.

\begin{defn}
Let $X_1\xrightarrow{\phi}Y_1\toto M$ be a Lie
groupoid extension. A 1-form $\alpha\in\Omega^1(X_1,\bt^*\liekernel )$ is said to be a \textbf{right connection 1-form} if $\alpha(\xi^{\btr})=\xi,\;\forall\xi\in\liekernel $ and $\ptr \alpha=0$. Similarly, a 1-form $\beta\in\Omega^1(X_1,\bs^*\liekernel )$ is said to be a \textbf{left connection 1-form} if $\beta(\xi^{\btl})=\xi,\;\forall\xi\in\liekernel $ and $\ptl \beta=0$.
\end{defn}

\begin{rmk}
\begin{enumerate}
\item
An $l$-form $\alpha\in\Omega^l(X_1,\bt^*\liekernel)$
 satisfies $\ptr \alpha=0$ if and only
if, for any composable pair $(x_1,x_2) \in X_2$ and for any $l$-tuples of
composable pairs $(u_1,v_1),\cdots,(u_l,v_l)$ in the tangent
groupoid $TX_1\toto TM$, with $u_i\in T_{x_1}X_1, v_i\in T_{x_2}X_1$
for all $i \in \{1,\cdots,l\}$, one has 
\begin{equation} \label{eq:connformintep}
 \alpha(u_1\cdot v_1,\cdots,u_l\cdot v_l)=
\alpha(v_1,\cdots,v_l)+\Ad_{x_2\inv }\alpha(u_1,\cdots,u_l) .\end{equation}
In particular, for $l=1$, one has  
\begin{equation}\label{eq:connformintepalpha} \alpha(u\cdot v)=\alpha(v)+
\Ad_{x_2\inv }\alpha(u) \end{equation}
for any composable $u\in T_{x_1}X_1$ and $v\in T_{x_2}X_1$.

Eq. \eqref{eq:connformintepalpha} can be interpreted as follows.
The tangent groupoid  $TX_1\toto TM$ acts  on $\liekernel\to M$
from the right  by $u\cdot k=Ad_{x}\inv k$ for any $u\in T_x X_1$
 and any $k\in \liekernel_{ {\bs}(x)}$. Then $\ptr \alpha=0$ if, and only if, $\alpha:TX_1\to \liekernel$
is a $1$-cocycle with respect to this action.
\item In the  case of $G$-extensions over a
 \v{C}ech groupoid, the condition $\ptr \alpha=0$ should  be
equivalent to the condition given by Breen-Messing in \cite[Eq.~(6.1.9)]{BreenMessing}. See \cite{BC}.
\end{enumerate}
\end{rmk}

\begin{lem}
\label{lem:pq}
Let $X_1\xrightarrow{\phi}Y_1\rightrightarrows M$ be a
Lie groupoid extension endowed with a Lie 
 groupoid extension connection $H$. 
Let $(p,q)$ be any point in $X_2$,
 and  $v_p\in T_p X_1$ and $v_q\in
T_q X_1$ any pair of  composable horizontal vectors.
\begin{enumerate}
\item If $\xi$ and $\eta$ are elements of the fibers
 of $\liekernel \to M$ at the points ${\bt}(p)$ and ${\bt}(q)$,
 respectively, then $(\xi^{\btr}_p+v_p)$ and $(\eta^{\btr}_q+v_q)$
are composable and
$$(\xi^{\btr}_p+v_p)\cdot (\eta^{\btr}_q+v_q)
={(\Ad_{q\inv }\xi +\eta )}^{\btr}_{pq}+(v_p\cdot v_q) .$$
\item If $\xi$ and $\eta$ are elements of the fibers of
$\liekernel \to M$ at the points $s(p)$ and $s(q)$ respectively, then 
$(\xi^{\btl}_p+v_p)$ and $ (\eta^{\btl}_q+v_q)$ are composable
and
$$(\xi^{\btl}_p+v_p)\cdot (\eta^{\btl}_q+v_q)
={(\xi+\Ad_{p}\eta)}^{\btl}_{pq}+(v_p\cdot v_q) .$$
\end{enumerate}
\end{lem}

\begin{proof} We prove (1) only since the argument for (2) is similar.
Choose a path $t\mapsto(\gamma(t),\delta(t))$
 in $X_2$ such that $\left.\tfrac{d}{dt}\right|_0\gamma(t)=v_p$
 and $\left.\tfrac{d}{dt}\right|_0\delta(t)=v_q$. Then,
$$ \xi^{\btr}_p+v_p = \left.\tfrac{d}{dt}\right|_0(\gamma(t)\exp(t\xi)), \ \ \ 
\eta^{\btr}_q+v_q = \left.\tfrac{d}{dt}\right|_0
(\delta(t)\exp(t\eta ) ). $$ 
Hence the result follows from the identity
$$\left.\tfrac{d}{dt}\right|_0
(\gamma(t)\exp(t\xi)\delta(t)\exp(t\eta) )=\left.\tfrac{d}{dt} \right|_0
\big( p\exp(t\xi)q\exp(t\eta)\big)+\left.\tfrac{d}{dt}\right|_0
(\gamma(t)\delta(t))  .$$
\end{proof}

\begin{thm}
Let $X_1\xrightarrow{\phi}Y_1\rightrightarrows M$ be a
Lie groupoid extension. Then the following
are all equivalent.
\begin{enumerate}
\item Lie  groupoid extension connections $H\subset TX_1$;
\item right connection 1-forms $\alpha \in \Omega^1(X_1,\bt^*\liekernel )$; and
\item left connection 1-forms $\beta \in \Omega^1(X_1,\bs^*\liekernel )$.
\end{enumerate}
\end{thm}

\begin{proof}
It suffices to prove the one-one correspondence
between (1) and (2). The equivalence between
(1) and (3) can be proved similarly.

\fbox{1$\Rightarrow$2} 
Define $\alpha\in \Omega^1(X_1,\bt^*\liekernel )$ by setting
 $$\alpha(\xi^{\btr})=\xi, \; \; \forall\xi  \in
\liekernel  \quad \text{and} \quad \alpha(v)=0, \; \; \forall v \in H.$$
From Lemma \ref{lem:pq}  and Eq.~\eqref{3}, it follows that
for all $(p, q)\in X_2$, $v_p \in H_p $ and $v_q \in H_q $
 with $v_p,v_q $ composable,
 $$\big((\xi^{\btr}_p+v_p)\cdot
 (\eta^{\btr}_q+v_q)\big)\ii \alpha=
\Ad_{q\inv }\xi +\eta
=\Ad_{q\inv }\big((\xi^{\btr}_p+v_p)\ii\alpha\big)+
(\eta^{\btr}_q+v_q)\ii\alpha.$$
Thus we have
$$ \ptr \alpha |_{(p, q)}
={(\epsilon_0^2)}^*\alpha - {(\epsilon_1^2)}^*\alpha+
\Ad_{q\inv }\rond{(\epsilon_2^2)}^*\alpha=0 .$$
Hence $\alpha$ is a right connection 1-form.

\fbox{2$\Rightarrow$1}
Set $H= \ker \alpha$. It is clear that $H$ is a
horizontal distribution  for the fiber bundle
$X_1\xrightarrow{\phi}Y_1$.
Since $\ptr \alpha=0$, $H$ is a subgroupoid of the tangent groupoid $TX_1\rightrightarrows TM$
by Eq.~\eqref{eq:connformintepalpha}. 
\end{proof}

The right and left connection 1-forms of a 
Lie groupoid extension connection are related
in a simple manner as described in the following

\begin{prop}
Let $X_1\xrightarrow{\phi}Y_1\rightrightarrows M$ be a Lie groupoid extension. And let $H\subset TX_1$ be a Lie 
groupoid extension connection, $\alpha\in\Omega^1(X_1,\bt^*\liekernel )$, and $\beta\in\Omega^1(X_1,\bs^*\liekernel )$ the right and left connection 1-forms associated to
$H$, respectively. Then $\alpha$ and $\beta$ are related by the following equation:
$\beta|_x=\Ad_x\rond\alpha|_x,\quad\forall x\in X_1$.
\end{prop}

\begin{proof}
This equation is obviously true on horizontal vectors. Moreover,
 the equation 
$\xi^{\btr}_x=(\Ad_x \xi )^{\btl}_x$,
$\forall \xi \in \liekernel $, 
implies the expected result on vertical vectors.
\end{proof}

\subsection{Connections on groupoid $G$-extensions}

In this subsection, we consider a Lie
 groupoid $G$-extension
$X\simplicial\xrightarrow{\phi}Y\simplicial$ whose  kernel
$\kernel$ is endowed with a trivialization
$\chi: M\times G\to \kernel$. Thus $X_1\xrightarrow{\phi}Y_1$
is a  principal  $G$-$G$ bibundle.
  Assume that there exists  a groupoid extension connection $H$
on $\phi$. A natural question is when this connection
is a connection for the right (or left) principal $G$-bundle
$X_1\xrightarrow{\phi}Y_1$.

According to Eq.~\eqref{horker},
 there is an induced connection $H^\kernel$ on
the kernel $\kernel\to M$. Via the trivialization $\chi$,  it in turn  induces
a connection on the group bundle $M\times G\to M$,
which is denoted by the symbol $H'$.
Then $H'$
is a connection for the trivial extension
$M\times G\to M\rightrightarrows M$.

For any $\xi \in \mathfrak{g}$,
$\xi^R$ denotes the right invariant vector field
on $G$ corresponding to $\xi$. 
Associated to $H'$, there exists a family $\{F_g\}_{g\in G}$
 of $\mathfrak{g}$-valued 1-forms on $M$ such that 
$(v_m,  {\big(F_g(v_m)\big)}^R_g) \in T_m M \times T_g G$ is the unique element
 of $H'_{(m,g)}$ whose projection in $T_m M$ is $v_m$. In other words, $\{F_g\}$ measures
the defect between  $H'_{(m,g)}$ and $T_{(m, g)}(M\times\{g\})$.
Note that, by definition
\begin{equation} \label{eq:consiFg}   F_g(v_m)  = -\beta{|_{(m,g)}} (g_* v_m),  \end{equation}
   where $g_*$ stands  for the differential of the constant section
 $ m \to (m,g) $  of $M\times G \to M$.

\begin{lem}
\label{lem:4.12}
The condition
 $$H'_{(m,g)}\cdot H'_{(m,h)}=H'_{(m,gh)}$$
 is equivalent to
 \begin{equation}
F_{gh}= F_g+ \Ad_g F_h, \ \ \forall g, h\in G \label{2} .
\end{equation}
I.e.
 $F: G\to \Omega^1(M)\otimes \mathfrak{g}$
 is a Lie  group 1-cocycle, where $G$ acts 
on $\Omega^1(M)\otimes \mathfrak{g}$ by adjoint action on the
second factor. 
In particular, $F_e=0$,
 where  $e$ denotes the unit element of $G$.
\end{lem}

\begin{proof}
It is simple to see that 
\begin{align*}
(v_m,  (F_g(v_m)\big)^R_g)\cdot (v_m,  (F_h(v_m))^R_h)&= 
(v_m ,  R_{h*} (F_g(v_m) )^R_g + L_{g*} ( F_h(v_m)^R_h )\\
&=(v_m , (F_g(v_m)+\Ad_g F_h (v_m))^R)_{gh}. 
\end{align*}
Thus the right hand side belongs to $H'_{(m,gh)}$ if, and only
if, Eq.~\eqref{2} is satisfied.
\end{proof}

\begin{cor}
\label{FS1}
Under the same hypothesis as in Lemma \ref{lem:4.12}, 
if $G=S^1$, then $F_g= 0$, $\forall g\in G$.
\end{cor}
\begin{proof}
 Since $G$ is abelian, Eq.~\eqref{2} becomes $F_{gh}=F_g+F_h$.
  Since there is no non-trivial group homomorphism between $S^1$
 and $\R$, $F$ must vanish.
In other words, for $S^1$-extensions,
$\{M\times \{e^{i\theta}\}\}_{\theta}$  must be
horizontal sections for the
group bundle $M\times S^1\to M$.
\end{proof}

\begin{prop}
\label{prop:4.14}
Let $X\simplicial\xrightarrow{\phi}Y\simplicial$ be a
Lie  groupoid $G$-extension
 whose  kernel $\kernel$ is endowed with a trivialization
$\chi: M\times G\to \kernel$. Assume that  $H$ is a  Lie
 groupoid extension connection 
with its associated right and left connection 1-forms
  $\alpha$ and $\beta$ respectively.
\begin{enumerate}
\item If $\xi_x \in H_x$, then 
$$R_{g*}\xi_x+{\big(\Ad_{g\inv } F_g({\bt}_* \xi_x)\big)}^{\btr}_{xg} \in H_{xg}$$ and 
$$L_{g*} \xi_x + {\big( F_g({\bs}_* \xi_x)\big)}^{\btl}_{gx}  \in H_{gx} .$$
\item The connection 1-forms satisfy 
$$R^*_g \alpha-\Ad_{g\inv }\alpha +{\bt}^*(\Ad_{g\inv } F_g)=0$$ and 
$$L^*_g \beta-\Ad_{g}\beta+  {\bs}^* F_g=0 .$$
\end{enumerate}
\end{prop}

\begin{proof}
\begin{enumerate}
\item It is simple to check that  $\forall \xi \in H_x$,
$$\xi_x\cdot \big({\bt}_* \xi_x,  R_{g*}F_g({\bt}_* \xi_x)\big)
=R_{g*}\xi_x+{\big( \Ad_{g\inv } F_g({\bt}_* \xi_x)\big)}^{\btr}_{xg},$$
where $\cdot$ on the left hand side stands for the groupoid
multiplication on $TX_1\toto TM$. It thus follows that
 $R_{g*}\xi_x+{\big(\Ad_{g\inv } F_g({\bt}_* \xi_x)\big)}^{\btr}_{xg} \in H_{xg}$.

Similarly,  $\forall \xi \in H_x$,
$$\big({\bs}_* \xi_x,  R_{g*}F_g({\bs}_* \xi_x)\big)\cdot \xi_x=
L_{g*}\xi_x+F_g({\bs}_* \xi_x)^{\btl}_{gx} . $$
Thus  $L_{g*} \xi_x + {\big( F_g({\bs}_* \xi_x)\big)}^{\btl}_{gx} 
 \in H_{gx}$.
\item Easy consequence of (1).
\end{enumerate}
\end{proof}

As an immediate consequence, we have the following

\begin{cor}
The following assertions are equivalent: 
\begin{enumerate}
\item $F_g=0, \; \; \forall g\in G$;
\item $M\times \{g\},\   \forall g\in G$, is horizontal. 
\item $H$ is right $G$-invariant ($R_{g*}H=H, \; \; \forall g \in G$),
i.e.  $\alpha\in \Omega^1(X_1 )\otimes \mathfrak{g}$ is
a connection one-form for the right $G$-principal bundle $X_1\to Y_1$;
\item $H$ is left $G$-invariant ($L_{g*}H=H, \; \; \forall g \in G$),
i.e.  $\beta\in \Omega^1(X_1 )\otimes \mathfrak{g}$ is
a connection one-form for the left $G$-principal bundle $X_1\to Y_1$.
\end{enumerate}
\end{cor}

By Corollary \ref{FS1}, we are led to the following result, as expected.

\begin{cor}
A connection on  a Lie  groupoid $S^1$-extension $X\simplicial\to Y\simplicial$
must be a principal left (and right) $S^1$-bundle connection.
\end{cor}

When the kernel of the groupoid extension 
$X_1\xrightarrow{\phi}Y_1\toto M$  is not 
identified to $M\times G$, Lemma \ref{lem:4.12} 
and Proposition \ref{prop:4.14} can be generalized by replacing the section $M\times\gendex{g}{} $ by any section $\sigma$ of $\kernel\to M$.

Define the transformations $R_{\sigma}:x_1\to x_1\cdot \sigma_{{\bt}(x_1)}$ and $L_{\sigma}:x_1\to \sigma_{ {\bs}(x_1)}\cdot x_1$ of $X_1$, and let $F_{\sigma}:=-\sigma^*\beta$. By construction, we have $F_{\sigma}\in\Omega^1(M,\liekernel )$.
According to Eq.~(\ref{eq:consiFg}), we recover the previous definition
 $F_g$ when $\kernel = M \times G$ and $\sigma(m)= (m,g)$, which justifies the notation.

\begin{lem}
The condition 
\begin{equation}
\hkernel_{\sigma_1}\cdot\hkernel_{\sigma_2}=\hkernel_{\sigma_1\cdot\sigma_2},
\quad\forall \sigma_1,\sigma_2\in\kernel\;\text{such that}\;\pi(\sigma_1)=\pi(\sigma_2)
\label{17}\end{equation} 
is equivalent to 
\begin{equation} 
F_{\sigma_1 \sigma_2}=F_{\sigma_1}+\Ad_{\sigma_1}F_{\sigma_2},
\quad\forall\sigma_1,\sigma_2\in\sect{\kernel} .
\label{22}\end{equation}
I.e. $F:\sect{\kernel}\to\Omega^1(M)\otimes\sect{\liekernel}$ is a 
group 1-cocycle, where $\sect{\kernel}$ acts on $\Omega^1(M)\otimes\sect{\liekernel}$ 
by the adjoint action on the second factor. 
In particular, $F_e=0$, where $e$ denotes the unit section of $\kernel$.
\end{lem}

\begin{proof}
The proof is similar to that of Lemma  \ref{lem:4.12} and is
omitted.
%
%Note that $\sigma_*v_m+\big(F_{\sigma}(v_m)\big)^{\btl}_{\sigma(m)} $ 
%is horizontal by construction, and is, indeed, the horizontal part of
%${\sigma}_*v_m$. Therefore, Eq.~\eqref{17} is satisfied if, and only if, 
%\begin{align*} & \Big({\sigma}_*v_m+\big(F_{\sigma}(v_m)\big)^{\btl}_{{\sigma}(m)}\Big)\cdot\Big({\sigma}'_*v_m+\big(F_{{\sigma}'}(v_m)\big)^{\btl}_{{\sigma}'(m)}\Big) \\
%=& ({\sigma}{\sigma}')_*v_m+R_{{\sigma}'*}\big(F_{\sigma}(v_m)\big)^{\btl}_{{\sigma}(m)}+L_{{\sigma}*}\big(F_{{\sigma}'}(v_m)\big)^{\btl}_{{\sigma}'(m)} \\
%=& ({\sigma}{\sigma}')_*v_m+\big(F_{\sigma}(v_m)+\Ad_{{\sigma}(m)}F_{{\sigma}'}(v_m)\big)^{\btl}_{({\sigma}{\sigma}')(m)}
%\end{align*}
%is horizontal; i.e. if, and only if, Eq.~\eqref{22} holds.
\end{proof}

Given $\sigma\in\sect{\kernel}$ and $x\in X_1$, we use the shorthand $\sigma x$ (resp. $x\sigma$) for $\sigma_{\bs(x)}\cdot x$ 
(resp. $x\cdot\sigma_{\bt(x)}$).

\begin{prop} \label{eq:rightequiv} Let $X_1\xrightarrow{\phi}Y_1\toto
  M$ be a groupoid $G$-extension endowed with a groupoid extension
  connection $H$ with associated right and left connection 1-forms
  $\alpha$ and $\beta$ respectively. Let ${\sigma}\in\sect{\kernel}$.
\begin{enumerate}
\item If $v_x\in H_x$, then
$$R_{{\sigma}*}v_x+{\big(\Ad_{{{\sigma}({\bt}(x))}\inv } F_{\sigma}({\bt}_*v_x)\big)}^{\btr}_{x{\sigma}}\in H_{x{\sigma}}$$ and
$$L_{{\sigma}*}v_x+{\big(F_{\sigma}({\bs}_*v_x)\big)}^{\btl}_{{\sigma}x}\in H_{{\sigma}x} .$$
\item The connection 1-forms satisfy
\begin{equation}R^*_{\sigma} \alpha-\Ad_{{\sigma}\inv }\alpha +{\bt}^*(\Ad_{{\sigma}\inv } F_{\sigma})=0 
 \label{fbeta} \end{equation}and
\begin{equation} L^*_{\sigma} \beta-\Ad_{{\sigma}}\beta+  {\bs}^* F_{\sigma}=0 . \label{fbeta'} \end{equation}
\end{enumerate}
\end{prop}

\begin{proof}
The proof is similar to that of Proposition  \ref{prop:4.14} and is
omitted.
%
%\begin{enumerate}
%\item \label{al1} It is simple to check that $\forall v_x\in H_x$,
%$$ v_x\cdot\Big({\sigma}_*({\bt}_*v_x)+\big(F_{\sigma}({\bt}_*v_x)\big)^{\btl}_{{\sigma}({\bt}(x))}\Big)=
%v_x\cdot\Big({\sigma}_*({\bt}_*v_x)+\big(\Ad_{{\sigma}({\bt}(x))}\inv F_{\sigma}({\bt}_*v_x)\big)^{\btr}_{{\sigma}({\bt}(x))}\Big) ,$$
%where $\cdot$ on both sides stand for the groupoid multiplication on $TX_1\toto TM$. 
%The left hand side of the previous equation being the product of two horizontal tangent vector
%is horizontal.
%It thus follows that $R_{{\sigma}*}v_x+{\big(\Ad_{{\sigma}\inv } F_{\sigma}({\bt}_*v_x)\big)}^{\btr}_{x{\sigma}}\in H_{x{\sigma}}$.

%Similarly, $\forall v_x\in H_x$,
%$$\Big({\sigma}_*({\bs}_*v_x)+\big(F_{\sigma}({\bs}_*v_x)\big)^{\btl}_{{\sigma}( {\bs}(x))}\Big)\cdot v_x=
%L_{{\sigma}*}v_x+{\big(F_{\sigma}({\bs}_*v_x)\big)}^{\btl}_{{\sigma}x}\in H_{{\sigma}x} .$$
%\item It is sufficient to check the formulas on horizontal and vertical vector fields separately.
%\end{enumerate}
\end{proof}

\subsection{Covariant derivative on the bundle of Lie algebras}
\label{sec:4.5}
Using  parallel transport on the Lie algebra
 bundle  $\liekernel\to M$, one can 
introduce  a covariant derivative on its space of sections: 
\begin{equation} (\nabla_{\gammadot}\xi)_{\gamma(0)}=\left.\tfrac{d}{dt}\tau^{\gamma}_{-t}(\xi_{\gamma(t)})\right|_0 \label{covder} ,\end{equation}
where $\gamma\in M^I$, $\xi\in\Gamma(\liekernel\to M)$ and $\tau^{\gamma}_{-t}$ denotes the parallel transport from $\liekernel_{\gamma(t)}$ to $\liekernel_{\gamma(0)}$ 
along the path $\gamma$.
The following proposition is obvious.

\begin{prop}\label{prop:covderLAB}
\begin{enumerate}
\item \label{POINT1} $(\nabla_{\gammadot}\xi)_{\gamma(0)}=\left.\tfrac{d}{dt}\right|_0 \left.\tfrac{d}{ds}\right|_0 \tau^{\gamma}_{-t}(\exp s \xi_{\gamma(t)})$;
\item $\nabla_X \xi$ is $C^{\infty}(M)$-linear in $X$ and $\R$-linear in $\xi$;
\item $\nabla_X (f\xi)=X(f)\cdot\xi+f\cdot\nabla_X \xi$; and
\item $\nabla_X\lie{\xi}{\eta}=\lie{\nabla_X \xi}{\eta}+\lie{\xi}{\nabla_X \eta}$,
\end{enumerate}
where $\xi,\eta\in \sect{\liekernel}$ and  $X\in \mathfrak{X}(M)$.
\end{prop}

\begin{prop} \label{prop:Fnabla}
$\big(\nabla_{X}\xi\big)_{x}=
\left.\tfrac{d}{ds}\right|_0 F_{\exp s\xi}(X_x)$.
\end{prop}
\begin{proof} First, one easily checks that 
$$\left.\tfrac{d}{dt}\tau^{\gamma}_{-t}(\exp s\xi_{\gamma(t)})\right|_0
=\Ver\big(\left.\tfrac{d}{dt}\exp s\xi_{\gamma(t)}\right|_0 \big)
=\big(\beta(\left.\tfrac{d}{dt}\exp s\xi_{\gamma(t)}\right|_0)\big)^{\btl}_{\exp s\xi_{\gamma(0)}} .$$
Second, one has 
\begin{align*}\big(\nabla_{\gammadot}\xi\big)_{\gamma(0)} &=
\left.\tfrac{d}{dt}\right|_0 \left.\tfrac{d}{ds}\right|_0 \tau^{\gamma}_{-t}(\exp s \xi_{\gamma(t)}) \\
&= \left.\tfrac{d}{dt}\right|_0 \left.\tfrac{d}{ds}\right|_0 \big( \tau^{\gamma}_{-t}(\exp s \xi_{\gamma(t)}) \cdot \exp(-s\xi_{\gamma(0)})\big)\cdot \exp(s\xi_{\gamma(0)}) \\ 
&= \left.\tfrac{d}{dt}\right|_0 \left( \left.\tfrac{d}{ds} \tau^{\gamma}_{-t}(\exp s \xi_{\gamma(t)}) \cdot \exp(-s\xi_{\gamma(0)}) \right|_0 + \xi_{\gamma(0)} \right) \\ 
%&= \left.\tfrac{d}{dt}\right|_0 \left.\tfrac{d}{ds}\right|_0 \tau^{\gamma}_{-t}(\exp s \xi_{\gamma(t)}) \cdot \exp(-s\xi_{\gamma(0)}) \\ 
&= \left.\tfrac{d}{ds}\right|_0 \left.\tfrac{d}{dt}\right|_0 \tau^{\gamma}_{-t}(\exp s \xi_{\gamma(t)}) \cdot \exp(-s\xi_{\gamma(0)}) \\ 
%&= \left.\tfrac{d}{ds}\right|_0 \left( R_{\exp(-s\xi_{\gamma(0)})_*} \left. \tfrac{d}{dt} \tau^{\gamma}_{-t}(\exp s\xi_{\gamma(t)}) \right|_0 \right) \\
&= \left.\tfrac{d}{ds}\right|_0 \left( R_{\exp(-s\xi_{\gamma(0)})_*} \big(\beta(\left.\tfrac{d}{dt}\exp s\xi_{\gamma(t)}\right|_0)\big)^{\btl}_{\exp s\xi_{\gamma(0)}} \right) \\
&= \left.\tfrac{d}{ds}\right|_0 \beta\left(\left.\tfrac{d}{dt}\exp s\xi_{\gamma(t)}\right|_0\right) \\
%&= \left.\tfrac{d}{ds}\right|_0 \beta\big(\exp(s\xi)_* (\gammadot(0))\big) \\
&= \left.\tfrac{d}{ds}\right|_0 F_{\exp s\xi}(\gammadot(0)) .
\end{align*}
\end{proof}

The connection $H^{\liekernel }$ defined 
on the Lie algebra
 bundle  $\liekernel\to M$ naturally induces a connection on the pull back
 bundle ${\bt}^*\liekernel\to X_1$,
% $$\xymatrix{ {\bt}^*\liekernel \ar[r]^{\hat{\bt}} \ar[d] & \liekernel \ar[d] \\ 
%X_1 \ar[r]_{\bt} & M ,}$$
 given by $H^{ {\bt}^*\liekernel}=(\hat{\bt}_*)\inv H^{\liekernel}$,
where $\hat{\bt}: {\bt}^*\liekernel\to \liekernel$ is
the projection.
The associated covariant derivatives are related by $$\nabla^{\bt}_v( {\bt}^*\xi)= {\bt}^*\big(\nabla_{{\bt}_* v}\xi\big), \quad \forall v\in TX_1,\; \forall \xi\in\sect{\liekernel},$$
where $\bt^* \xi  \in \Gamma(\bt^* \liekernel \to M)$ denotes 
the pull-back through $\bt$ of a section $\xi \in \Gamma(\liekernel \to M)$.
Similarly, we have the pull back connection $H^{ {\bs}^*\liekernel}$ on
$ {\bs}^*\liekernel \to X_1$, whose covariant derivative is denoted by
$\nabla^{\bs}$.
 
\begin{prop}
\begin{enumerate}
\item 
For every $\eta\in\sect{ {\bt}^*\liekernel\to X_1}$, we define an associated vertical vector field on $X_1$: 
$$\eta^{\wtr}_x:= \left.\tfrac{d}{d\tau}\right|_0   x\cdot \exp\big(\tau\eta(x)\big).$$
Then \begin{equation} \nabla^{\bt}_X\eta=(\derlie{\eta^{\wtr}}\alpha)(X)+\lie{\eta}{\alpha(X)}. \label{eq:sablier} \end{equation}
\item
For every $\eta\in\sect{ {\bs}^*\liekernel\to X_1}$, we define an
associated vertical vector field on $X_1$:
$$\eta^{\wtl}_x:= \left.\tfrac{d}{d\tau}\right|_0 \exp\big(\tau\eta(x)\big)\cdot x .$$
Then \begin{equation} \nabla^{\bs}_X\eta=(\derlie{\eta^{\wtl}}\beta)(X)-\lie{\eta}{\beta(X)} \label{eq:sablier'} . \end{equation}
\end{enumerate}
Here the  notations $\eta^{\wtr}_x$ and $\eta^{\wtl}_x$ generalize
those in   Eq.~\eqref{eq:selection}. More precisely, for any $\eta\in\Gamma(\liekernel\to M)$, we have
$(t^*\eta)^\wtr=\eta^\btr $ and $(s^*\eta)^\wtl=\eta^\btl$.
\end{prop}
\begin{proof}
We will  prove (1).  The argument for (2) is similar.

 Let $\xi$ be a section of $\liekernel\to M$.
Setting $\sigma=\exp(u\xi)$ in Eq.~\eqref{fbeta} and evaluating on a tangent vector $X_x$, we get 
$$\alpha(R_{\exp(\tau \xi)*}X_x)-\Ad_{\exp(\tau \xi)}\inv
\alpha(X_x)+\Ad_{\exp(\tau \xi)}\inv
F_{\exp(\tau \xi)}({\bt}_*X_x)=0 .$$ Differentiating with respect to $\tau$ at $u=0$ and using Proposition \ref{prop:Fnabla}, we 
obtain $$-(\derlie{\xi^{\btr}}\alpha)(X_x)+\lie{ {\bt}^*\xi}{\alpha(X_x)}
+(\nabla_{{\bt}_*X}\xi)_{{\bt}(x)}=0 ,$$
where  $\xi^{\btr}_x:=\left.\tfrac{d}{d\tau}\right|_0x\cdot \exp\big(\tau\xi({\bt}(x))\big)$.
Hence $\nabla^{\bt}_X {\bt}^*\xi =(\derlie{\xi^{\btl}}\alpha)(X)
+\lie{ {\bt}^*\xi}{\alpha(X)}$.

Now, for any function $f\in\cty(X_1)$, 
we have 
\begin{align*} \nabla^{\bt}_X(f\cdot {\bt}^*\xi)
=& X(f)\cdot {\bt}^*\xi+f \nabla^{\bt}_X( {\bt}^*\xi) \\
=& X(f)\cdot {\bt}^*\xi+f \big(\derlie{\xi^{\btr}}(\alpha(X))-\alpha(\lie{\xi^{\btr}}{X})+\lie{ {\bt}^*\xi}{\alpha(X)}\big) \\ 
=& \derlie{f\xi^{\btr}}\big(\alpha(X)\big)-\alpha(\lie{f\xi^{\btr}}{X})+\lie{f\cdot\bt^*\xi}{\alpha(X)} \\
=& (\derlie{f\xi^{\btr}}\alpha)(X)+[f\cdot\bt^*\xi,\alpha(X)]. \\
\end{align*}

The result follows from the fact that any section 
$\eta\in\sect{\bt^*\liekernel\to X_1}$
is a linear combination of sections of the type $f\cdot\bt^*\xi$
for $\xi\in\sect{\liekernel \to M}$ and $f\in\cty(X_1)$.
\end{proof}

\begin{rmk}
If $X$ is horizontal, $\alpha(X)=0$ and $\beta(X)=0$ and  therefore we have
\begin{equation} \nabla^{\bt}_X\eta=\alpha(\lie{X}{\eta^{\wtr}}), \qquad \forall\eta\in\Gamma(t^*\liekernel\to X_1)  \label{eq:parthor} \end{equation}
and
\begin{equation} \nabla^{\bs}_X\eta=\beta(\lie{X}{\eta^{\wtl}}), \qquad \forall\eta\in\Gamma(s^*\liekernel\to X_1) . \label{eq:parthor'} \end{equation}
\end{rmk}

\subsection{Ehresmann curvature}

In this subsection, we study the curvature of a connection $H$ on a 
Lie groupoid extension. We denote the horizontal and vertical 
parts of a vector $v\in TX_1$ by $\Hor(v)$ and $\Ver(v)$
respectively. Recall that the Ehresmann  curvature of a horizontal distribution $H$, on the bundle $X_1\xrightarrow{\phi}Y_1$, is the $2$-form on $X_1$, valued
 in the vertical space $\ker\phi_*$, which is  defined by
\begin{equation} \label{eq:defehrcurv}
\omega(u,v)=-\Ver\big(\lie{\Hor(\tilde{u})}{\Hor(\tilde{v})}\big),
 \end{equation}
 where $u,v\in T_x X_1$ and $\tilde{u}, \tilde{v}$ are vector fields on $X_1$ such that $\tilde{u}_x=u$ and $\tilde{v}_x=v$.
It is easy to check that the right hand side of Eq.~\eqref{eq:defehrcurv}
is well defined, i.e. independent of the choice of the vector fields $\tilde{u},\tilde{v}$.

Using the right (resp. left) action of $\kernel$ on $X_1$, one can identify the vertical space $\ker\phi_*\subset T_x X_1$ of the groupoid extension $X_1\xrightarrow{\phi}Y_1$ with $\bt^*\liekernel$ (resp. $\bs^*\liekernel$) via the right (resp. left) action of $\kernel$ on $X_1$. Therefore, the curvature of the connection $H$ can be seen as a 2-form either in $\Omega^2(X_1,\ker\phi_*)$, $\Omega^2(X_1,\bt^*\liekernel)$ or $\Omega^2(X_1,\bs^*\liekernel)$.
Note that for any Lie algebra bundle $\liekernel\to M$
 over a smooth manifold $M$,  there is a graded Lie bracket
$\Omega^k(M,\liekernel)\otimes\Omega^l(M,\liekernel) \xrightarrow{[\cdot,\cdot]} \Omega^{k+l}(M,\liekernel)$ given  by 
$[\omega_1\otimes a_1,\omega_2\otimes a_2]=(\omega_1\wedge\omega_2)\otimes[a_1,a_2]$, $\forall a_1, a_2 \in \Gamma(\liekernel),\; \omega_1 \in \Omega^k (M), \omega_2 \in \Omega^l(M)$.
In particular, for any $\alpha \in \Omega^1(M,\liekernel)$,
and any vector fields $X,Y\in\mathfrak{X}(M)$, we have
$\half[\alpha,\alpha](X,Y)=[\alpha(X),\alpha(Y)]$.

%Note also that this bracket satisfies the graded Jacobi identity:
%\begin{equation} \label{eq:gradedJacobi} 
%(-1)^{km}[\alpha,[\beta,\gamma]]+  (-1)^{lk}[\beta,[\gamma,\alpha]]+ (-1)^{ml}[\gamma,[\alpha,\beta]]=0 ,
%\end{equation}
%where $\alpha \in \Omega^k(M,\liekernel),\beta \in \Omega^l(M,\liekernel)$ 
%and $\gamma \in \Omega^m(M,\liekernel)$.

The following lemma  indicates that,
similar to the case of principal bundles,
the Ehresmann curvature can be computed using the
standard formula.

\begin{lem} \label{lem:curvexp}
\begin{enumerate}
\item
If $\alpha\in\Omega^1(X_1,\bt^*\liekernel)$ is a  left
 connection 1-form of a groupoid extension connection, 
then its curvature 2-form $\omega\in\Omega^2(X_1,\bt^*\liekernel)$
 is given by 
$$\omega=d^{\nabla^{\bt}}\alpha+\half \lie{\alpha}{\alpha} ,$$
where  $d^{\nabla^{\bt}}: \Omega\degree (X_1,\bt^*\liekernel)\to 
\Omega^{\scriptscriptstyle\bullet+1} (X_1,\bt^*\liekernel)$
 denotes the exterior covariant derivative induced by $\nabla^{\bt}$.
\item
Similarly, if $\beta\in\Omega^1(X_1,\bs^*\liekernel)$ is a  left
 connection 1-form of a groupoid extension connection, then its curvature 2-form $\omega\in\Omega^2(X_1,\bs^*\liekernel)$ is given by 
$$\omega=d^{\nabla^{\bs}}\beta-\half \lie{\beta}{\beta} ,$$
where  $d^{\nabla^{\bs}}: \Omega\degree (X_1,\bs^*\liekernel)\to
\Omega^{\scriptscriptstyle\bullet+1} (X_1,\bs^*\liekernel)$ denotes the exterior covariant
derivative induced by $\nabla^{\bs}$.
\end{enumerate}
\end{lem}
\begin{proof}
The proof is straightforward and is omitted.
\end{proof}

Now let us recall  a general fact regarding  horizontal distributions.
I.e. the Ehresmann curvature is the holonomy of infinitesimal loops
(see \cite[page 118]{Walschap}).
For any point $m$ in a manifold $N$ and any tangent vectors
$u,v \in T_m N$, we say that a smooth map $ C$ from a neighborhood of
$0 $ in $\R^2$ to a neighborhood of $m \in N$ is \emph{adapted to $(u,v)$} if
\[ C_* \big(\tfrac{\partial}{\partial x}_{|_0}\big)=u \qquad \text{and} \qquad C_* \big(\tfrac{\partial}{\partial y}_{|_0}\big)=v ,\]
where $x,y$ are the standard coordinates on $\R^2$.

Consider now, for any 
small enough $\epsilon_1,\epsilon_2$, the loop
$R_{\epsilon_1,\epsilon_2} : [0,1] \to \R^2$ obtained by turning
anti-clockwise along the rectangle with edges $(0,0),
(\epsilon_1,0),(\epsilon_1,\epsilon_2),(0,\epsilon_2) $. 
Define, for any (small enough) $ \epsilon_1,\epsilon_2$, a family of loops  
$L^C_{\epsilon_1,\epsilon_2}$   on $N$ by
\begin{equation} \label{eq:deflc}
L^C_{\epsilon_1,\epsilon_2}  = C \smalcirc  R_{\epsilon_1,\epsilon_2}.
\end{equation}

The curvature can be computed from the holonomy according to the
following formula:

 \begin{prop}
\label{pro:4.25}
 For any  $C$ adapted to $(u,v) $, we have
\begin{equation}
\label{eq:loopcurv} 
\tfrac{\partial^2}{\partial \epsilon_1 \partial \epsilon_2}_{|_{\epsilon_1=\epsilon_2=0}} 
Hol \big( L^C_{\epsilon_1,\epsilon_2}\big)   = \omega (u,v).
\end{equation}
\end{prop}

From now on, we shall work  only with the right-connection 1-form $\alpha$.
Below we list  two important identities  that $\omega$ satisfies,
which we call \textbf{Bianchi identities}.
They will turn out to be of fundamental importance in the sequel.

\begin{thm}\label{thm:bianchies}
Let $\omega\in \Omega^2 (X_1, {\bt}^*{ \mathfrak K})$ 
be the Ehresmann curvature for a Lie  groupoid
extension connection  $ \alpha \in \Omega^1 ( X_1,  {\bt}^* \liekernel )$. 
Then we have the Bianchi identities:
\begin{gather}
\label{eq:Bianchi1} d^{\nabla^t}\omega + [\alpha,\omega] = {\bt}^* \omega^{\liekernel} (\alpha) \\ 
\label{eq:Bianchi2} \ptr  \omega = 0 ,
\end{gather}
where  $\omega^{\liekernel} \in \Omega^2(M,\End(\liekernel))$
is the curvature of the induced connection $\nabla$ on the Lie algebra
bundle $\liekernel\to M$, and
$ {\bt}^* \omega^{\liekernel} (\alpha)\in  
\Omega^3(X_1, {\bt}^* \liekernel )$ is the ${\bt}^* \liekernel$-valued
3-form on $X_1$ obtained by composing  $  {\bt}^* \omega^{\liekernel} \in
\Omega^2 ( X_1,\End({\bt}^* \liekernel) )$ with
$\alpha \in \Omega^1 ( X_1,{\bt}^* \liekernel)$ under the natural
pairing 
$  \Omega^2 \big( X_1,\End({\bt}^* \liekernel) )
 \otimes   \Omega^1 \big( X_1,{\bt}^* \liekernel ) \to
 \Omega^3 ( X_1,  {\bt}^* \liekernel )$.
\end{thm}

%We need the following lemma.

%\begin{lem} \label{lem:basicdiffgeo}
%\begin{enumerate}
%\item Let $\liekernel \to M$ be a Lie algebra bundle over a 
%smooth manifold $M$, and $\nabla $ a connection
% of this Lie algebra bundle.  The covariant derivative $d^{\nabla}$ of this connection satisfies
%  \begin{equation} \label{eq:dercorLiealg}  d^{\nabla} \lie{\xi}{\eta}=
%\lie{d^{\nabla}\xi}{\eta}+(-1)^{k}\lie{\xi}{d^\nabla \eta} \end{equation}
%for any $\xi \in \Omega^k(M,{\liekernel}), \eta \in \Omega^l(M,{\liekernel})  $.
%\item For any connection $\nabla$ of a vector bundle $E \to M$, the square of the covariant derivative
%is given by $(d^{\nabla})^2 = \omega$, where $\omega \in \Omega^2 (M,\End(E))$
%is the curvature of $\nabla $, i.e.
%    $$ \omega(u,v)= \nabla_u \nabla_v -\nabla_v \nabla_u - \nabla_{[u,v]} .$$
%\item Let $E \to M$ be a vector bundle endowed with a connection $\nabla$,
%and $\pi: N \to M $ a  smooth map. The curvature 
%of the pulled-back connection $\nabla^{\pi}$ on $\pi^* E \to N$ is the pull-back of
%the curvature of the connection $\nabla$. 
%\end{enumerate}
%\end{lem}

\begin{proof}
Let us prove  Eq.~\eqref{eq:Bianchi1} first.
By Lemma \ref{lem:curvexp}(1), we have 
\begin{multline} \label{eq:bi1} 
(d^{\nabla^{\bt}} + \ad_{\alpha}) \omega = (d^{\nabla^{\bt}} + \ad_{\alpha})\smalcirc (d^{\nabla^{\bt}} + \half \ad_{\alpha})(\alpha) \\
=(d^{\nabla^{\bt}})^2 \alpha + \half d^{\nabla^{\bt}} [\alpha,\alpha]  +[\alpha, d^{\nabla^{\bt}} \alpha ]
+ \half [\alpha,[\alpha,\alpha]] . 
\end{multline}
The  graded Jacobi identity implies that $ [\alpha,[\alpha,\alpha]]=0$.
Since $\nabla^{\bt}$ is a connection
on the Lie algebra bundle $\bt^* \liekernel \to X_1 $ we have 
$\half d^{\nabla^{\bt}} [\alpha,\alpha] = [d^{\nabla^{\bt}} \alpha ,\alpha]= -[\alpha, d^{\nabla^{\bt}} \alpha ]  $.
Eq.~\eqref{eq:bi1} then becomes 
$(d^{\nabla^{\bt}}-\ad_{\alpha})\omega=(d^{\nabla^{\bt}})^2\alpha$.
On the other hand, the curvature of $\nabla^{\bt}$ is ${\bt}^*\omega^{\liekernel}$. 
Hence we have $(d^{\nabla^{\bt}})^2 \alpha={\bt}^*\omega^{\liekernel}(\alpha)$.
This concludes the proof of Eq.~\eqref{eq:Bianchi1}.

Let us prove Eq.~\eqref{eq:Bianchi2}. 
Denote by $p_1,m,p_2$ the three face maps from $X_2 $ to $X_1 $ 
(previously denoted by $\epsilon_0^1, \epsilon_1^1, \epsilon_2^1 $):
$p_1(x_1,x_2)=x_1, \  m(x_1,x_2)=x_1 x_2 , \   p_2(x_1,x_2)=x_2$.
We denote by $(u,v)$, with $u \in T_{x_1}X_1$, $v \in T_{x_2} X_1$
and ${\bt}_* u = {\bs}_*  v$, 
an element in $ T_{(x_1 , x_2)} X_2$. Note that $m_* (u,v)= u \cdot v$, 
where the dot on the right hand side  stands for the
multiplication in the tangent groupoid $TX_1 \toto TM $.

For any composable pair $(x_1,x_2) \in X_2$, and any $(u_1,u_2), (v_1,v_2) \in T_{(x_1,x_2)}X_2$, 
let us choose a smooth map $c$ from a neighborhood $\mathcal{U}$ of $0$ in
$\R^2$ to a neighborhood of $(x_1,x_2) \in X_2$ adapted to $ \big( (u_1
,u_2), (v_1,v_2)\big) $. Then $p_1 \smalcirc c$, $m \smalcirc c$ and $p_2 \smalcirc c$ are smooth maps from
 $\mathcal{U}$ to $X_1$ adapted to $(u_1,v_1)$ in $T_{x_1}X_1$,
 $(u_1  u_2 , v_1  v_2)$ in $T_{x_1\cdot x_2}X_1$ and $(u_2,v_2)$ in $T_{x_2}X_1$, respectively.

For any small enough $\epsilon_1,\epsilon_2$, the loops
$L^{p_1\smalcirc c}_{\epsilon_1,\epsilon_2}$ and
$L^{p_2\smalcirc c}_{\epsilon_1,\epsilon_2}$, as defined in Eq.~\eqref{eq:deflc}, 
are compatible and their product is precisely $L^{m\smalcirc c}_{\epsilon_1,\epsilon_2}$. 
According to Proposition \ref{pro:holo}, we have 
$$\Hol(L^{m\smalcirc c}_{\epsilon_1,\epsilon_2}) = 
\AD_{x_2\inv } \Hol(L^{p_2\smalcirc c}_{\epsilon_1,\epsilon_2}) 
\cdot \Hol(L^{p_1\smalcirc c}_{\epsilon_1,\epsilon_2}) .$$

Applying $\tfrac{\partial^2}{\partial\epsilon_1
\partial\epsilon_2}\big|_{\epsilon_1=\epsilon_2=0} $ to this equation
 and using Eq.~\eqref{eq:loopcurv}, one obtains
$$\omega (u_1 \cdot u_2, v_1 \cdot v_2 ) = \omega(u_2,v_2)+ \Ad_{x_2\inv }  \omega(u_1,v_1) .$$
The result now follows.
\end{proof}

%\begin{rmk}
%Eq.~\eqref{eq:Bianchi2} has the following nice interpretation. Since
%$H \toto TM$ is a  Lie groupoid, so is $\wedge^2 H \toto \wedge^2 TM $.
%The $2$-form $\omega$ can be considered as a vector bundle morphism
%\[\xymatrix{
%\wedge^2 H \ar[r]^{\omega^{\flat}} \ar[d]_{\bt_*} & \liekernel \ar[d] \\
%\wedge^2 TM \ar[r]_p & M}\]
%
%Define a right action of $\wedge^2 H \toto \wedge^2 TM $
%on $p^* \liekernel \to \wedge^2 TM$ by $\xi\cdot u=\Ad_{x\inv } \xi$, $\forall
%u\in \wedge^2 H|_x$, $x\in X_1$ and $ \xi \in p^* \liekernel_{\bs(x)}$. Then
%Eq.~\eqref{eq:Bianchi2} simply means that $\omega^{\flat}$ is a 1-cocycle.
%\end{rmk}

\begin{rmk}
The relation between
Theorem \ref{thm:bianchies} and Eqs.~(6.1.12)-(6.1.15) in \cite{BreenMessing}
is investigated in \cite{BC}.
\end{rmk}

\section{Induced connections on the band}

\subsection{Induced horizontal distributions on the band}

Let $X_1\xrightarrow{\phi}Y_1\toto M$ be a Lie 
groupoid $G$-extension endowed with a connection $H\subset TX_1$. The purpose of this subsection is to construct an induced connection on the band out of the connection $H$.
First of all, let us recall some basic notions.

\begin{defn}[\cite{LTX}]
Let $\Gamma_1\rightrightarrows\Gamma_0$ be a Lie groupoid with Lie algebroid $A$.
A connection on a principal $G$-bundle $P\xrightarrow{J}\Gamma_0$ over 
 $\Gamma_1\rightrightarrows\Gamma_0$ is a $G$-invariant horizontal distribution $H\subset TP$ satisfying the following two conditions:
\begin{enumerate}
\item for each $p\in P$, we have the inclusion $\hat{A}_p\subset H_p$, 
where $\hat{A}_p$ denotes the subspace of $T_p P$ generated by the infinitesimal action of the Lie algebroid $A\to \Gamma_0$;
\item the distribution $\gendex{H_p}{p\in P}$ on $P$ is preserved under the action of $U_{\loc}(\Gamma\simplicial)$,
 the  pseudo-group of local bisections \cite{Mackenzie} of
$\Gamma_1\rightrightarrows\Gamma_0$, which naturally acts on $P$ locally.
\end{enumerate}
\end{defn}

First we consider the right $\Aut(G)$-principal bundle $\Iso(\kernel,G)\to M$.
For any $g\in G$, let 
\[\ev_g:\Iso(\kernel,G)\to\kernel:f\mapsto f(g)\] 
be the evaluation map. Differentiating it with respect to $f$ yields, for all fixed
 $g\in G$, a map $$T_f\Iso(\kernel, G)\xrightarrow{\ev_{g*}}T_{f(g)}\kernel .$$

We define a horizontal distribution $\hiso$ on $\Iso(\kernel,G)\to M$ by 
\begin{equation} \label{eq:dis}
\hiso_f=\genrel{v\in T_f \Iso(\kernel,G)}{\ev_{g*}v\in\hkernel_{f(g)},\ \forall g\in G}\subset T_f \Iso(\kernel,G),
\end{equation}
where, as defined in Eq.~\eqref{horker},
 $\hkernel$ denotes the induced connection on the kernel $\kernel\to M$. 
It is obvious that $\hiso$
is invariant under the  $\Aut(G)$-action on $\Iso(\kernel, G)$.

Set $$\hout=\rho_*\hiso,$$
 where $\rho:\Iso(\kernel,G)\to\Out(\kernel,G)$ denotes the projection. It is clear that $\hout$ is a horizontal distribution
on the bundle $\Out(\kernel,G)\to M$.

The following lemma is immediate.

\begin{lem}
\label{lem:5.2}
\begin{enumerate}
\item The horizontal paths in $\Iso(\kernel,G)$ are the
 paths $f_t$ such that, for any $g\in G$, $\ev_g f_t=f_t(g)$ is a horizontal path in $\kernel$.
\item The horizontal paths in $\Out(\kernel,G)$ are the images 
of horizontal paths in $\Iso(\kernel,G)$ under the projection $\rho:\Iso(\kernel,G)\to\Out(\kernel,G)$.
\end{enumerate}
\end{lem}

\begin{prop}
The horizontal distribution $\hout$ defines a connection on the band.
\end{prop}
\begin{proof}
First, as $\hiso$ is invariant under the right $\Aut(G)$-action, $\hout$ is invariant under the right $\Out(G)$-action.

Second, $\forall m \in M$, let $f$ be any element in $\Iso(\kernel,G)|_m$ and $\bbt\mapsto\gamma(\bbt)
$ a path in $Y_1$ lying in the target fiber  $t\inv (m)$ over $m$. Consider the horizontal
 lift $\bbt\mapsto \bar{\gamma}(\bbt)$ of $\gamma$ in $X_1$. Fix $f\in \Iso(\kernel,G)|_m$. 
For all $g\in G$, the paths $\bbt\mapsto \bar{\gamma}(\bbt)f(g)(\bar{\gamma}(\bbt))\inv $ in 
$\kernel$ are horizontal paths. Therefore,
 the path $\bbt\mapsto\rho(\AD_{\bar{\gamma}(\bbt)}\rond f)$ in $\Out(\kernel,G)$ is horizontal by 
Lemma \ref{lem:5.2}.
 Since any element of  $\hat{A}_f$ is tangent to such a path at its origin, 
we have $\hat{A}_f\subset\hout_f$.

For any $f\in\Iso(\kernel,G)|_m$, let $m_\bbt$ be any path in $M$ starting at the point
 $m$. For any $g\in G$, by $\bar{m}^{f(g)}_\bbt$, we denote the horizontal lift of $m_\bbt$
 in $\kernel$ starting at the point $f(g)\in\kernel_m$. Let $f_\bbt$ be a path in $\Iso(\kernel,G)$ defined by $f_\bbt (g)
=\bar{m}^{f(g)}_\bbt, \ \forall g\in G$. By definition, $f_\bbt$ is the horizontal lift of $m_{\bbt}$ 
in $\Iso(\kernel,G)$ through the point $f$. Hence, $\rho(f_\bbt)$ is the horizontal lift of $m_\bbt$ 
in $\Out (\kernel , G)$ starting at the point $\rho(f)$.

To show that $\hout$ is preserved under the action of $U_{\loc}(Y\simplicial)$, it suffices to show that, 
for any $L\in U_{\loc}(Y\simplicial)$, $L\cdot\rho(f_\bbt)$ is still a horizontal path in $\Out(\kernel,G)$.
 For this purpose, let $\sigma_\bbt$ be the unique path in $L$ such that
 ${\bt}\rond\sigma_\bbt=m_\bbt$, and let $\bar{\sigma}_{\bbt}$ be any of its horizontal lifts on $X_1$. By definition, 
$$L\cdot\rho(f_\bbt)=\sigma_\bbt\cdot\rho(f_\bbt)=\rho\big(\AD_{\bar{\sigma}_\bbt}\rond f_\bbt\big),$$
which is clearly still a horizontal path since the paths 
$$\big(\AD_{\bar{\sigma}_\bbt}\rond f_\bbt \big)(g)=\bar{\sigma}_\bbt\cdot\bar{m}^{f(g)}_\bbt\cdot(\bar{\sigma}_\bbt)\inv$$ are horizontal in $\kernel$, for all $g\in G$.
This concludes the proof.
\end{proof}

\subsection{Connection 1-forms on the band}

A connection on a principal bundle over a Lie groupoid can be equivalently described by a 1-form, called the \textbf{connection 1-form}. 

\begin{defn}[\cite{LTX}] \label{def:LTX}
Let $P\xrightarrow{J}\Gamma_0$ be a principal $G$-bundle over a Lie 
groupoid  $\Gamma_1\rightrightarrows\Gamma_0$. A connection 1-form is
 a usual connection 1-form $\theta\in\Omega^1(P)\otimes\mfg$ of the   principal $G$-bundle
 $P\to \Gamma_0$ (ignoring the groupoid action),
% i.e.
%\begin{gather*}
%  \theta(\xi^{\btr})=\xi,\quad\forall\xi\in\mfg, \\ R_g^*\theta=\Ad_{g\inv }\theta ,\end{gather*}
 satisfying  the additional  equation
 ${\bt}^*\theta- {\bs}^*\theta=0$.
Here
% $R_g$ denotes the right $G$-action on the principal bundle $P\to\Gamma_0$, 
%and $\mfg\to\mathfrak{X}(P):\xi\mapsto\xi^{\btr}$ is the corresponding 
%infinitesimal action,
 ${\bs}$ and ${\bt}$ are the source and target maps of the transformation groupoid $\Gamma_1\times_{{\bt},\Gamma_0,J}P\toto P$ associated to the $\Gamma\simplicial$-action  on $P$.
\end{defn}

Just like the
 usual principal $G$-bundles, we have the following (see Proposition 3.6 in \cite{LTX}):

\begin{prop} \label{prop:4.26} 
For a principal $G$-bundle over a Lie groupoid, a connection is equivalent to a connection 1-form.
\end{prop}

The purpose of this subsection is to construct the connection 1-form
for the induced connection on the band, and to prove directly
 that it satisfies the conditions in Definition \ref{def:LTX}. Hence, this subsection can be considered
 as an alternative approach to obtain, from a connection on a groupoid $G$-extension, an induced connection on its band.

Let $X_1\xrightarrow{\phi}Y_1\toto M$ be a 
Lie groupoid $G$-extension endowed with a groupoid $G$-extension connection. Let $\alpha\in\Omega^1(X_1,\bt^*\liekernel)$ be its corresponding right connection 1-form and $\alpha^{\kernel}\in\Omega^1(\kernel,\bt^*\liekernel)$ the induced connection 1-form on the group bundle $\kernel\xrightarrow{\pi}M$ obtained by restricting $\alpha $ to $\kernel$.

Denote by $\partial$ the Lie group cohomology
 differential
 of $G$ with values in its Lie algebra $\mfg$, where $G$ acts on $\mfg$ by the adjoint action.
 In particular, for any $\xi\in\mfg$, $\partial\xi$ is the $\mfg$-valued function 
on $G$ given by $(\partial\xi)(g)=\xi-Ad_{g\inv }\xi, \ \forall g\in G$.

The following lemma is standard.

\begin{lem}\label{lem:autG}
Let $G$ be a Lie group with Lie algebra $\mfg$. 
\begin{enumerate}
\item The Lie algebra $\Lie\Aut(G)$ of $\Aut(G)$ is naturally identified with the space of 1-cocycles $Z^1(G,\mfg)$ with the bracket:
$$ [z_1,z_2](g) = z_{2*} (z_1(g)) - z_{1*} (z_2(g))-[z_1 (g), z_2(g)],   \quad \forall 
z_1,z_2 \in Z^1(G,\mfg), \ \ g\in G,$$
where for any $z \in Z^1 (G,\mfg)$, $z_*:\mfg\to\mfg$ is the differential of $z$ at the identity.

The isomorphism from $\Lie\Aut(G)$ to $Z^1(G,\mfg)$ is given as follows.
Let $f_t$ be any $C^1$-path in $\Aut(G)$ with $f_0=\idn$. Then the 
element $z(g)$ in $Z^1(G,\mfg)$ corresponding to
$\tfrac{df_t}{dt}|_{t=0}\in \Lie\Aut(G)$ is
 \begin{equation} \label{eq:derivative} z(g)=(L_{g\inv })_* \left.\frac{d f_t(g)}{dt}\right|_{t=0} .\end{equation}

\item The Lie algebra $\Lie\Inn(G)$ of $\Inn(G)$ is naturally identified with the space of 1-coboundaries $B^1(G,\mfg)$.
\item Let $\Ad:\mfg\to B^1(G,\mfg)$ be the Lie algebra morphism
 corresponding to the Lie group morphism $G\to\Inn(G)$ given by $x\to\AD_x$.
Then $\Ad\xi=\partial\xi, \;\forall\xi\in\mfg$.
\item The Lie algebra $\Lie\Out(G)$ of $\Out(G)$ is naturally identified with the
first cohomology group $H^1(G,\mfg)$.
\end{enumerate}
\end{lem}

Define a $C^\infty(G,\mfg)$-valued 1-form  $\alphaiso$
on $\Iso(\kernel,G)$ by
$$(v\ip\alphaiso)(g)=f_*\inv (\ev_{g*}(v)\ip \alphakernel),
 \quad f\in\Iso(\kernel,G), \text{ and } v\in T_f\Iso(\kernel,G),$$ 
where $f_*:\mfg\to\liekernel_m$ is the Lie algebra isomorphism corresponding
to the Lie group isomorphism $f: G\to\kernel_m$,
and $\alphakernel \in \Omega^1 (\kernel, \liekernel)$  is the  pull back
of the right connection 1-form $\alpha$ on $\kernel$.

\begin{lem}
The 1-form $\alphaiso$ is $Z^1(G,\mfg)$-valued. I.e. $\alphaiso\in\Omega^1(\Iso(\kernel,G))\otimes Z^1(G,\mfg)$.
\end{lem}

\begin{proof}
Denote by $m:\kernel\times_M\kernel\to\kernel$ the multiplication in the Lie group bundle $\kernel$ and by $p_1$ and $p_2$ the projections of $\kernel\times_M\kernel$ onto its first and second factors respectively. The relation $\ptr \alpha=0$ 
implies that 
$\ptr _{\kernel}\alphakernel=0$, where $\ptr _{\kernel}$ stands for the restriction of $\ptr $ to the simplicial manifold associated to the groupoid $\kernel\toto M$. 
Therefore, for any $(k_1,k_2)\in\kernel\times_M\kernel$,
\begin{equation}\label{eq:gpalpha} m^*\alphakernel=p_2^*\alphakernel+\Ad_{k_2\inv }p_1^*\alphakernel .\end{equation}
For any $g_1,g_2\in G$ and $f\in\Iso(\kernel,G)$, differentiating with respect to $f$ the relation $\ev_{g_1 g_2}f=m(\ev_{g_1}f,\ev_{g_2}f)$, one obtains that for any $v\in T_f\Iso(\kernel,G)$, $$\ev_{g_1 g_2*}(v)=m_*(\ev_{g_1*}v,\ev_{g_2*}v) .$$
Eq.~\eqref{eq:gpalpha} implies that 
$$\ev_{g_1 g_2*}(v)\ip \alphakernel=\ev_{g_2*}(v)\ip \alphakernel+\Ad_{(f(g_2))\inv}
(\ev_{g_1*}(v)\ip \alphakernel) .$$
Applying $f_*\inv $ yields that
$$(v\ip\alphaiso)(g_1 g_2)=(v\ip\alphaiso)(g_2)+\Ad_{g_2\inv}(v\ip\alphaiso)(g_1).$$
Therefore $\alphaiso$ takes its values in the Lie algebra $Z^1(G,\mfg)$.
\end{proof}

\begin{prop}
\label{pro:bandconnection}
The 1-form $\alphaiso\in\Omega^1(\Iso(\kernel,G))\otimes Z^1(G,\mfg)$
defines a connection on the $\Aut(G)$-principal bundle $\Iso(\kernel,G)\to M$.
\end{prop}

\begin{proof}
For any $\eta\in Z^1(G,\mfg)$, denote by
$\eta^{\btr}_f$ the tangent vector in $T_f\Iso(\kernel,G)$
induced by the infinitesimal action of
the Lie algebra $Z^1(G,\mfg)$.

Differentiating the
relation $ \ev_g (f \cdot \psi) =
 (f \rond \psi) (g)$ with respect to $\psi \in \Aut(G)$ at the identity,
one obtains the
relation $\ev_{g*}(\eta^{\btr}_f)=\big(f_*\eta(g)\big)^{\btr}$ 
for any $\eta\in Z^1(G,\mfg)$ and $g\in G$. Hence $\alphakernel(\ev_{g*}
\eta^{\btr})= f_*\eta(g)$. Applying $f_*\inv $, we obtain, for any $\eta
\in Z^1(G,\mfg)$,
\begin{equation} \label{eq:co1} \eta^{\btr}\ip\alphaiso=\eta .\end{equation}

Let us now check that $\alphaiso$ 
is $\Aut(G)$-equivariant. Fix $\psi\in\Aut(G)$. Differentiating the relation $\ev_g(f\cdot\psi)=(f\smalcirc\psi)(g)$ with respect to $f$, one obtains the relation $\ev_{g*}(R_\psi)_* v =\ev_{\psi(g)*}v$ for all $v \in T_f \Iso (\kernel, G)$. Applying $(f\smalcirc\psi)_*\inv \smalcirc\alphakernel$ yields 

\begin{equation}\label{eq:co2}
 (\psi_* v\ip\alphaiso )(g)=\psi_*\inv  (v\ip \alphaiso \ \big( \psi(g) \big) = 
\big(\Ad_{\psi\inv }
(v\ip\alphaiso)\big)(g).
\end{equation}

Eqs.~\eqref{eq:co1}-\eqref{eq:co2} imply that $\alphaiso$ is a connection 1-form.
\end{proof}

\begin{prop} \label{prop:aut}
The following are equivalent characterizations of a connection
on the $\Aut(G)$-principal bundle $\Iso(\kernel,G)\to M$:
\begin{enumerate}
\item \label{w1} a connection 1-form $\alphaiso\in\Omega^1
\big(\Iso(\kernel,G)\big)\otimes Z^1(G,\mfg)$ as in Proposition \ref{pro:bandconnection};
\item \label{w2} a distribution $\hiso\subset T\Iso(\kernel,G)$ as in Eq.~\eqref{eq:dis};
\item \label{w3}  horizontal paths in $\Iso(\kernel,G)$ are those paths $f_t$ such that, 
for any $g\in G$, $\ev_g f_t=f_t(g)$ is a horizontal path in $\kernel$.
\end{enumerate}
\end{prop}

\begin{proof}
\fbox{\ref{w1}$\Leftrightarrow$\ref{w2}} 
It is straightforward to check 
that $\ker(\alphaiso)=\hiso$. Hence $\hiso$ defines the same connection as
 $\alphaiso$. In particular, it is an $\Aut(G)$-invariant horizontal distribution.

\fbox{\ref{w2}$\Leftrightarrow$\ref{w3}} 
A path $f_t $ is horizontal if, and only if, its tangent
 vectors are in $\hiso$. In other words, the tangent vectors of the paths $f_t(g)$ are
 in $\hkernel$ for all $g\in G$. Therefore the horizontal paths of the connection defined 
by $\hiso$ or $\alphaiso$ are the paths in $\Iso(\kernel, G)$ such that, for any $g\in G$, $\ev_g(f_t)=f_t(g)$ is a horizontal path in $\kernel$. This completes the proof.
\end{proof}

Assume that we are given a short exact sequence of Lie groups $1\to
 R\to G\to H\to 1$. Then, for any principal $G$-bundle $P\to M$,
 $P/R\to M$ is a
 principal $H$-bundle. Given a connection 1-form
 $\theta\in\Omega^1(P)\otimes\mfg$, then 
$\bar{\theta}:=\pr( \theta)$
 is a $\mathfrak{h}$-valued
1-form on $P$ which is $R$-basic. Here
$\pr: \mathfrak{g}\to \mathfrak{h}$ is the natural projection.
 Therefore, it descends to a $\mathfrak{h}$-valued 1-form on $P/R$, which is a connection 1-form on the principal $H$-bundle $P/R\to M$.

Applying this construction to the particular case of the exact sequence $1\to\Inn(G)/Z(G)\to\Aut(G)\to\Out(G)\to 1$ and the connection 1-form $\alphaiso$ on $\Iso(\kernel,G)$,
 we obtain a connection 1-form $\alphaout\in\Omega^1(\Out(\kernel,G))\otimes H^1(G,\mfg)$. 

Recall that $X_1\toto M$ acts on $\Iso(\kernel,G)\xrightarrow{\pi}M$ by conjugation, so $\Iso(\kernel,G)\xrightarrow{\pi}M$ is indeed an $\Aut(G)$-torsor over $X_1\toto M$. One can consider the transformation groupoid: $$X_1\times_M \Iso(\kernel,G)\toto\Iso(\kernel,G).$$
Similarly, $Y_1\toto M$ acts on $\Out(\kernel,G)\xrightarrow{\pi}M$, so $\Out(\kernel,G)\to M$ is an $\Out(G)$-torsor over $Y_1\toto M$. And one can consider the transformation groupoid $$Y_1\times_M \Out(\kernel,G)\toto\Out(\kernel,G).$$

\begin{prop}
\label{prop:coband2}
The connection 1-form $\alphaout\in\Omega^1(\Out(\kernel,G))\otimes H^1(G,\mfg)$
 is compatible with the $Y\simplicial$-action.
More precisely, the following relation holds: $$\bs^*\alphaout-\bt^*\alphaout=0.$$
Here $\bs$ and $\bt$ denote the source and target maps of the transformation groupoid 
$Y_1\times_M\Out(\kernel,G)\toto\Out(\kernel,G)$.
\end{prop}

The above proposition is an immediate consequence of the following lemma.

\begin{lem} \label{pro:st}
For all $(x,f)\in X_1\times_M\Iso(\kernel,G)$, we have
\begin{equation} \label{eq:delalpha}
\bs^*\alphaiso-\bt^*\alphaiso=\partial\big(p^*(f_*\inv \rond\alpha_x)\big), 
\end{equation}
\end{lem}

where both sides are  1-forms on $X_1\times_M\Iso(\kernel,G)$ with values in $Z^1(G,\mfg)$.
Here, $\bs$ and $\bt$ are the source and target maps of the transformation groupoid $X_1\times_M\Iso(\kernel,G)\toto\Iso(\kernel,G)$ and $p:X_1\times_M\Iso(\kernel,G)\to X_1$ is the projection on the first component.
Note that $\bs^*\alphaiso$ and $\bt^*\alphaiso$ are $1$-forms
 on $X_1\times_M\Iso(\kernel,G)$  with values in $Z^1(G,\mfg)$.
On the right hand side, the pull-back via $p$ of the composition of the 
$\liekernel_{\bt(x)}$-valued covector $\alpha_x$ on $X_1$ and the isomorphism 
$f_*\inv :\liekernel_{{\bt}(x)}\to\mfg$ is a covector on $X_1\times_M\Iso(\kernel,G)$ at $(x, f)$ with values in $\mfg$. Therefore, 
$\partial\big(p^*(f_*\inv \rond\alpha_x)\big)$ is a covector of $X_1\times_M\Iso(\kernel,G)$ at $(x, f)$ with values in $B^1(G,\mfg)$.

\begin{proof}
The tangent space
 $T_{(x,f)}\big(X_1\times_M\Iso(\kernel,G)\big)$ consists of the 
pairs $(u,v)\in T_x X_1\times T_f\Iso(\kernel,G)$ such that
 $\bt_*u=\pi_*v$, and therefore is the direct sum of the following 
three vector spaces:
\begin{align*}
E_1 &:= \genrel{(0,\eta^{\btr}_f)}{\eta\in Z^1(G,\mfg)}, \\
E_2 &:= \genrel{(\xi^{\btr}_x,0)}{\xi\in\liekernel_{{\bt}(x)}}, \\
E_3 &:= \genrel{(u,v)}{u\in H_x,v\in\hiso_f,\bt_*u=\pi_*v}. 
\end{align*}

Below we check that Eq.~\eqref{eq:delalpha} holds on each of the direct summands.

\begin{itemize}
\item By construction, $$(0,\eta^{\btr})\ip\big(\bs^*\alphaiso-\bt^*\alphaiso\big)=\eta-\eta=0,$$
and $$(0,\eta^{\btr})\ip p^*\alpha=0 .$$
Therefore, Eq.~\eqref{eq:delalpha} holds on $E_1$.
\item Fixing $(x,f) \in X_1\times_M \Iso (\kernel, G)$ and differentiating the relation 
$${\bs}(x\cdot k,f)=\AD_{xk} \smalcirc f=(\AD_x  \smalcirc f)\cdot \AD_{f\inv (k)}\in\Iso(\kernel,G)$$ with respect to $k\in\kernel_{{\bt}(x)}$ at the identity,
we obtain 
$$ {\bs}_*(\xi^{\btr},0)|_{(x,f)}=
\big(\partial(f_*\inv \xi)\big)^{\btr}_{\AD_x  \smalcirc f}.$$
Here $\AD_{f\inv (k)}$ is considered as an element in $\Inn(G)\subset\Aut(G)$ 
and the dot refers to the right $\Aut(G)$-action on $\Iso(\kernel,G)$.
Hence, we have
$$(\xi^{\btr},0)\ip\big( {\bs}^*\alphaiso-\bt^*\alphaiso\big)=\partial(f_*\inv \xi),$$
while, on the other hand,
$$(\xi^{\btr},0)\ip\partial\big(p^* f_*\inv (\alpha_x)\big)=\partial(f_*\inv \xi).$$
Therefore, Eq.~\eqref{eq:delalpha} holds on $E_2$.
\item Let $(u,v) \in E_3$. We have
\begin{equation} \label{eq:cas3.1} 
(u,v)\ip\partial\big(p^*(f_*\inv \rond\alpha)\big)=\partial\big(f_*\inv (v\ip\alpha_x)\big)=0 
\end{equation}
and
\begin{equation}\label{eq:cas3.2} 
(u,v)\ip\bt^*\alphaiso=v\ip\alphaiso=0.
\end{equation}
Let $x_t$ and $f_t$ be horizontal paths in $X_1$ and $\Iso (\kernel, G)$, respectively, such that $\pi\smalcirc f_t=\bt\smalcirc x_t$, for all $t$, and $\left.\tfrac{d}{dt}(x_t,f_t)\right|_{t=0}=(u,v)$. For any $g\in G$, $f_t(g)$ is an horizontal path in $\kernel$ by Proposition \ref{prop:aut}(3). Since a path in $\kernel$ is horizontal in $\kernel$ if, and only if, it is horizontal when considered as a path in $X_1$, and the product of horizontal paths in $X_1$ is horizontal by Proposition \ref{pro:para}, $x_t\cdot f_t(g)\cdot x_t\inv $ is an horizontal path in $\kernel$ for all $g\in G$.
Since $\bs(x,f)=\AD_x\rond f$, it thus follows that $ {\bs}_*(u,v)=\left.\tfrac{d}{dt}(x_t\cdot f_t(g)\cdot x_t\inv )\right|_{t=0}$ is a horizontal vector in $\Iso(\kernel,G)$.
Together with Eqs.~\eqref{eq:cas3.1}-\eqref{eq:cas3.2}, this implies that Eq.~\eqref{eq:delalpha} holds on $E_3$.
\end{itemize}
\end{proof}

Proposition \ref{prop:coband2} means that $\alphaout$ is indeed a connection 1-form for the principal $\Out(G)$-bundle $\Out(\kernel,G)\to M$ over $Y\simplicial$, i.e. a connection 1-form on the band.

In summary, we have the following

\begin{thm}
The 1-form $\alphaout\in\Omega^1\big(\Out(\kernel,G)\big)\otimes H^1(G,\mfg)$ is a connection 1-form on the band of the $G$-extension, whose corresponding horizontal distribution is $\hout$.
\end{thm}

\begin{rmk} 
As shown in \cite{LTX}, connections behave well with respect to Morita equivalence. More precisely, there is a 1-1 correspondence between connections on principal bundles over Morita equivalent groupoids. Therefore, a connection in our sense indeed yields a connection on the corresponding torsor over the stack. This implies that a connection on a gerbe induces a connection on its band.
\end{rmk}

\subsection{Curvature on the band}

This section is devoted to describing
the relation between
the Ehresmann curvature $\omegakernel$
on the group bundle $\kernel\to M$ and the curvature $\omegaiso$ 
on the $\Aut(G)$-principal bundle
 $\Iso(\kernel,G)\xrightarrow{\pi}M$, as well as  
the curvature $\omegaout$ on the band $\Out(\kernel,G)\xrightarrow{\pi}M$.

Since $\omegakernel$ 
is an horizontal form, it can be  considered as a form on the manifold $M$.
More precisely, one defines 
$\omegakernel \in \Omega^2(M,Z^1(\kernel, \liekernel)) $
 as follows. For any $k \in \kernel_{m}$ and    $u,v \in T_m M$, 
 
\begin{equation}\label{eq:bothomega} 
 (i_{ u \wedge v } \omegakernel) (k) =
  i_{\tilde{u}\wedge \tilde{v}}  (\omega |_{k} ),
    \end{equation}
where $\tilde{u},\tilde{v} \in T_k \kernel$ are any tangent vectors 
such that $\phi_* (\tilde{u}) =u  $ and $\phi_* (\tilde{v})=v $.
This definition requires some justification.

\begin{itemize}
\item  First, it is clear that
 $(i_{ u \wedge v } \omegakernel) (k)$  is independent of the 
choice of $\tilde{u},\tilde{v}$
because $ \omega|_{\kernel}$ is an horizontal $2$-form, and is 
therefore well defined.
\item Second, we need  to check that, for any
fixed $u, \ v$, the map
$ k \to (i_{ u \wedge v } \omegakernel) (k)$ is  
an element in $Z^1(\kernel, {\mathfrak K})$.
For any $k_1,k_2 \in \kernel_{m}$, let 
$ \tilde{u}_1, \tilde{v}_1 \in T_{k_1}\kernel$
and $  \tilde{u}_2, \tilde{v}_2 \in T_{k_2} \kernel$ be 
any tangent vectors such that
$\phi_* (\tilde{u}_1)= \phi_* (\tilde{u}_2)= u$
 and $\phi_* (\tilde{v}_1)= \phi_* (\tilde{v}_2)= v $.
Then  $\tilde{u}_1 \cdot \tilde{u}_2$ and $\tilde{v}_1 \cdot \tilde{v}_2$ 
(where the dot stands for the product in the tangent groupoid 
 $TX_1 \toto TM$) are elements in $T_{k_1 k_2}\kernel$ such that
 $\phi_*( \tilde{u}_1 \cdot \tilde{u}_2)= u $
and $ \phi_* ( \tilde{v}_1 \cdot \tilde{v}_2) = v$. 
The Bianchi identity Eq.~\eqref{eq:Bianchi2} implies that
  $ \omega(  \tilde{u}_1 \cdot \tilde{u}_2, \tilde{v}_1 \cdot \tilde{v}_2)
= \omega (\tilde{u}_2, \tilde{v}_2)+
\Ad_{k_2\inv } \omega(\tilde{u}_1, \tilde{v}_1)$. Hence
  $$  (i_{u \wedge v}  \omegakernel) ( k_1  k_2)=
(i_{u \wedge v}\omegakernel )(  k_2) +  \Ad_{k_2\inv } (i_{u \wedge v} \omegakernel)
 ( k_1)  .$$
\end{itemize}

Note that any $f \in \Iso(\kernel,G)|_m$  induces  an
isomorphism   $T_f :   Z^1(\kernel_m,{\mathfrak K}_m)  
\to  Z^1(G,\mfg)$ given  by,
 $\forall z \in Z^1(\kernel_m,{\mathfrak K}_m)$,
 \begin{equation}
\label{eq:Tf}
 T_f (z) (g) =  f_*\inv  z\big( f(g) \big),  \ \ \ \ \forall g \in G .
\end{equation}
Tensoring with $\wedge^2 T^*_mM$, 
this extends naturally to a map (denoted
 again by the same symbol by abuse of notation)
$$ T_f: \wedge^2 T_m^* M \otimes  Z^1(\kernel_m,{\mathfrak K}_m)
\lto   \wedge^2 T^*_{m}M \otimes Z^1(G,\mfg). $$

\begin{prop} \label{prop:autker}
For any $f\in \Iso(\kernel,G)|_m$, we have the following relation
$$ \omegaiso |_f = \pi^* T_f  \omegakernel|_m,    $$
where 
$ \omegaiso |_f \in \wedge^2  T_f^* 
\Iso (\kernel,G) \otimes Z^1(G,\mfg)$ is the curvature
of the connection $ \alphaiso$,
$\omegakernel|_m \in \wedge^2 T^*_{m}M \otimes Z^1(\kernel,{\mathfrak K})$ is as in Eq.~\eqref{eq:bothomega} and
$\pi:  \Iso(\kernel,G)\to M$ is the projection.
\end{prop}

\begin{proof}
For any path $f_t : [0,1] \to \Iso (\kernel,G)$,
let $\overline{f_t} $ be the unique horizontal path
in $\Iso (\kernel,G)$ starting at the same point $f_0$
and satisfying $\pi \smalcirc \gamma =\pi \smalcirc\bar{\gamma}$.
Then the holonomy $\Hol(f_t)\in \Aut(G)$ of the path $f_t$ is given by
 \begin{equation} \label{eq:hogamma} 
  f_1\cdot \Hol(f_t) = \overline{f}_1  .
\end{equation} 
The curvature  $\omegaiso $ can be expressed by
\begin{equation} \label{eq:hogamma2}
\omegaiso(u,v)=\left.\tfrac{\partial^2}{\partial\epsilon_1\partial\epsilon_2}\right|_{\epsilon_1=\epsilon_2=0} \Hol(L^C_{\epsilon_1,\epsilon_2})
\end{equation}
for any smooth map $C$ from $\R^2$ to
$\Iso (\kernel,G)$ adapted to the tangent vectors $u,v \in T_f \Iso(\kernel,G)$ 
(see Eq.~\eqref{eq:deflc} for the notation  $L^C_{\epsilon_1,\epsilon_2}$).
 According to Lemma \ref{lem:autG}(1), Eq.~\eqref{eq:hogamma2} can be written 
as
\begin{equation} \label{eq:hogamma2withg}
\omegaiso(u,v)(g)=\left.\tfrac{\partial^2}{\partial\epsilon_1\partial\epsilon_2}\right|_{\epsilon_1=\epsilon_2=0}
 g\inv  \Hol(L^C_{\epsilon_1,\epsilon_2})(g), 
\qquad \forall g \in G \end{equation}

By Lemma \ref{lem:5.2}, for any path $f_t$ in $\Iso (\kernel,G)$ and any $g \in G$,
$\ev_g(\overline{f_t})$ is an horizontal path 
starting at $\ev_g(f_0)$.
Hence $\ev_g\overline{f_t}=\overline{\ev_g f_t}$
holds for all $t \in [0,1]$. 
In particular, for $t=1$, one obtains
$$\ev_g \big(f_1\cdot \Hol(f_t) \big) = (\ev_g f_1) \cdot \Hol(\ev_g f_t), \qquad \forall g \in G, $$
and therefore
\begin{equation}\label{eq:relholo} 
 f_1 (g\inv ) \cdot  [f_1\cdot \Hol(f_t)] (g) =\Hol (\ev_g f_t)  .
\end{equation}
The latter can be re-written as
\begin{equation}
\label{eq:relholo1} 
 f_1 (g\inv  \Hol(f_t) (g)) = \Hol (\ev_g f_t)  .
\end{equation}

For any smooth map $C$ from a neighborhood $\mathcal{U} \subset  \R^2 $
of $0$  to a neighborhood of  $f \in \Iso (\kernel , G)$
 adapted to  the tangent vectors $u,v \in T_f \Iso (\kernel,G)$, the smooth map
$\ev_g \smalcirc C \in C^\infty(\mathcal{U},\kernel)$
is adapted to $   \ev_{g*} u , \ev_{g*} v \in T_{f(g)} \kernel $.
Take $f_t$ to be the loop $L^C_{\epsilon_1,\epsilon_2}$. Then
$\ev_g f_t$ is the loop $L^{\ev_g \smalcirc C}_{\epsilon_1,\epsilon_2}$. 
Now, according to Eq.~\eqref{eq:relholo1},  one obtains
 that, for any $\epsilon_1, \epsilon_2$, 
$$f\big(g\inv \Hol(L^C_{\epsilon_1,\epsilon_2})(g)\big) 
=\Hol(L^{\ev_g \smalcirc C}_{\epsilon_1,\epsilon_2}), \qquad\forall g\in G.$$
Applying $\left.\tfrac{\partial^2}{\partial\epsilon_1\partial\epsilon_2}\right|_{\epsilon_1=\epsilon_2=0}$ to 
both sides and using Proposition \ref{pro:4.25} and Eq.~\eqref{eq:hogamma2withg}, we have
$$f_*(\omegaiso(u,v)(g))=\omega(\ev_{g*}u,\ev_{g*}v)|_{f(g)}, \qquad\forall g\in G.$$
By  Eq.~\eqref{eq:bothomega}, we obtain 
$$f_*(\omegaiso(u,v)(g))=\omegakernel(\pi_* u,\pi_* v)(f(g)), \qquad\forall g\in G.$$
This concludes the proof.
\end{proof}

Denote by $[\omegakernel] \in \Omega^2\big(M,H^1(\kernel,\liekernel)\big)$
 the class of $\omegakernel$.
Since for any  $f \in \Iso(\kernel,G)|_m$, the map $T_f$ defined in Eq.~\eqref{eq:Tf} 
 maps $B^1(\kernel_m,\liekernel_m)$ to $B^1(G,\mfg)$,
it induces  a map $T_f: H^1(\kernel_m,\liekernel_m)\to H^1(G,\mfg)$.
It is simple to check that $T_f$ only depends on the class of
$\overline{f}\in\Out(\kernel,G)|_m$. Therefore, for any $\overline{f}\in\Out(\kernel,G)|_m$,
 we have a well defined map
$ T_{\overline{f}}: H^1(\kernel_m,\liekernel_m)\to H^1(G,\mfg)$.
Alternatively, one may also obtain $ T_{\overline{f}}$ as follows.
For any $\overline{f}\in\Out(\kernel,G)|_m$, there is an induced
group homomorphism $\Out(\kernel)|_m \to\Out(G)$ 
by taking the conjugation. Then $T_{\overline{f}}:H^1(\kernel_m,\liekernel_m)\to
 H^1(G,\mfg)$
is its corresponding Lie algebra homomorphism.

Tensoring  with $\wedge^2 T^*_mM$ as before,  $T_{\overline{f}}$
 extends naturally to a map
 $$ T_{\overline{f}}: \wedge^2 T_m^* M \otimes  H^1(\kernel_m,\liekernel_m) 
\to  \wedge^2 T_m^* M \otimes H^1(G,\mfg) .$$

Proposition \ref{prop:autker} immediately implies the following

\begin{thm} 
\label{prop:outker}
For any $\overline{f}\in \Out(\kernel,G)|_m $,
the curvature on the band $\omegaout|_{\overline{f}} \in \wedge^2 T^*_{m}\Out(\kernel,G)
\otimes H^1(G,\mfg)$ and the class 
$[\omegakernel] \in \wedge^2 T^*_{m}M \otimes H^1(\kernel,\liekernel)$ 
are related by  the following equation
$$\omegaout|_{\overline{f}} = \pi^* T_{\overline{f}}  [\omegakernel], $$
where $\pi :\Out(\kernel,G)\to M$ is the bundle projection.
\end{thm}

As an immediate consequence, we have the following

\begin{cor}
\label{cor:4.46}
The band of a $G$-extension $X\simplicial\xrightarrow{\phi}Y\simplicial$ is flat if and only if $[\omegakernel]=0$.
\end{cor}

To end this section, we describe an important relation
 between the curvatures  $\omega \in \Omega^2 (X_1, {\bt}^* {\mathfrak K} )$
and $\omegakernel$.

Let $\omega^{\liekernel}\in\Omega^2\big(M,\Lie\Aut(\liekernel)\big)$
be the curvature of the induced connection $\nabla$ on the
Lie algebra bundle $\liekernel\to M$.
Note that $\Lie\Aut(\liekernel)$ can be naturally identified
with the space of derivations of $\liekernel$. 
They can also be viewed
as Lie algebra $1$-cocycles with values in the Lie algebra
itself with the action being the adjoint action.   
Therefore we have an identification
$\Lie\Aut(\liekernel)\isomorphism Z^1(\liekernel,\liekernel)$.

First we need to give a relation between 
$\omega^{\liekernel}$ and the curvature $\omegakernel$
on the group bundle $\kernel\to M$.
For any Lie group $G$ with Lie algebra $\mfg$,
differentiating a Lie group cocycle in
$Z^1(G,\mfg)$ at the identity, one obtains an element 
in $Z^1(\mfg,\mfg)$. 
More generally, for any group bundle $\kernel\to M$, 
this allows us to obtain a map
$Z^1(\kernel,\liekernel)\to Z^1(\liekernel,\liekernel)$,
which extends to a map
$D:\Omega^2\big(M,Z^1(\kernel,\liekernel)\big)\to
\Omega^2\big(M,Z^1(\liekernel,\liekernel)\big)$.
The following can be easily verified:

\begin{lem}
\begin{equation}\label{eq:curvLieK}
D\omegakernel=\omega^{\liekernel}.
\end{equation}
\end{lem}

Recall that for any Lie group $G$ with Lie algebra $\mfg$,
we have 
$\partial:\mfg\to B^1(G,\mfg):\xi\mapsto\big(g\mapsto\xi-\Ad_{g\inv }\xi\big)$.
Applying $\partial$ fiberwise, one obtains a map
$\Gamma(\bt^*\liekernel)\to\Gamma({\bt}^* B^1(\kernel,\liekernel))\subset\Gamma({\bt}^* Z^1(\kernel,\liekernel))$ 
that we denote by $\partial_{\kernel}$.
 By abuse of notation, we denote
by $ \partial_{\kernel} $ again the induced map 
$$\partial_{\kernel}:\Omega^2(X_1,{\bt}^* \liekernel) \to \Omega^2(X_1,{\bt}^* B^1(\kernel,\liekernel))
\subset \Omega^2(X_1,{\bt}^* Z^1(\kernel,\liekernel)) .$$ 
We can now prove the following:

\begin{prop}\label{prop:curvuptocenter}
The following relations hold, for any $x \in X_1 $,
\begin{equation}\label{eq:delalpha2}
\begin{gathered} 
\Ad_{x\inv }{{\bs}}^*{\omegakernel}|_{{\bs}(x)}-{\bt}^*{\omegakernel}|_{{\bt}(x)}=-\partial_{\kernel}\omega|_x \\
\Ad_{x\inv }{{\bs}}^*\omega^{\liekernel}|_{{\bs}(x)}-{\bt}^*\omega^{\liekernel}|_{{\bt}(x)}=\ad_{\omega_x} ,
\end{gathered}
\end{equation}
where $\Ad_{x\inv }$ stands for the natural isomorphism
from $Z^1(\kernel_{\bs(x)},\liekernel_{\bs(x)})$
to $Z^1(\kernel_{\bt(x)},\liekernel_{\bt(x)})$, 
and from $Z^1(\liekernel_{\bs(x)},\liekernel_{\bs(x)})$ 
to $Z^1(\liekernel_{\bt(x)},\liekernel_{\bt(x)})$
induced by $\AD_x:\kernel_{\bt(x)}\to\kernel_{\bs(x)}$.
Note that $\omega|_x \in \wedge^2 T_x^* X_1 \otimes \liekernel_{\bt(x)}$, 
$\partial \omega|_x \in \wedge^2 T_x^* X_1\otimes Z^1 (\kernel_{\bt(x)},\liekernel_{\bt(x)})$, 
$\omegakernel|_{\bt(x)} \in \wedge^2 T_{\bt(x)}^* M \otimes 
Z^1(\kernel_{\bt(x)},\liekernel_{\bt(x)})$, and $\omegakernel|_{\bs(x)}   
\in \wedge^2 T_{\bs(x)}^* M \otimes Z^1(\kernel_{\bs(x)},\liekernel_{\bs(x)})$.
Similarly, $\ad_{\omega_x}\in \wedge^2 T_x^* X_1 \otimes Z^1(\liekernel_{\bt(x)},\liekernel_{\bt(x)})$.
\end{prop}

\begin{proof}
The second identity in Eq.~\eqref{eq:delalpha2} follows from the
first one by applying $D$ and using Eq.~\eqref{eq:curvLieK}.
Let us now prove the first one. This equation is equivalent to
\begin{equation}
\label{eq:delalphabis}   
\Ad_{x\inv }{\bs}^*\omegakernel(\AD_x k)-\bt^*\omegakernel(k) = -\omega_x+ \Ad_{k\inv }\omega_x \in
\wedge^2 T_{x}^* X_1\otimes\liekernel_{\bt(x)}.
\end{equation}

Let $u_1,u_2 \in T_x X_1$ be
 any two tangent vectors at the point $x$, and
$\epsilon_1,\epsilon_2 \in T_k \kernel$ any 
two tangent vectors of $\kernel$ at the point $k$ 
such that $\bt_* u_i = \pi_* \epsilon_i$, $i\in\{1,2\}$,
(i.e. $(u_i,\epsilon_i)$ is a composable pair
in the tangent groupoid $TX_1\toto TM$ for $i\in\{1,2\}$).

By using Eq.~\eqref{eq:Bianchi2} repeatedly, we obtain
\[ \omega_{\AD_x k}(u_1 \epsilon_1 u_1\inv , u_2 \epsilon_2 u_2\inv )
= \Ad_{x k\inv }\omega_x(u_1,u_2)+\Ad_{x}\omega_{k}(\epsilon_1,\epsilon_2)-\Ad_{x}\omega_{x}(u_1,u_2) .\]
Hence
\begin{equation} \label{eqdelalp1} 
\Ad_{x\inv }\omega_{\AD_x k}(u_1 \epsilon_1 u_1\inv ,u_2 \epsilon_2 u_2\inv )
=\Ad_{k\inv }\omega_x(u_1,u_2)-\omega_x(u_1,u_2)+\omega_{k}(\epsilon_1,\epsilon_2).
\end{equation}

Since $\epsilon_1,\epsilon_2 \in T_k \kernel$, 
$u_i \epsilon_i u_i\inv  \in T_{\AD_x k} \kernel$ 
for $i \in \{1,2\}$.
According to Eq.~\eqref{eq:bothomega},
\begin{equation} \label{eqdelalp2}
\omega_k(\epsilon_1,\epsilon_2)
=\omegakernel_k(\phi_*\epsilon_1,\phi_*\epsilon_2)
=\omegakernel_k(\bt_* u_1,\bt_* u_2)
\end{equation}
and
\begin{equation} \label{eqdelalp3}
\omega_{\AD_x k}(u_1 \epsilon_1 u_1\inv ,u_2 \epsilon_2 u_2\inv )
=\omegakernel_{\AD_x k}(\phi_*(u_1 \epsilon_1 u_1\inv ),\phi_*(u_2 \epsilon_2 u_2\inv ))
=\omegakernel_{\AD_x k}(\bs_* u_1,\bs_* u_2) .
\end{equation}

The result now follows from Eqs.~\eqref{eqdelalp1}-\eqref{eqdelalp3}.
\end{proof}

\begin{rmk} \label{rmk:uptocenter}
For the Lie algebra $\mfg$ of any connected Lie group $G$, denote
by $Z(\mfg)$ its center. Since $G$ is connected,
one has $Z(\mfg)=\genrel{\xi\in \mfg }{\Ad_g\xi=\xi, \;\forall g\in G}=\ker(\partial)$.

Similarly, let $\liecenter=\coprod_{m \in M}\liecenter_m$
be the vector bundle over $M$ obtained 
by taking the center of each fiber, and $\bt^* \liecenter \to X_1$, its pull back by the target map $\bt$.
The space of sections of the vector bundle $\bt^* \liecenter \to X_1$ 
is precisely the kernel of $\partial_{\kernel}$.

Proposition \ref{prop:curvuptocenter} implies that the image of $\omega$
under $\partial_{\kernel}$ is entirely determined by $\omegakernel$.
Equivalently, this means that the class of $\omega$ in $\Omega^2(X_1,\bt^* \tfrac{\liekernel}{Z(\liekernel)})$ 
is entirely determined by $\omegakernel$.
\end{rmk}

\section{Cohomology theory of connections}
\label{sec:6}
The purpose of this section is to develop  a  cohomology theory
for groupoid extensions,  which appears naturally
 while studying connections and curvings.

 First of all, in Section \ref{sec:coho1}, 
we introduce the cohomological object in question, which we call 
horizontal cohomology of a Lie groupoid extension. We also discuss in Section
\ref{sec:coho1} how and when the 
horizontal cohomology can be pulled back by a morphism of Lie groupoid extensions.
In Section \ref{sec:coho2}, we show that the horizontal cohomology
 can actually be in fact defined for a gerbe, i.e.
it is invariant under Morita equivalence. 
In Section \ref{sec:coho3},
 we introduce an obstruction class $[\obs]$ in the horizontal cohomology, 
which characterizes the existence of connections.
As a consequence, we show that, if a groupoid extension admits a connection, 
then any Morita equivalent groupoid extension admits a connection as well. 
In the language of stacks and gerbes, it means that the existence of a connection ``goes down'' to a ``gerbe'' notion.
In Section \ref{sec:coho4}, we compute the horizontal cohomology for $G$-gerbes over a manifold, 
i.e. for a $G$-extension of a \v{C}ech groupoid. We show that, in this case, 
the class $[\obs]$ vanishes, and therefore, any $G$-extension 
of a \v{C}ech groupoid (i.e. any $G$-gerbe over a manifold) 
admits a connection. Section \ref{sec:6.5} is devoted to the
study of flat gerbes.   In Section \ref{sec:6.6}, we study 
connections and curvings on central extensions. 

\subsection{Horizontal cohomology} \label{sec:coho1}

Recall that, for any fiber bundle $T\xrightarrow{\phi}B$, 
with $T$ and $B$ being finite dimensional manifolds, an $l$-form 
$\xi$ on $T$ (possibly valued in some vector bundle over $T$) 
is said to be \textbf{horizontal} if $v\ii\xi =0$ for any 
$v\in\ker(\phi_*)$.
In other words, a form is horizontal if it vanishes when contracted with a 
vector tangent to the fiber of $\phi$.

\begin{rmk}
\begin{enumerate}
\item
The reader should not confuse horizontal forms with basic forms. 
Basic forms are simply obtained by pulling-back forms on the base space, 
while horizontal forms form a much bigger space in general.
\item
Given  a fiber bundle $T\xrightarrow{\phi}B$, 
horizontal $l$-forms on $T$ can also be considered  as  $l$-forms 
on the base manifold  $B$ with values in, for any $b\in B$, 
$C^\infty (\phi\inv(b))$.
We will use both viewpoints all along this section.
\end{enumerate}
\end{rmk}

We now introduce the horizontal forms relevant to our situation. 
Given a Lie groupoid extension $X_1\xrightarrow{\phi}Y_1\toto M$
with kernel $\kernel\to M$, for any $n\in\N$,
the manifold $X_n$ of $n$-tuples of composable arrows of $X_1\toto M$ is
a fiber bundle over $Y_n$, with respect to the projection
$$\phi_n:\big(x_1,\cdots,x_n\big)\mapsto\big(\phi(x_1),\cdots,\phi(x_n)\big) .$$
A typical fiber over a point
 $(y_1,\cdots,y_n)\in Y_n$ is isomorphic to 
$\kernel_{\bt(y_1)}\times\kernel_{\bt(y_2)}
\times\cdots\times\kernel_{\bt(y_n)}$.
By a \textbf{horizontal $\liekernel$-valued $l$-form on $X_n$}, 
we mean a $\liekernel$-valued $l$-form on $X_n$ 
(i.e. a section of the vector bundle $\wedge^l T^*X_n \otimes \bt^* \liekernel \to X_n$, 
which is horizontal with respect to the fiber bundle $\phi_n: X_n \to Y_n$).
From now on, for any $n,l\in\N$, we denote by $\Omega_{\hor}^l(X_n,\bt^*\liekernel)$ 
the space of horizontal $\liekernel$-valued $l$-form on $X_n$.

We list below two important examples of horizontal forms.

\begin{prop} \label{ex:hor1form} 
The Ehresmann curvature $\omega$ of a connection on a groupoid extension
$X_1\xrightarrow{\phi}Y_1\toto M$ is an horizontal $2$-form.
\end{prop}

\begin{proof}
This follows immediately from Eq.~\eqref{eq:defehrcurv}.
\end{proof}

\begin{prop}\label{ex:hor2form}
Let $X_1\xrightarrow{\phi}Y_1\toto M$ be a Lie groupoid extension
with kernel $\kernel\to M$, whose corresponding Lie algebra bundle is $\liekernel\to M$.
Let $\theta\in\Omega^1(X_1,\bt^*\liekernel)$
be a $\liekernel$-valued 1-form  on $X_1$ satisfying
          \begin{equation}\label{eq:rightaction}
           \theta(\xi^{\btr})=\xi ,
   \hspace{1cm} \forall \xi \in \liekernel,
          \end{equation}
  where, as in Section \ref{sec:4.3},
 for any $\xi \in \liekernel$ we denote by
 $ \xi^{\btr}$ the fundamental
vector field on $X_1$ corresponding
to  the infinitesimal right action of $\kernel\to M$ on $X_1$.
Then $\ptr  \theta \in \Omega^1 (X_2,
\bt^*\liekernel)$ is an horizontal 1-form on $X_2$.
\end{prop}

\begin{proof}
Let us introduce some notations.
The group bundle $\kernel \to M $ acts on $X_2$ 
from the right in two different ways.
The first one is given, for all $(x_1,x_2) \in X_2$ and 
$k \in \kernel_{\bt(x_1)}$, by $(x_1,x_2)\cdot k = (x_1\cdot k, x_2)$; 
the second one is given, for all $(x_1,x_2) \in X_2$ and 
$k \in \kernel_{\bt(x_2)}$, by $(x_1,x_2)\cdot k = (x_1,x_2\cdot k)$.
For any $(x_1,x_2) \in X_2$, and any $\xi \in \liekernel_{\bt(x_1)}$ 
(resp. $\xi \in \liekernel_{\bt(x_2)}$), we denote by 
$\xi^{\btr}_1$ (resp. $\xi^{\btr}_2$) 
the tangent vectors in $T_{(x_1,x_2)} X_2$ corresponding 
infinitesimally to these two actions.
Denote by $p_1,m,p_2$ the three face maps from $X_2 $ to 
$X_1 $. For any $\xi \in \liekernel$, we have
\begin{align*} 
p_{1*} \xi^{\btr}_1 &= \xi^{\btr}, 
& m_* \xi^{\btr}_1 &= (\Ad_{x_2\inv} \xi)^{\btr},
& p_{2*} \xi^{\btr}_1 &= 0, \\
p_{1*} \xi^{\btr}_2 &= 0, 
& m_* \xi^{\btr}_2 &= \xi^{\btr}, 
& p_{2*} \xi^{\btr}_2 &= \xi^{\btr} .
\end{align*}
Thus it is routine to check that the two relations
$\xi^{\btr}_1 \ii \ptr \theta = 0$ and 
$\xi^{\btr}_2 \ii \ptr \theta = 0$
hold, which implies that $\ptr \theta$ is
horizontal.
\end{proof}

\begin{defn}
The \textbf{$\ptr $-cohomology of
a Lie groupoid extension} $X\simplicial\xrightarrow{\phi}Y\simplicial$ 
is the cohomology of the cochain complex 
$\big((\Omega^l(X_n,\bt^*\liekernel))_{n\in\N},\ptr \big)$.
It is denoted by $\pcoho{n}{l}{X}{Y}$.
\end{defn}

The following proposition can be verified directly.

\begin{prop}
Let $X_1\xrightarrow{\phi}Y_1\toto M$ be a Lie groupoid extension.
Then for any fixed $l\in\N$, $\Omega^l_{\hor}(X\simplicial,\bt^*\liekernel)$ 
is stable under $\ptr $, i.e.
$$\ptr  \big(\Omega^l_{\hor}(X_n,\bt^*\liekernel)\big)
\subset \Omega^l_{\hor}(X_{n+1},\bt^*\liekernel) .$$
\end{prop}

Therefore, for any fixed $l\in\N$, the horizontal $l$-forms 
$\big((\Omega^l_{\hor}(X_n,\bt^*\liekernel))_{n\in\N},\ptr \big)$ 
form a cochain subcomplex.

\begin{defn}
The \textbf{horizontal cohomology of a Lie groupoid extension} 
$X\simplicial\xrightarrow{\phi}Y\simplicial$ 
is the cohomology of the cochain complex 
$\big((\Omega^l_{\hor}(X_n,\bt^*\liekernel))_{n\in\N},\ptr \big)$.
It is denoted by $\hcoho{n}{l}{X}{Y}$.
\end{defn}

\subsection{Horizontal cohomology and strict homomorphisms}

Let $f$ be a homomorphism of Lie groupoid extensions
from $X\simplicial'\to Y\simplicial'$ to $X\simplicial\to Y\simplicial$. 
As usual, denote by $\liekernel\to M$ (resp. $\liekernel'\to M'$)
the Lie algebra bundle associated to the groupoid extension $X_1\to Y_1\toto M $ 
(resp. $X'_1\to Y'_1\toto M'$). In general, there is
no straightforward way to make sense of the pull-back through $f$ of
a $\liekernel$-valued form on $X_1$ (or, more generally on $X_n$
for some $n\in\N$) such that the pull-back form is $\liekernel'$-valued.
Indeed, if $\zeta$ is an $l$-form on $X_n$ with values in $\bt^*\liekernel$, 
the usual pull-back $f^*\zeta$ takes its values in $(f^*\smalcirc\bt^*)\liekernel$ 
rather than in ${\bt'}^{*}\liekernel'$.
To build up a form on $X_n'$ with values in ${\bt'}^*\liekernel'$, 
one needs to identify $(f^*\smalcirc\bt^*)\liekernel$ with ${\bt'}^{*}\liekernel'$ 
at any point $(x_1',\cdots,x_n')\in X_n'$.
When $f$ is a strict homomorphism of Lie groupoid extensions, in particular a Morita morphism, 
$f$ establishes an isomorphism between $\kernel_{\bt(x'_n)}$ and 
$\kernel_{\bt(f (x_n ))}$. 
Hence its differential at the identity $f_{\bt(x'_n)*}$ 
gives the required identification.
This leads to the following definition of the pull back 
$\Omega^l(X_n,\bt^*\liekernel)\xrightarrow{f^*}\Omega^l(X'_n,{\bt'}^*\liekernel')$:
\begin{equation} \label{eq:defpull}
(f^*\zeta)(e_1,\cdots,e_l) = f_{\bt(x'_n)*}\inv
\big(\zeta(f_* e_1,\cdots,f_* e_l)\big), \quad\forall e_1,\cdots,e_l\in T_{(x_1,\cdots,x_n)}X_n .  
\end{equation}

The following lemma can be verified directly.

\begin{lem} \label{lem:pb}
Let $f$ be a strict homomorphism of Lie groupoid extensions
from $X\simplicial'\to Y\simplicial'$ to $X\simplicial\to Y\simplicial$. 
Then
\begin{enumerate}
\item the pull-back map
$\Omega^l(X_n,\bt^*\liekernel)\xrightarrow{f^*}\Omega^l(X'_n,{\bt'}^*\liekernel')$ 
is a chain map with respect to $\ptr $;
\item the pull-back of an horizontal form is again horizontal;
\item $\theta\in\Omega^1(X_1,\bt^*\liekernel)$ satisfies 
$\xi^\btr\ii\theta=\xi$, for all $\xi\in\liekernel$, 
if, and only if, $f^*\theta\in\Omega^1(X'_1,\bt^*\liekernel')$ satisfies 
${\xi'}^\btr\ii f^*\theta=\xi'$, for all $\xi'\in\liekernel'$.
\end{enumerate}
\end{lem}

As an immediate consequence, we have the following

\begin{cor} \label{cor:5.8}
Let $f$ be a strict homomorphism of Lie groupoid extensions
from $X\simplicial'\to Y\simplicial'$ to $X\simplicial\to Y\simplicial$. 
Then the pull-back map $f^*$ gives rise to morphisms
\begin{gather} \label{eq:pbinco2}
f^*:\pcoho{n}{l}{X}{Y}\lto \pcoho{n}{l}{X'}{Y'}
\\ \label{eq:pbinco}
f^*:\hcoho{n}{l}{X}{Y}\lto \hcoho{n}{l}{X'}{Y'} ,
\end{gather}
for all $n,l\in\N$.
\end{cor}

\subsection{Morita invariance} \label{sec:coho2}

In this subsection, we show that both the horizontal cohomology 
and the $\ptr $-cohomology are invariant under Morita
equivalence. In other words, they turn out to be cohomological objects associated to a gerbe.
Our proof relies on the notion of generalized morphisms introduced in Section \ref{sec:gen}, 
and is similar to the proof of Proposition 2.5 in \cite{TXL}.

\begin{thm} \label{prop:MoriHor}
\begin{enumerate}
\item Morita equivalent groupoid extensions have isomorphic $\ptr $-cohomologies.
\item Morita equivalent groupoid extensions have isomorphic horizontal cohomologies. 
\item In particular, if $f$ is a Morita morphism of groupoid extensions, 
the maps $f^*$ defined in Eq.~\eqref{eq:pbinco2}-\eqref{eq:pbinco} are isomorphisms.
\end{enumerate}
\end{thm}

We need two lemmas.

\begin{lem}
\label{lem:6.10}
If $f$ is a Morita morphism of Lie groupoid extensions
from $X'_1\to Y'_1\toto M'$ to $X_1\to Y_1\toto M$
such that $f:M'\to M$ is an \'etale map, then
the pull-back map $f^*$ as in Eqs.~\eqref{eq:pbinco2}-\eqref{eq:pbinco}
induces isomorphisms on both horizontal and
$\ptr $-cohomologies.
\end{lem}

\begin{proof}
Note that $U\mapsto\Omega_{\hor}^l(U,\bt_{|U}^*\liekernel)$ 
(resp. $U\mapsto\Omega^l(U,\bt_{|U}^*\liekernel$) is a sheaf over 
the simplicial manifold $(X_n)_{n\in\N}$ (see \cite{tu:05}).
And the horizontal cohomology $\hcoho{n}{l}{X}{Y}$ 
(resp. the $\ptr $-cohomology
$\pcoho{n}{l}{X}{Y}$) is the corresponding sheaf cohomology.
When $f$ is an \'etale map, it is clear that the pull-back sheaf 
$f^*\big(\Omega_{\hor}^l(-,\bt^*\liekernel)\big)$ 
(resp. $f^*\big(\Omega^l(-,\bt^*\liekernel)\big)$) 
is isomorphic to 
$\Omega_{\hor}^l(-,{\bt'}^*\liekernel')$
(resp. $\Omega^l(-,{\bt'}^*\liekernel')$).
According to Theorem 8.1 in \cite{tu:05}, 
$f^*$ induces an isomorphism in cohomology:
\begin{gather} 
\label{eq:JLT} f^*:\hcoho{n}{l}{X}{Y}\xrightarrow{\isomorphism}\hcoho{n}{l}{X'}{Y'} \\
f^*:\pcoho{n}{l}{X}{Y}\xrightarrow{\isomorphism}\pcoho{n}{l}{X'}{Y'} . 
\end{gather}          
\end{proof}

\begin{prop} \label{pro:gen}
\begin{enumerate}
\item Any generalized homomorphism of Lie groupoid extensions $F$ 
from $X'\simplicial\xrightarrow{\phi'}Y'\simplicial$
to $X\simplicial\xrightarrow{\phi}Y\simplicial$
induces canonical homomorphisms
\begin{gather}
\label{eq:pbinc2} F^*:\pcoho{n}{l}{X}{Y}\lto\pcoho{n}{l}{X'}{Y'} \\
\label{eq:pbinc} F^*:\hcoho{n}{l}{X}{Y}\lto\hcoho{n}{l}{X'}{Y'} ,
\end{gather}
for all $n,l\in\N$. 
\item In particular, if $F$ is a generalized homomorphism of Lie groupoid extensions
induced by a strict morphism of groupoid extensions $f$, then $F^*$ coincides with
the pull-back map $f^*$ in Corollary \ref{cor:5.8}.
\item Moreover, the relation
$$(F_1\smalcirc F_2)^*=F_2^* \smalcirc F_1^*$$
holds for any composable pair $F_1,F_2$ of generalized homomorphisms 
of Lie groupoid extensions.
\end{enumerate}
\end{prop}

\begin{proof}
\begin{enumerate}
\item Assume that $F$ is given by the bimodule
$M'\xleftarrow{g}B\xrightarrow{f}M$.
According to Proposition~\ref{prop:decom}, $F$ is the composition 
of the canonical Morita equivalence between $X'\simplicial\to Y'\simplicial$ 
and $X'\simplicial[B]\to Y'\simplicial[B]$,
with a strict homomorphism of groupoid extensions
from $X'\simplicial[B]\to Y'\simplicial[B]$ to $X\simplicial\to Y\simplicial$.
Since $B\xrightarrow{g}M'$ is a surjective submersion, it admits
local sections. Hence, there exists an open cover $(U_i)_{i\in I}$ of $M'$ 
and local sections $\sigma_i$ of $B\xrightarrow{g}M'$ relative to this open cover.
The sections $(\sigma_i)$ induce a strict homomorphism $\psi_{\UU}$
from $X'\simplicial[\UU]\to Y'\simplicial[\UU]$ to $X\simplicial\to Y\simplicial$
by composing the natural strict homomorphism 
from $X'\simplicial[\UU]\to Y'\simplicial[\UU]$ to $X'\simplicial[B]\to Y'\simplicial[B]$ 
with the one from $X'\simplicial[B]\to Y'\simplicial[B]$ to $X\simplicial\to Y\simplicial$.
Denote by $\chi_\UU$ the Morita morphism 
from $X'\simplicial[\UU]\to Y'\simplicial[\UU]$ to $X'\simplicial\to Y'\simplicial$.

We now define $F^*$ as the composition
$$\hcoho{n}{l}{X}{Y}\xrightarrow{\psi^*_\UU}\hcohop{n}{l}{X'}{Y'}{\UU}\underset{\isomorphism}{\xrightarrow{(\chi^*_\UU)\inv}}\hcoho{n}{l}{X'}{Y'}$$
and similarly
$$\pcoho{n}{l}{X}{Y}\xrightarrow{\psi^*_\UU}\pcohop{n}{l}{X'}{Y'}{\UU}\underset{\isomorphism}{\xrightarrow{(\chi^*_\UU)\inv}}\pcoho{n}{l}{X'}{Y'} .$$
Here we have used the fact that $\chi^*_\UU$ is an isomorphism by Lemma~\ref{lem:6.10}.

We need to check that $F^*$ does not depend on
the choice of the open cover and the local sections.
Let $\VV:=(V_j,\tilde{\tau}_j)_{j \in J}$
be another choice of such local sections. 
The union of $\UU$ and $\VV$ is another such open cover of $M'$ 
together with a set of local sections of $B\to M'$.
It is simple to see that we have the following commutative diagram, 
where all the arrows are strict homomorphisms of groupoid extensions:
$$\xymatrix{
X'\simplicial \ar@2{-}[r] & X'\simplicial \ar@2{-}[r] & X'\simplicial \\ 
X'\simplicial[\UU] \ar[u]_{\chi_\UU} \ar[d]_{\psi_\UU} \ar[r] 
& X'\simplicial[\UU\cup\VV] \ar[u]_{\chi_{\UU\cup\VV}} \ar[d]_{\psi_{\UU\cup\VV}} 
& X'\simplicial[\VV] \ar[u]_{\chi_\VV} \ar[d]_{\psi_\VV} \ar[l] \\
X\simplicial \ar@2{-}[r] & X\simplicial \ar@2{-}[r] & X\simplicial
}$$
It follows from a diagram chasing argument that
$$(\chi_\UU^*)\inv \smalcirc\psi_\UU^*=
(\chi_{\UU\cup\VV}^*)\inv \smalcirc\psi_{\UU\cup\VV}^*=
(\chi_\VV^*)\inv \smalcirc\psi_\VV^* .$$
Hence $F^*$ is indeed well-defined.
\item If $F$ is induced from a strict homomorphism of Lie groupoid 
extensions $f$, 
in other words, if $F$ is given by the bimodule $M'\xleftarrow{\idn}M'\xrightarrow{f}M$, 
then $\chi=\idn$ and $\psi$ coincides with the map $f^*$ defined in Corollary~\ref{cor:5.8}.
\item Let $X\simplicial\to Y\simplicial$, $X'\simplicial\to Y'\simplicial$ and 
$X''\simplicial\to Y''\simplicial$ be Lie groupoid extensions with base manifolds $M$, $M'$ and $M''$, respectively.
Let $F_1$ and $F_2$ be generalized homomorphisms 
from $X'\simplicial\to Y'\simplicial$ to $X\simplicial\to Y\simplicial$ and 
from $X''\simplicial\to Y''\simplicial$ to $X'\simplicial\to Y'\simplicial$,
 respectively.
According to the first part of the proof, there exists an open covering
$\UU =(U_k)$ of $M'$ and a strict homomorphism $\psi_\UU$ of groupoid extensions 
from $X'\simplicial[\UU]\to Y'\simplicial[\UU]$ to $X\simplicial\to Y\simplicial$ 
such that $F_1^*=(\chi_\UU^*)\inv \smalcirc\psi_{\UU}$. 
The same holds for $F_2$. However, we can choose an open covering
$\VV=(V_l)_{l\in L}$ of $M''$ and a strict homomorphism of groupoid extensions 
$\psi_{\VV}$ such that $\psi_{\VV}(V)$ is contained in an open set $U\in\UU$, for all $V\in\VV$.
Composing $\psi_{\VV}:\coprod_l V_l\to M'$ with these inclusions yields a map 
$j:\coprod_l V_l\to\coprod_k U_k$.
We can then construct a strict homomorphism $\mu_j$ of groupoid extensions
from $X''\simplicial[\VV]\to Y''\simplicial[\VV]$ to $X'\simplicial[\UU]\to Y'\simplicial[\UU]$:
$$ \left\{ \begin{array}{l} 
\mu_j(v_1,x,v_2)=\big(j(v_1),\psi_{\VV}(v_1,x,v_2),j(v_2)\big) \\
\mu_j(v_1,y,v_2)=\big(j(v_1),\psi_{\VV}(v_1,y,v_2),j(v_2)\big) 
\end{array} \right. $$
for all $x\in X_1''$, $y \in Y_1''$, $v_1,v_2\in\coprod_l V_l$ 
with $\chi_{\VV}(v_1)=\bs(x)=\bs(y)$ and $\chi_{\VV}(v_2)=\bt(x)=\bt(y)$.
We then have the following commutative diagram of strict homomorphisms
of Lie groupoid extensions:
\begin{equation} \label{eq:commdiag} 
\xymatrix{ X''\simplicial & & X'\simplicial & & X\simplicial \\
& X''\simplicial[\VV] \ar[ul]^{\chi_\VV} \ar[ur]^{\psi_\VV} \ar[rr]_{\mu_j} & & 
X'\simplicial[\UU] \ar[ul]_{\chi_\UU} \ar[ur]_{\psi_\UU} & } 
\end{equation}

The composition $\psi_{\UU}\smalcirc\mu_j$ is a strict homomorphism of Lie groupoid extensions 
from $X''\simplicial[\VV]\to Y''\simplicial[\VV]$ to $X\simplicial\to Y\simplicial$. 
It is obvious from Diagram~\eqref{eq:commdiag} that 
$$X''\simplicial\xleftarrow{\chi_\VV}X''\simplicial[\VV]\xrightarrow{\psi_\UU\rond\mu_j}X\simplicial$$
is a bimodule representing the generalized homomorphism $F_1\rond F_2$. 
Therefore, 
\[ (F_1\rond F_2)^* = (\chi_{\VV}^*)\inv\rond(\psi_{\UU}\rond\mu_j)^* = (\chi_{\VV}^*)\inv\rond\mu_j^*\rond\psi_{\UU}^* 
= (\chi_{\VV}^*)\inv\rond(\psi_{\VV})^*\rond(\chi_{\UU}^*)\inv\rond\psi_{\UU}^* = F_2^*\rond F_1^* .\] 
\end{enumerate}
\end{proof}

\begin{proof}[Proof of Theorem \ref{prop:MoriHor}]
If $F$ is a Morita equivalence
from $X'\simplicial\to Y'\simplicial$ to $X\simplicial\to Y\simplicial$ 
given by the bitorsor $M'\xleftarrow{f}B\xrightarrow{g}M$, 
then, according to Proposition \ref{pro:com}, $F\inv\rond F$ is a Morita equivalence 
from $X\simplicial\to Y\simplicial$ to itself, 
whose corresponding bitorsor is given by $M\leftarrow(B\times_{M'}B)/X'_1\rightarrow M$. 
The later is canonically isomorphic to $M\leftarrow X_1\rightarrow M$, 
which is indeed the generalized morphism corresponding to the identity strict homomorphism. 
Thus it follows from Proposition~\ref{pro:gen} that $(F\inv\rond F)^*=\idn$. 
Therefore the conclusion follows from Proposition~\ref{pro:gen} immediately.
\end{proof}

\begin{rmk}
Note that the isomorphism between the  horizontal
cohomologies  of Morita equivalent groupoid extensions
is not canonical. It depends on the  choice of the Morita  equivalence bitorsors, as seen in the above proof.
\end{rmk}

\subsection{Obstruction class to the existence of  connections}
\label{sec:coho3}

The main purpose of this subsection is to introduce a
characteristic class $[\obs]$ in $\hcoho{2}{1}{X}{Y}$,
which measures the existence of connections 
on a Lie groupoid extension $X\simplicial\xrightarrow{\phi}Y\simplicial$.

Forgetting about the groupoid structure, consider any horizontal distribution 
$H$ on the fiber bundle $X_1\xrightarrow{\phi}Y_1$.
Let $\theta \in \Omega^1(X_1,\bt^*\liekernel)$ be its associated
$\liekernel$-valued 1-form, i.e. the unique 1-form that
vanishes on $H$ and satisfies Eq.~\eqref{eq:rightaction}. 
According to Proposition~\ref{ex:hor2form},
$\ptr \theta$ is a $\liekernel$-valued horizontal
1-form on $X_2$. Moreover, since $(\ptr)^2=0$, 
$\ptr\theta$ is a 2-cocycle of the horizontal cohomology, 
and hence it defines a class in $\hcoho{2}{1}{X}{Y}$ that we call the
\textbf{obstruction class} and denote simply by $[\obs]$. 
This terminology is justified by the following

\begin{prop} \label{pro:exist}
\begin{enumerate}
\item The class $[\obs]\in\hcoho{2}{1}{X}{Y}$ does not depend on the particular choice 
of an horizontal distribution $H$ on the fiber bundle $X_1\xrightarrow{\phi}Y_1$.
\item The class $[\obs]$ vanishes if, and only if, the groupoid extension 
admits a groupoid extension connection.
\end{enumerate}
\end{prop}

\begin{proof}
\begin{enumerate} 
\item Let $H$ and $H'$ be any two horizontal distributions on
the fiber bundle $X_1\xrightarrow{\phi}Y_1$, 
and $\theta,\theta'\in\Omega^1(X_1,\bt^*\liekernel)$
their associated $\liekernel$-valued 1-forms. 
By construction, $\theta-\theta'$ is an horizontal 1-form on $X_1$.
Therefore $\ptr \theta$ and
$\ptr \theta'$ differ by a horizontal coboundary.
Hence $[\ptr \theta]=[\ptr \theta']$. 
This proves (1).
\item Assume that $[\obs]=0$. 
The $\liekernel$-valued 1-form $\theta\in\Omega^1(X_1,\bt^*\liekernel)$ 
associated to any horizontal distribution $H$ on the fiber bundle 
$X_1\xrightarrow{\phi}Y_1$ gives a representative $\ptr\theta$ for the obstruction class $[\obs]$.
Hence, by assumption, there exists an horizontal $\liekernel$-valued 1-form $\zeta$ on $X_1$ such that $\ptr\theta=\ptr\zeta $. It is simple to check that $\alpha=\theta-\zeta\in\Omega^1(X_1,\bt^*\liekernel)$ 
is indeed a right connection 1-form for the groupoid extension 
$X\simplicial\xrightarrow{\phi}Y\simplicial$.
\end{enumerate}
\end{proof}

Combining Proposition~\ref{pro:exist} with the Morita invariance of
the horizontal cohomology obtained in Proposition~\ref{prop:MoriHor}, 
we are led to the following main result of this subsection.

\begin{thm} 
\label{th:ObsMor}\begin{enumerate}
\item Given a Morita morphism $f$ of Lie groupoid
extensions  from $X'\simplicial \to Y'\simplicial$
to $X\simplicial \to Y\simplicial$ as in Eq.~\eqref{eq:morita}, 
let $[\obs']$ and $[\obs]$ be their obstruction classes. 
Then $[\obs']=f^*[\obs]$, where 
$f^*:\hcoho{2}{1}{X}{Y}\to\hcoho{2}{1}{X'}{Y'}$ 
is the homomorphism  given by Eq.~\eqref{eq:pbinco}.
\item If a Lie groupoid extension admits a connection, 
so does any Morita equivalent Lie groupoid extension.
\end{enumerate}
\end{thm}

\begin{proof} 
Observe that (2) is a trivial consequence of (1).
Now let $\theta$ be a $\liekernel$-valued 1-form on $X_1$ satisfying Eq.~\eqref{eq:rightaction}. 
Then, according to Lemma~\ref{lem:pb}(3), the pull-back $f^*\theta\in\Omega^1(X'_1,{\bt'}^*\liekernel')$ 
satisfies Eq.~\eqref{eq:rightaction} as well.
By Lemma~\ref{lem:pb}(1), the relation $\ptr(f^*\theta)=f^*(\ptr\theta)$ holds, 
which implies that $[\obs']=f^*[\obs]$. 
\end{proof}

\subsection{Moduli space of connections on a Lie groupoid extension}

We now study the space of right connection 1-forms of a Lie groupoid extension
$X_1\xrightarrow{\phi}Y_1\toto M$ assuming that the obstruction class vanishes.
Denote by $Z^1_{\hor}(X_1,\bt^*\liekernel)$ (resp. $Z^1_{\ptr}(X_1,\bt^*\liekernel)$) 
the space of horizontal 1-forms in $\Omega^1_{\hor}(X_1,\bt^*\liekernel)$
(resp. 1-forms in $\Omega^1(X_1,\bt^*\liekernel)$) which are $\ptr$-closed.

\begin{prop} \label{prop:affine1} 
For any Lie groupoid extension $X_1\xrightarrow{\phi}Y_1\toto M$ with vanishing obstruction class,
the space of right connection $1$-forms is an affine subspace of $Z^1(X_1,\bt^*\liekernel)$
with underlying vector space $Z^1_{\hor}(X_1,\bt^*\liekernel)$.
\end{prop}

\begin{proof}
Since the obstruction class vanishes, there exists a right connection form 
$\alpha\in\Omega^1(X_1,\bt^*\liekernel)$.
It is easy to see that a 1-form $\alpha'\in\Omega^1(X_1,\bt^*\liekernel)$ 
is a right connection 1-form (i.e. is $\ptr$-closed and satisfies Eq.~\eqref{eq:rightaction}) 
if, and only if, $\alpha-\alpha'\in Z^1_{\hor}(X_1,\bt^*\liekernel)$.
\end{proof}

Two right connection 1-forms $\alpha$ and $\alpha'$ on a Lie groupoid extension
$X_1\xrightarrow{\phi}Y_1\toto M$ are said to be equivalent if, and only if,
$\alpha-\alpha'=\ptr\beta$ for some $\beta\in\Omega^1(M,\liekernel)$. 
By $\MM$ we denote the \textbf{moduli space of right connection
 1-forms on a Lie groupoid extension}\footnote{Readers should not confuse
this with moduli space of flat connections in gauge theory. Here there  are
no gauge groups involved.}
 $X_1\xrightarrow{\phi}Y_1\toto M$. The following proposition is thus immediate.

\begin{prop} \label{prop:affine2}
For a Lie groupoid extension $X_1\xrightarrow{\phi}Y_1\toto M$
with vanishing obstruction class, the moduli space of right connection 1-forms $\MM$ 
is an affine subspace of $\pcoho{1}{1}{X}{Y}$ with underlying vector space $\hcoho{1}{1}{X}{Y}$.
\end{prop}

We now describe the relation between the moduli spaces of right connection 
1-forms of Morita equivalent extensions. We start with a lemma.

\begin{lem}
\label{lem:pbalp}
Assume that $f$ is a Morita morphism of Lie groupoid extensions 
from $X'\simplicial\xrightarrow{\phi'}Y'\simplicial$ to $X\simplicial\xrightarrow{\phi}Y\simplicial$.
A 1-form $\alpha\in\Omega^1(X_1,\bt^*\liekernel)$ is a right connection 1-form of the Lie groupoid extension
$X\xrightarrow{\phi}Y\toto M$ if, and only if, $f^*\alpha\in\Omega^1(X'_1,{\bt'}^*\liekernel')$
is a right connection 1-form of the Lie groupoid extension $X'\simplicial\xrightarrow{\phi}Y'\simplicial$.
\end{lem}

\begin{proof}
For any $\alpha\in\Omega^1(X_1,\bt^*\liekernel)$, according to Lemma~\ref{lem:pb}(1), we have 
$f^*\ptr\alpha=\ptr f^*\alpha=0$. Hence $\alpha$ is $\ptr$-closed if, and only if, so is $f^*\alpha$. 
According to Lemma~\ref{lem:pb}(3), Eq.~\eqref{eq:rightaction} holds for $\alpha$ if, and only if, 
it holds for $f^*\alpha$. 
Hence it follows that $\alpha$ is a right connection 1-form on $X_1\xrightarrow{\phi}Y_1\toto M$ 
if, and only if, $f^*\alpha$ is a right connection 1-form on $X'\xrightarrow{\phi}Y'\toto M'$.
\end{proof}

\begin{thm} \label{th:11corr}
\begin{enumerate}
\item For any Morita morphism of Lie groupoid extensions $f$ from $X'_1\xrightarrow{\phi'}Y'_1\toto M'$
to $X_1\xrightarrow{\phi}Y_1\toto M$, the pull-back map
$f^*:\pcoho{1}{1}{X}{Y}\to \pcoho{1}{1}{X'}{Y'}$
induces an isomorphism of their corresponding moduli spaces 
of right connection 1-forms $\MM$ and $\MM'$.
\item The moduli spaces of right connection 1-forms of Morita equivalent extensions are isomorphic.
\end{enumerate}
\end{thm}

\begin{proof}
By Theorem~\ref{prop:MoriHor}, the map $f^*:\pcoho{1}{1}{X}{Y}\to\pcoho{1}{1}{X'}{Y'}$ 
is an isomorphism, which, according to Lemma~\ref{lem:pbalp}, sends $\MM$ into $\MM'$. 
It is thus enough to show that the restriction of $f^*$ to $\MM$ is surjective onto $\MM'$.
Assume that $[\alpha']\in\MM'$.
Since $f^*:\pcoho{1}{1}{X}{Y}\to \pcoho{1}{1}{X'}{Y'}$ is surjective 
by Theorem~\ref{prop:MoriHor}(1), there exists $\alpha\in Z_{\ptr}^1(X_1,\bt^*\liekernel)$
such that $f^*[\alpha]=[\alpha']$.
That is, $f^*\alpha=\alpha'+\ptr\beta'$ for some $\beta'\in\Omega^1(M',\liekernel')$.
However, since $\alpha'+\ptr\beta'$ is also a right connection 1-form
on the groupoid extension $X_1'\xrightarrow{\phi}Y_1'\toto M'$,
according to Lemma~\ref{lem:pbalp}, $\alpha$ must be a right connection 1-form
on the groupoid extension $X_1\xrightarrow{\phi}Y_1\toto M$.
This concludes the proof.
\end{proof}

\subsection{Connections for $G$-gerbes over manifolds} \label{sec:coho4}

The main goal of this subsection is to compute the horizontal
cohomology of a groupoid $G$-extension of a \v{C}ech groupoid, 
i.e. of a $G$-gerbe over a manifold, where $G$ is a connected Lie group. 
As a consequence, we show the existence of connections.

First we need to introduce some preliminary constructions.
%More precisely, we introduce a vector bundle whose sections
%represent classes of the horizontal cohomology. 
We do not have to restrict ourselves to the case of $G$-gerbes 
over manifolds at this point, although this is the only case that
we are concerned about in applications.

Let $X_1\xrightarrow{\phi}Y_1\toto M$ be a Lie groupoid $G$-extension.
Its kernel is the Lie group bundle $\kernel\to M$
with Lie algebra bundle $\liekernel\to M$. For any $n\in\N$, consider 
the fiberwise group cohomology $\coprod_m H^n(\kernel_m,\liekernel_m)$
of $\kernel$ with values in $\liekernel$ relative to the adjoint action,
which is clearly a vector bundle over $M$, denoted by $H^n(\kernel,\liekernel)\to M$.

As we have seen in Section~\ref{sec:2.1}, any element $x\in X_1$ defines a
Lie group isomorphism $\AD_x:\kernel_{\bt(x)}\to \kernel_{\bs(x)}$, 
and thus a Lie algebra isomorphism $\Ad_x:\liekernel_{\bt(x)}\to \liekernel_{\bs(x)}$. 
These isomorphisms induce an isomorphism of cochain complexes mapping 
$f\in\cty(\kernel_{\bs(x)}\times\cdots\times\kernel_{\bs(x)},\liekernel_{\bs(x)})$ 
to $(k_1,\cdots,k_n)\mapsto \Ad_{x\inv}f(\Ad_x k_1,\cdots,\Ad_x k_n)$ 
in
 $\cty(\kernel_{\bt(x)}\times\cdots\times\kernel_{\bt(x)},\liekernel_{\bt(x)})$.
Therefore it induces an isomorphism 
$H^n(\kernel_{\bs(x)},\liekernel_{\bs(x)})\xrightarrow{\sim}H^n(\kernel_{\bt(x)},\liekernel_{\bt(x)})$,
which is still denoted by $\AD_x$ by abuse of notation. In other
words, $H^n(\kernel,\liekernel)\to M$ is a vector bundle
over $M$ on which the groupoid $X_1\toto M$ acts from the left.
That is, $H^n(\kernel,\liekernel)\to M$ is a left vector
bundle over the groupoid $X_1\toto M$.
In fact, this action descends and induces an action of the groupoid
$Y_1\toto M$ since any element in the kernel $\kernel$ of 
$X\simplicial\xrightarrow{\phi}Y\simplicial$ acts trivially on
$H^n(\kernel,\liekernel)\to M$. 
Indeed for a connected Lie group $G$, it is well known that the inner 
automorphisms of $G$ preserve the cohomology classes of $H^n(G,\mfg)$. 
Thus we have proved the following

\begin{prop} \label{pro:Hn}
For all $n\in\N$, $H^n(\kernel,\liekernel)\to M$ is a vector bundle 
over the groupoid $Y_1\toto M$.
\end{prop}

%The sections of this vector bundle are related to the horizontal cohomology cocycles. 
Considering $\kernel$ as a Lie subgroupoid of $X_1\toto M$, one obtains a
monomorphism of groupoid extensions $i$ from $\kernel\to M\toto M$
to $X_1\to Y_1\toto M$. 
The horizontal cohomology of $\kernel\to M \toto M$ 
can be easily described by the following 

\begin{lem} \label{lem:horcomBG}
We have 
\begin{equation} \label{eq:isomomo} 
\hcoho{n}{l}{\kernel}{M} \isomorphism \Omega^{l}\big(M,H^n(\kernel,\liekernel)\big) ,
\end{equation} 
where $H^n(\kernel,\liekernel)\to M$ is the vector bundle constructed above 
(ignoring the $Y\simplicial$-action).
\end{lem}

\begin{proof}
The natural identification $\XX_{\hor}^l(\kernel_n)\isomorphism\sect{\pi^*\wedge^l TM}$ 
given by the differential of $\pi:\kernel_n\to M$ induces an isomorphism 
$\Omega_{\hor}^l(\kernel_n,\bt^*\liekernel)\isomorphism\Omega^l(M,E_n)$, 
where $E_n$ is the vector bundle $\coprod_{m\in M}\cty\big((\kernel|_m)^n,\liekernel_m\big)\to M$.
It is simple to check that, under this isomorphism,
the differential 
$\ptr:\Omega^l_{\hor}(\kernel_n,\bt^*\liekernel)\to\Omega^l_{\hor}(\kernel_{n+1},\bt^*\liekernel)$ 
becomes the operator $\partial:\Omega^l(M,E_n)\to\Omega^l(M,E_{n+1})$,
which is the fiberwise group cohomology differential 
$\cty\big((\kernel|_m)^n ,\liekernel_m\big)\to\cty\big((\kernel|_m)^{n+1},\liekernel_m\big)$. 
Taking its cohomology, one obtains an isomorphism 
$\hcoho{n}{l}{\kernel}{M}\isomorphism\Omega^{l}\big(M,H^n(\kernel,\liekernel)\big)$.
\end{proof}

For all $n\in\N$, from the induced map $i:\kernel_n\to X_n$, we obtain a chain map 
$i^*:\Omega_{\hor}^l(X_n,\bt^*\liekernel)\lto\Omega_{\hor}^l(\kernel_n,\bt^*\liekernel)$ 
and hence a map in cohomology 
\begin{equation} \label{eq:inclusion}
i^*:\hcoho{n}{l}{X}{Y}\lto\hcoho{n}{l}{\kernel}{M} .
\end{equation}

Composing the maps given by Eq.~\eqref{eq:inclusion} with Eq.~\eqref{eq:isomomo}, 
we are led to

\begin{prop} \label{pro:i*}
Let $X_1\to Y_1\toto M$ be a Lie groupoid $G$-extension, with kernel 
being the Lie group bundle $\kernel\to M$ with corresponding Lie algebra bundle 
$\liekernel\to M$. For any $l\in\N$, the inclusion map 
$i:\kernel\simplicial\to X\simplicial$ induces a natural restriction map 
\begin{equation} \label{eq:46}
i^*:\hcoho{n}{l}{X}{Y}\lto\Omega^{l}\big(M,H^n(\kernel,\liekernel)\big) ,
\end{equation}
where $H^n(\kernel,\liekernel)\to M$ is the vector bundle constructed above 
(ignoring the $Y\simplicial$-action).
\end{prop} 

It is simple to see that this restriction map $i^*$ is stable with respect to Morita morphisms of Lie groupoid extensions.

\begin{lem} \label{lem:if}
Let $f$ be a Morita morphism of Lie groupoid extensions 
from $X'\simplicial\xrightarrow{\phi'}Y'\simplicial$ to $X\simplicial\xrightarrow{\phi}Y\simplicial$ 
as in Diagram~\eqref{eq:morita}.
Let $\kernel$ and $\kernel'$ denote the kernels of $\phi$ and $\phi'$ respectively.
Then the following diagram commutes:
\[ \xymatrix{ \hcoho{n}{l}{X}{Y} \ar[d]_{i^*} \ar[r]^{{f}^*} &
\hcoho{n}{l}{X'}{Y'} \ar[d]^{i'^*} \\ 
\Omega^l\big(M,H^n(\kernel,\liekernel)\big) \ar[r]_{f^*} &
\Omega^l\big(M',H^n(\kernel',\liekernel')\big) . } \] 
\end{lem}

Let us turn our attention to the case
of a $G$-extension of a \v{C}ech groupoid
$\coprod_{ij}U_{ij}\toto \coprod_{i}U_i$ 
associated to an open covering $(U_j)_{j\in J}$ of a manifold $N$. 
In this case, vector bundles over the \v{C}ech groupoid 
$\coprod_{ij}U_{ij}\toto \coprod_{i}U_i$ 
are in 1-1 correspondence with ordinary vector bundles over the
manifold $N$. For all $n\in\N$, we denote by
$H^n(\kernel,\liekernel)_N$ the vector bundle over $N$
corresponding to the vector bundle $H^n(\kernel,\liekernel)$
over $\coprod_{ij}U_{ij}\toto \coprod_{i}U_i$ obtained 
as in Proposition~\ref{pro:Hn}.

\begin{thm} \label{th:manifcase}
Let $X_1\xrightarrow{\phi}Y_1\toto M$ denote a Lie groupoid $G$-extension 
$X_1\to\coprod_{ij}U_{ij}\toto\coprod_{j}U_{j}$ of a \v{C}ech groupoid 
$\coprod_{ij}U_{ij}\toto\coprod_{j}U_{j}$ associated to an open covering 
$\gendex{U_j}{j\in J}$ of a manifold $N$. Then, for any pair $n,l\in\N$, 
the map in Eq.~\eqref{eq:46} factorizes through
$\Omega^{l}\big(N,H^n(\kernel,\liekernel)_N\big)\to\Omega^{l}\big(M,H^n(\kernel,\liekernel)\big)$:
\begin{equation} \label{eq:uneetoile} 
\xymatrix{ & \hcoho{n}{l}{X}{Y} \ar[d]^{i^*} \ar[dl]_{\tilde{\imath}^*} \\
\Omega^{l}\big(N,H^n(\kernel,\liekernel)_N\big) \ar[r]_{p^*} & \Omega^{l}\big(M,H^n(\kernel,\liekernel)\big) } 
\end{equation}
Moreover, $\tilde{\imath}^*:\hcoho{n}{l}{X}{Y}\lto\Omega^{l}\big(N,H^n(\kernel,\liekernel)_N\big)$ is an isomorphism 
and the restriction map $i^*:\hcoho{n}{l}{X}{Y}\to\Omega^{l}\big(M,H^n(\kernel,\liekernel)\big)$ is injective.
\end{thm}

We need some preliminaries.
Any smooth function $f$ on the manifold $N$ can be pulled back to
a function $\tilde{f}_n$ on $X_n$ using any appropriate 
composition of consecutive face maps together with the projection
$\coprod_i U_i\xrightarrow{p}N$.
In other words, $\tilde{f}_n=p_n^*f$, where $p_n$ is the composition 
$X_n\xrightarrow{\phi}\coprod_{j_1,\cdots,j_n}U_{j_1\cdots j_n}\to N$.
The following result is straightforward.

\begin{lem} \label{lem:locality}
For any $\zeta\in\Omega^l(X_n,\bt^*\liekernel)$ 
and any function $f\in\cty(N)$, we have
\[ \ptr(\tilde{f}_n\zeta) = \tilde{f}_{n+1}\ptr\zeta .\]
\end{lem}

For any open set $U\subset N$, we denote by $X_1[p\inv U]\toto p\inv(U)$ 
and $Y_1[p\inv U]\toto p\inv(U)$ the restriction of the groupoids $X_1\toto M$
and $Y_1\toto M$ to the open submanifold $p\inv(U)$ of the unit space $M$,
respectively. It is clear that $X_1[p\inv U]\to Y_1[p\inv U]\toto p\inv(U)$ is a $G$-extension.

For any fixed integers $n$ and $l$, we define a pre-sheaf $\EEE{n}{l}$ over $N$ 
by $U\mapsto\hcohop{n}{l}{X}{Y}{p\inv U}$,  
the restriction maps $\rr{V}{U}:\hcohop{n}{l}{X}{Y}{p\inv U}\to\hcohop{n}{l}{X}{Y}{p\inv V}$ 
being the pull-backs of the natural inclusion of 
$X\simplicial[p\inv V]\to Y\simplicial[p\inv V]$ into $X\simplicial[p\inv U]\to Y\simplicial[p\inv U]$ 
for any open subsets $V\subset U$ of $N$.

\begin{lem} \label{pro:fine}
For any $n,l\in\N$, the presheaf $\EEE{n}{l}$ is a sheaf.
\end{lem}

\begin{proof}
First, let $(U_i)_{i\in I}$ be open subsets of $N$, $U=\cup_{i\in I}U_i$ and $(f_i)_{i\in I}$
a partition of unity of $U$ with $\supp(f_i)\subset U_i$ for all $i\in I$
(which exists since we consider paracompact manifolds only).
First, let $[\omega_i]\in\hcohop{n}{l}{X}{Y}{U_i}$
be the cohomology class of some horizontal cocycle $\omega_i\in\Omega^l(X_n[U_i],\bt^*\liekernel)$,
for all $i\in I$.
If $\rr{U_i\cap U_j}{U_i}[\omega_i]=\rr{U_i\cap U_j}{U_j}[\omega_j]$ for all $i,j\in I$, then,
thanks to Lemma~\ref{lem:locality}
$\omega=\sum_{i\in I}\widetilde{(f_i)}_n \ \omega_i$ is a 
horizontal cocycle in $\Omega^l(X_n[U],\bt^*\liekernel)$.
By construction, its class $[\omega]\in\hcohop{n}{l}{X}{Y}{U}$
satisfies $\rr{U_i}{U}[\omega]=[\omega_i]$, for all $i\in I$.

Secondly, let $[\omega]\in\hcohop{n}{l}{X}{Y}{U}$ be the class of 
some horizontal cocycle $\omega\in\Omega^l(X_n,\bt^*\liekernel)$. 
If $\rr{U_i}{U}[\omega]=0$ for all $i\in I$, then there exists 
an horizontal form $\eta_i\in\Omega^l(X_{n-1}[U_i],\bt^*\liekernel)$
such that $\ptr\eta_i =\omega_{|U_i}$ for all $i\in I$.
Set $\eta=\sum_{i\in I}\widetilde{(f_i)}_{n-1} \ \eta_i$.
By Lemma~\ref{lem:locality}, we have
$$\ptr\eta = \sum_{i\in I}\widetilde{(f_i)}_{n} \ \ptr\eta_i = \sum_{i\in I}\widetilde{(f_i)}_{n} \ \omega_{|U_i} = \omega .$$
As a consequence, we have $[\omega]=0$.
Therefore, $\EEE{n}{l}$ is a sheaf.
\end{proof}

We define two other sheaves over $N$: 
\begin{gather*}
\FFF{n}{l}:U\mapsto\Omega^l\big(U,H^n(\kernel,\liekernel)_N\big) , \\
\GGG{n}{l}:U\mapsto\Omega^l\big(p\inv U,H^n(\kernel,\liekernel)\big) . 
\end{gather*}
Moreover, the covering $M\xrightarrow{p}N$ induces a morphism of sheaves 
$\FFF{n}{l}\xrightarrow{p^*}\GGG{n}{l}$ and, by Proposition~\ref{pro:i*}, the 
inclusion $\kernel\simplicial\xrightarrow{i}X\simplicial$ induces a morphism 
of sheaves $\EEE{n}{l}\xrightarrow{i^*}\GGG{n}{l}$.

Now we need a general fact about sheaves over a manifold.

\begin{lem} \label{lem:sheaf}
Let $\EE$, $\FF$ and $\GG$ be sheaves over a manifold $N$,
and let $h:\EE\to\GG$ and $p:\FF\to\GG$ be morphisms of sheaves with $p$ injective.
Assume that for any contractible open set $U$, the
map $h_U:\EE(U)\to\GG(U)$ factorizes through 
\begin{equation} \label{eq:sheafcondition} 
\xymatrix{ & \EE(U) \ar[d]^{h_U} \ar[dl]_{\tilde{h}_U} \\ 
\FF(U) \ar[r]_{p_U} & \GG(U) } 
\end{equation}
where the map 
$\tilde{h}_U:\EE(U)\to\FF(U)$ is an isomorphism.
Then the morphism of sheaves 
$h:\EE\to\GG$ factorizes through 
$$\xymatrix{ & \EE \ar[d]^{h} \ar[dl]_{\tilde{h}} \\ 
\FF \ar[r]_{p} & \GG }$$
where $\tilde{h}:\EE\to\FF$ is an isomorphism of sheaves.
\end{lem}

We can now turn to the  proof of the our main results in this section.

\begin{proof}[Proof of Theorem~\ref{th:manifcase}]
It remains to show that the morphism of sheaves $i^*$ factorizes 
through the monomorphism $p^*$ as an isomorphism of sheaves $\tilde{\imath}^*$: 
\begin{equation} \label{eq:deuxetoiles}
\xymatrix{ & \EEE{n}{l} \ar[d]^{i^*} \ar[dl]_{\tilde{\imath}^*} \\ 
\FFF{n}{l} \ar[r]_{p^*} & \GGG{n}{l} }
\end{equation}
Indeed, evaluation of the presheaves $\EEE{n}{l}$, $\FFF{n}{l}$ and $\GGG{n}{l}$ in 
Diagram~\eqref{eq:deuxetoiles} on $N$, seen as an open subset of itself, yields 
Diagram~\eqref{eq:uneetoile}.

As a consequence of Lemma~\ref{lem:sheaf}, since $\EEE{n}{l}$, $\FFF{n}{l}$ and $\GGG{n}{l}$ 
are sheaves, it suffices to prove that, for any contractible open subset $U$ of $N$, 
there exists a group isomorphism $\tilde{\imath}^*_U$ such that 
\begin{equation} \label{eq:cdg} 
\xymatrix{ & \hcohop{n}{l}{X}{Y}{p\inv U} \ar[d]^{i^*_U} \ar[dl]_{\tilde{\imath}^*_U} \\
\Omega^{l}\big(U,H^n(\kernel,\liekernel)_N\big) \ar[r]_{p^*_U} & \Omega^{l}\big(p\inv U,H^n(\kernel,\liekernel)\big) } 
\end{equation}
commutes.

However, according to Corollary~\ref{lem:refinement}, the $G$-extension 
$X\simplicial[p\inv U]\to Y\simplicial[p\inv U]$ of the \v{C}ech groupoid 
$\coprod_{i,j} U_{ij}\cap U\toto\coprod_i U_{i}\cap U$ of the contractible 
open manifold $U$ has a refinement which is isomorphic to a trivial $G$-extension. 
In other words, there exists a refinement $\VV=\coprod_{\alpha\in A}V_{\alpha}$ of 
the open covering $\coprod_{i\in I}U_i\cap U$ of $p\inv U$ together with a 
Morita equivalence 
\[ \xymatrix{ 
U\times G \ar[d] & \coprod V_{\alpha\beta}\times G \ar[l] \ar[d] \ar[r] & X_1[p\inv U] \ar[d] \\ 
U \dar[d] & \coprod V_{\alpha\beta} \ar[l] \dar[d] \ar[r] & Y_1[p\inv U] \dar[d] \\ 
U & \coprod V_{\alpha} \ar[l] \ar[r] & p\inv(U) .
} \]

Therefore, applying Lemma~\eqref{lem:if} to each Morita morphism, one obtains the commutative diagram 
\[ \xymatrix{ 
\hcoho{n}{l}{(U\times G)}{U} \ar[r]^{(2)} \ar[d]_{(3)} & \hcohop{n}{l}{(U\times G)}{U}{\VV} \ar[d] 
& \hcohop{n}{l}{X}{Y}{p\inv U} \ar[l]_{(1)} \ar[d]^{i^*_U} \\
\Omega^l\big(U,H^n(G,\mfg)\big) \ar[r]_{(4)} & \Omega^l\big(\VV,H^n(G,\mfg)\big) & \Omega^l\big(p\inv U,H^n(G,\mfg)\big) \ar[l]^{(5)} 
} \]
where, the restrictions of $H^n(\kernel,\liekernel)$ to contractible open sets have been identified with the trivial vector bundles 
with fiber $H^n(G,\mfg)$.
Here $(3)$ is an isomorphism by Lemma~\ref{prop:MoriHor}, $(1)$ and $(2)$ are isomorphism by Theorem~\ref{prop:MoriHor} 
and the injection $(4)$ factorizes through the injection $(5)$ to yield $p^*_U$. Thus one obtains the commutative Diagram~\eqref{eq:cdg}.
\end{proof}

We can now state the main corollary of this result.

\begin{cor}\label{cor:exi}
 Any groupoid $G$-extension of a \v{C}ech groupoid admits a
connection.
\end{cor}
\begin{proof}
It suffices to prove that the obstruction 
class $[\obs] $ of any groupoid
$G$-extension of a \v{C}ech groupoid vanishes. Since the vanishing
of $[\obs]$  is preserved under Morita equivalence, we can assume
that $(U_j)_{j\in J}$ is a good open covering and 
$X_1\isomorphism \coprod_{ij}U_{ij}\times G$,
$Y_1\isomorphism  \coprod_{ij}U_{ij}$.
In particular  $\kernel\isomorphism \coprod_{i}U_i
\times G$ is the trivial group bundle.
 Consider the $ \mfg$-valued 1-form
on $X_1$ given by
 $\theta = p_2^* \theta_G $,
  where $p_2: \coprod_{ij} U_{ij} \times G \to G $
is the projection onto the second component and $\theta_G $ is the
right Maurer-Cartan form on the Lie group $G$. Since
$\theta$ satisfies Eq.~\eqref{eq:rightaction}, we have
$[\obs]=[\ptr\theta] $. Moreover,
by construction, one has
$i^*[\obs]=i^*[\ptr\theta] =[\partial\theta_G ]=0$, 
where $i^*:\hcoho{2}{1}{X}{Y}\to\hcoho{2}{1}{\kernel}{M}$
 is the restriction map.
By Theorem~\ref{th:manifcase}, $i^*$ is injective.
Hence the obstruction class $[\obs]$
vanishes.
\end{proof}

Note that for bundle gerbes, the existence of connections was proved by Murray \cite{Murray}.
Another consequence of Theorem \ref{th:manifcase}
is the following

\begin{cor}
\label{cor:flatband} 
Let $X\simplicial\xrightarrow{\phi}Y\simplicial$ be a Lie
 groupoid $G$-extension $X_1 \to
\coprod_{ij} U_{ij}  \toto \coprod_{j} U_{j} $
 of a \v{C}ech groupoid $\coprod_{ij} U_{ij}  \toto \coprod_{j} U_{j}$
associated  to   an open covering
$(U_j)_{j\in J}$ of a manifold $N$. Assume that
$\alpha$ is a right connection 1-form with Ehresmann curvature
$\omega \in \Omega^2( X_1, \bt^* \liekernel)$. 
\begin{enumerate}
\item \label{mn1} The class of $\omega$ in $\hcoho{1}{2}{X}{Y}$ is zero if, and only if,
 the band is flat.
\item \label{mn2} If $N$ is simply connected  and
 $[\omega]$ vanishes in $\hcoho{1}{2}{X}{Y}$, then the band is
trivial. And therefore the extension is central.
\end{enumerate}
\end{cor}

\begin{proof}
\begin{itemize}
\item[\ref{mn1}.] 
According to Corollary~\ref{cor:4.46}, the band is flat if, and only if, 
$[\omega^{\kernel}]=0$ in $\hcoho{1}{2}{\kernel}{M}$ or, equivalently, if, and only if, 
$[\omega]=0$ in $\hcoho{1}{2}{X}{Y}$ since $i^*:\hcoho{1}{2}{\kernel}{M}\to\hcoho{1}{2}{X}{Y}$ 
is injective.
\item[\ref{mn2}.] The band can be identified with an $\Out (G)$-principal
bundle over $N$. Hence, since $N$ is simply connected, it must be
trivial if it is flat.
\end{itemize}
\end{proof}

\subsection{Flat gerbes}
\label{sec:6.5}
First of all, we introduce the following 
definition which generalizes the same notion
in the case of abelian gerbes \cite{BX1,brylinski1, Hitchin, Murray}.

\begin{defn}
Given a Lie groupoid extension $X_1\xrightarrow{\phi}Y_1 \toto M$,
and a  right connection 1-form $\alpha \in \Omega^1 (X_1, \bt^* \liekernel  )$ with
Ehresmann curvature $\omega \in \Omega^2 (X_1, \bt^* \liekernel  )$,
\begin{enumerate}
\item a \textbf{curving} is a two-form $B \in \Omega^2 (M,{\liekernel})$ such that $\ptr B =\omega$;
\item and  given $(\alpha,B)$, $\Omega = d^{\nabla} B \in
\Omega^3 (M,{\liekernel})$ is called the \textbf{3-curvature} of
$(\alpha, B)$, where $d^{\nabla}: \Omega\degree (M,{\liekernel})\to
\Omega^{\scriptscriptstyle\bullet+1} (M,{\liekernel})$ is the exterior 
covariant derivative with respect to the induced
connection $\nabla$ on the Lie algebra bundle $\liekernel\to M$ as
in Section \ref{sec:4.5}.
\end{enumerate}
\end{defn}

\begin{prop}
\begin{enumerate}
\item Given a Lie groupoid extension $X_1\xrightarrow{\phi}Y_1\toto M$, and a right connection 1-form $\alpha
\in \Omega^1 (X_1, \bt^* \liekernel)$ with Ehresmann curvature $\omega \in \Omega^2 ( X_1,\bt^* \liekernel)$,
 a curving exists if and only if $[\omega]\in\hcoho{1}{2}{X}{Y}$ vanishes.
\item If $X_1\xrightarrow{\phi}Y_1\toto M$
is the groupoid extension of a \v{C}ech groupoid,
then a curving exists if and only if the band is flat.
\end{enumerate}
\end{prop}
\begin{proof}
The first assertion is straightforward. And the second follows from the first and Corollary \ref{cor:flatband}.
\end{proof}

\begin{rmk}
The above proposition indicates that the existence of both
connections and curvings on a $G$-gerbe over a manifold
would force it to be close to being a $G$-bound gerbe (or
an abelian gerbe).
\end{rmk}

The following lemma will be useful.
We denote by $\Ad_{x\inv }$ the isomorphism
$ \Gamma( \bs^* \liekernel \to M)\xrightarrow{\isomorphism}\Gamma( \bt^* \liekernel \to M)  $
obtained by mapping a section $\sigma \in  \Gamma( \bs^* \liekernel \to M)$
to the section $x \to \Ad_{x\inv } \sigma(x)$
in  $\Gamma( \bt^* \liekernel \to M)$.

\begin{lem}\label{lem:com}
For any $\eta \in \Omega^k(M,{\liekernel})$,
we have
$$  (d^{\nabla^{\bt}}  \ptr  - \ptr   d^{\nabla}) \eta  
=  [\alpha, \Ad_{x\inv } {\bs}^* \eta] .$$
\end{lem}
\begin{proof}
Since $\nabla^{\bt}$ is the pull-back of $\nabla$ via ${\bt}$, one has
$ (d^{\nabla^{\bt}}   {\bt^*} - {\bt^*} d^{\nabla}  ) \eta =0   $.
Therefore, it suffices to prove the following 
relation
 \begin{equation}\label{eq:sufficient}  
(d^{\nabla^{\bt}}  \Ad_{x\inv } {\bs}^* - \Ad_{x\inv } {\bs}^*  d^{\nabla}) \eta  
=  [\alpha, \Ad_{x\inv } {\bs}^* \eta] .
    \end{equation}
The latter is equivalent to:
  \begin{equation} \label{eq:sufficient1}  
\nabla^{\bt}_u \Ad_{x\inv } {\bs}^* \sigma - \Ad_{x\inv } {\bs}^*  \nabla_{\bs_* u} \sigma =
  [\alpha(u), \Ad_{x\inv } {\bs}^* \sigma], 
   \ \ \ \  \forall u \in T_x X_1 \end{equation}
where  $\sigma $ is any local section of ${\liekernel} \to M$ 
in a neighborhood of ${\bs}(x) $.

If $ u=\xi^{\btr}  $ is vertical, then 
$$  \nabla^{\bt}_u \Ad_{x\inv } {\bs}^* \sigma =  \left.\tfrac{d}{d\tau}\right|_{\tau =0}
    \Ad_{(x\cdot exp(\tau \xi))\inv  } {\bs}^* \sigma  = [\xi, \Ad_{x\inv } {\bs}^* \sigma  ] 
  = [ \alpha (u) , \Ad_{x\inv } {\bs}^* \sigma ]  .$$
Hence Eq.~\eqref{eq:sufficient1} holds.

On the other hand, if $ u \in H_x $ is horizontal, and $X \in\XX(X_1)$
is a horizontal vector field through $u$, then
according to Eq.~\eqref{eq:parthor}, we have 
  $  \nabla^{\bt}_u (\Ad_{x\inv  } \bs^* \sigma) =   \alpha([X, (\Ad_{x\inv  } \bs^* \sigma)^{\btr}])|_{x}$.
Since $(\Ad_{x\inv  } \bs^* \sigma)^{\btr}=  \bs^*\sigma^{\btl}$
and $\alpha = \Ad_{x\inv }\smalcirc 
 \beta$, we have 
   $  \nabla^{\bt}_u (\Ad_{x\inv  } \bs^* \sigma) =  
  \Ad_{x\inv } \smalcirc \beta ( [X,\bs^* \sigma^{\btl}])|_{x}$.
According to Eq.~\eqref{eq:parthor'}, we have
  $  \nabla^{\bt}_u (\Ad_{x\inv  } \bs^* \sigma) =
 \Ad_{x\inv }  \big(\nabla^{\bs}_u (\bs^* \sigma)\big)  =   \Ad_{x\inv } \bs^* \nabla_{\bs_* u} \sigma $.
Eq.~\eqref{eq:sufficient1} thus holds in this case. 
This concludes the proof.
 \end{proof}

  Let $ \liecenter \to M$ be the subbundle
of $\liekernel \to M $ defined in Remark \ref{rmk:uptocenter}.
We have the following:

\begin{thm}\label{theo:partial3curv}
Assume that $X_1\xrightarrow{\phi}Y_1\toto M$
is a Lie groupoid extension with connection and curving.
Let $\Omega \in \Omega^3(M,\liekernel) $
be the corresponding $3$-curvature. Then 
$\Omega \in  \Omega^3( M, Z(\liekernel) )$.
Moreover, we have
\begin{equation}
\label{eq:consis3curv}
d^{\nabla}\Omega= 0,   \ \ \ \ \ \ptr   \Omega  = 0. 
\end{equation}
\end{thm}
\begin{proof}
Denote the curving by $B \in \Omega^2(M,\liekernel)$.
When being restricted to the group bundle $\kernel\xrightarrow{\pi}M$,
the relation $\ptr   B =\omega$
yields that
$ \pi^* B - \Ad_{k\inv } \pi^* B = \omega^{\kernel}(k), \ \ \ \ \forall k \in \kernel$. 
Differentiating this relation with respect to $k$ at the identity, and using Eq.~\eqref{eq:curvLieK}, 
one obtains that 
 \begin{equation}\label{eq:labcurv}  -\ad_B = \omega^{\liekernel} .
\end{equation}
Recall that $\omega^{\liekernel} \in \Omega^2(M,\End(\liekernel))$
 is the curvature of the covariant derivative
on the Lie algebra bundle ${\liekernel} \to M$.
By the Bianchi identity,
the relation $d^{\nabla}\omega^{\liekernel}=0  $ holds.
We therefore have
   $$ 0= d^{\nabla} \omega^{\liekernel} = -d^{\nabla} \ad_B= -\ad_{d^{\nabla} B} = -\ad _{\Omega} .$$
As a consequence, the $3$-curvature takes its values in 
the subbundle $Z({\liekernel}) \to M$.
 
%According to Lemma \ref{lem:basicdiffgeo}(2),
 We have
  $d^{\nabla} \Omega = (d^{\nabla})^2 B= -\omega^{\liekernel}(B) = 
-\ad_B (B)= -[B,B]=0$
since $B$ is a $2$-form.
According to Eq.~\eqref{eq:labcurv}, the Bianchi identity (\ref{eq:Bianchi1})
reads 
$ d^{\nabla^{\bt }}  \ptr   B  
+ [\alpha, \ptr   B ] = -\ad_{{\bt}^*B} (\alpha) = [\alpha, {\bt^*} B]$.
The latter is equivalent to
  $ d^{\nabla^{\bt }}  \ptr   B   - [\alpha, \Ad_{x\inv } {\bs}^* B ]=0$. 
According to Lemma \ref{lem:com}, we have $ \ptr  \Omega =
\ptr  d^{\nabla} B = 0$.
This concludes the proof. 
\end{proof}

\begin{defn}
A Lie groupoid extension $X_1\xrightarrow{\phi}Y_1\toto M$
is called \textbf{flat} if there exists a right connection
1-form $\alpha \in \Omega^1 (X_1, \bt^* \liekernel  )$, a curving $B\in
\Omega^2 (M,\liekernel)$ such that the 
induced connection on the group
bundle $\kernel\xrightarrow{\pi}M$ is flat, 
and the 3-curvature $\Omega$ vanishes.
\end{defn}

Recall that, from Corollary \ref{cor:4.46},
 the  band is flat if and only if $ [\omega^\kernel]=0$ in 
$\hcoho{1}{2}{\kernel}{M}$.
Hence when a Lie  groupoid extension is flat, its band
must be flat. In fact, requiring band to be flat
is slightly weaker than requiring the group bundle $\kernel\xrightarrow{\pi}M$ to be flat.

The following proposition is obvious.

\begin{prop}
    Given a  connection on a Lie  groupoid $G$-extension
 $X_1\xrightarrow{\phi}Y_1\toto M$, where $M$ is a disjoint union of contractible manifolds,
then the following are equivalent
\begin{itemize}
\item   the group bundle is flat, i.e. ${\omegakernel}=0 $;
\item  there exists a trivialization $\chi : \kernel  \isomorphism  M\times G $
such that $ F_g =0  $ for all $g \in G$.
\end{itemize}
\end{prop}
%\begin{proof}
%Since  $M$ is a disjoint union of contractible manifolds,
%the group bundle is flat if and only if it admits horizontal sections.
%\end{proof}

Let $X_1\xrightarrow{\phi}Y_1\toto M$ be a $G$-extension with
 kernel $\kernel\xrightarrow{\pi}M$.
Denote 
 by $Z(\kernel)\xrightarrow{\pi}M$ its bundle of centers.
Assume that  $G$ 
 is a connected Lie group with compact center $Z(G)$,
i.e.    $Z(G)$ is isomorphic to a tori.
Hence the automorphism group $\Aut (Z(G))$ is 
discrete. It thus follows that, for any connection on the extension
$X_1\xrightarrow{\phi}Y_1\toto M $, the induced connection
on the kernel $\kernel\xrightarrow{\phi}M$ defines a flat
connection on this subgroup bundle $Z(\kernel)\xrightarrow{\phi}M$.  
Therefore, the induced connection   on the Lie algebra bundle  ${\liecenter} \to M$ is  also flat, which implies that
its pull-back connection on $\bt^* {\liecenter}  \to Y_n$, 
$\forall n\in\N$, is  flat as well.
 By   $${\covder}:  \Omega^k (Y_n, \bt^* {\liecenter} )
\lto \Omega^{k+1} (Y_n, \bt^* {\liecenter} ), $$
we denote its exterior  covariant differential.

It is simple to see that the adjoint action of
 $X_1 \toto M$ on $\liekernel \to M$ induces an action of $Y_1 \toto M$
 on $\liecenter \to M $. Thus we have a differential
\begin{equation}
\label{eq:par}
{\pt} :\Omega^k \big( Y_n ,\bt^* {\liecenter} \big) 
\lto  \Omega^k \big(Y_{n+1} , \bt^* {\liecenter} \big) .
\end{equation}

\begin{lem}
\label{lem:commut}
 The following relations hold:
\[ \pt^2=0,  \qquad  \covder^2=0,   \qquad  [\pt,\covder]=0 .\]
\end{lem}

\begin{proof}
The first relation is a general fact
for any Lie groupoid representation.
Let us prove the second one.

Since the fibers of $Z(\kernel) \to M $ are tori, for any
 $m \in M$,  the inverse image of
the exponential map
$exp: \liecenter_m \to Z(\kernel)_m $ of the identity element in
$Z(\kernel)_m $ 
 is a lattice $L_m\subset\liecenter_m$ of maximal rank.
This lattice is smooth, i.e. it is generated by smooth sections. 
 The identity section in $Z(\kernel)\to M$ is horizontal and 
 by Lemma~\ref{lem:expopreserveshoriz}, 
a section of $\liecenter\to M$ is horizontal if and only if its image 
under  the 
exponential map is horizontal.
Hence it follows that 
 any smooth section of the lattice $L\to M$ is horizontal.
As a consequence, the pull back $\bt^* l$ of a section $l$ 
of the lattice $L$ is a parallel section of $\bt^* \liecenter\to Y_n$.

This fact allows us to give an explicit local expression of $D$ with the help of
sections of the lattice $L \to M$ as follows. 
Choose  $ l_1,\cdots, l_k $ local smooth generators of the lattice $L \to M$
on some open set $U \subset M$. Any
 $\omega \in \Omega^m(\bt\inv (U),\bt^* {\liecenter}) $
can then be written uniquely as
 $ \omega = \sum_{i=1}^k \omega_i \otimes \bt^* l_i$
(where $\omega_1,\cdots, \omega_k $ are $ m$-forms on  $\bt\inv (U) $ and 
$\bt^* l_i $ denotes the pull-back of $l_i$ through $\bt$, for all 
$i \in \{1,\cdots,k\}$).
Hence, the covariant differential of $ \omega$ is given by
\begin{equation}\label{eq:covederZ}     {\covder}     \omega = \sum_{i=1}^k (d \omega_i ) \otimes \bt^* l_i   , \end{equation}
where $d$  stands for the usual de Rham differential.
Hence it follows that $D^2=0$.

Next we give an explicit  expression of $\pt$ with the help of sections of the lattice $L \to M$.
Consider any point  $(y_1,\cdots,y_n) \in Y_n$. Let $(l_1,\cdots, l_k) $ 
be smooth generators of the lattice $L $ in a neighborhood $U$
of $ {\bt}(y_n)$ and $l_1',\cdots, l_k' $  smooth generators in a neighborhood
$V$ of $ {\bs}(y_n)$.
The adjoint action of any $y\in {\bs}\inv (V) \cap {\bt}\inv (U)  \subset  Y_1 $ must preserve the lattice $L$. That is,
 for any $ l \in L_{{\bs}(y)} $, we have $(\Ad_y)\inv  l \in  L_{\bt(y)}$.
%%and, as a consequence, we also have $(\mbox{Ad}_{y_n})\inv  l \in  {\bt}^* L_{(y_1,\cdots,y_n)}$
%for any $ l \in {\bt}^* L_{(y_1,\cdots,y_n)} $, 
As a consequence, there exists a matrix $(M_a^j)_{a,j=1}^k$ with constant integral coefficients
  such that, on ${\bs}\inv (V) \cap {\bt}\inv (U)$,
    $ \Ad_{y\inv } l_j'  = \sum_{a=1}^k M_a^ j   l_j$.
Taking the pull-back under $\bt$, one obtains that    
  \begin{equation}\label{eq:ReprderZ}  \Ad_{y_n\inv } \bt^* l_j'  = \sum_{a=1}^k M_a^ j  \bt^* l_j .\end{equation}
The restrictions of $\omega \in \Omega^m (Y_n,\bt^*\liecenter)$
to $ \bt\inv (U)$ and $\bt\inv  (V) $,
can be written in an unique way as $   \omega =  \sum_{j=1}^k \omega_j  \otimes \bt^ * l_j  $ and
$    \omega =  \sum_{j=1}^k \omega_j ' \otimes \bt^ * l_j'  $,
respectively, where $\omega_1,\cdots,\omega_k \in
\Omega^m\big(\bt\inv (U)\big)$  and $ \omega_1',\cdots, \omega_i ' \in
\Omega^m\big(\bt\inv (V)\big) $.
By the  definition of $\pt$,
 for any $(y_1,\cdots,y_n) \in Y_n$ 
with $y_n \in {\bs}\inv (V) \cap {\bt}\inv (U) $, we have: 
\begin{equation} \label{eq:ReprderZ2}
\begin{aligned} \pt\omega =& \sum_{i=0}^{n-1}(-1)^i(\epsilon^n_i)^*\omega+(-1)^n\Ad\inv_{y_n}(\epsilon^n_n)^*\omega \\ 
=& \sum_{i=0}^{n-1}\sum_{j=1}^{k}(-1)^{i}(\epsilon^n_i)^*\omega_j\otimes\bt^* l_j+(-1)^{n}\sum_{j=1}^{k}(\epsilon^n_n)^*\omega_j'
\otimes\Ad_{y_n}\inv\bt^*l_i' \\
=& \sum_{i=0}^{n-1}\sum_{j=1}^{k}(-1)^{i}(\epsilon_i^n)^*\omega_j\otimes\bt^* l_j+(-1)^{n}\sum_{j=1}^{k}\sum_{a=1}^{k}
 M_a^j  (\epsilon_n^n)^*   \omega_j'  \otimes \bt^* l_j .
\end{aligned} \end{equation}
The relation $\big[\pt,\covder\big]=0$ thus follows 
immediately from Eqs.~\eqref{eq:covederZ}-\eqref{eq:ReprderZ2}
together with the fact that the de Rham differential commutes
with the pull-back maps
 $(\epsilon_i^n)^* $ for all $i\in\{0,\cdots,n\}$.
\end{proof}

If $Y_1\toto M$ is a \v{C}ech groupoid 
associated to an open covering of a  manifold $N$,
the Lie algebra  bundle $Z(\liekernel) \to M$ 
over $Y_1\toto M$ corresponds to  a Lie algebra
bundle over $N$, denoted by  $Z(\liekernel)_N \to N$.

\begin{lem}
\label{lem:par}
Assume that   $Y_1\toto M$ is a  \v{C}ech groupoid
 associated to an open covering of a  manifold $N$.
Then the cohomology of the cochain complex
\eqref{eq:par}
is given by
\begin{equation}
H^n ( \Omega^k ( Y\simplicial ,\bt^* {\liecenter} ) , \pt )
= \left\{  \begin{array}{ccc}
\Omega^k (N, Z (\liekernel)_N)  &  & n=0  \\
  0 &  & n>0 
  \end{array} \right. \end{equation}
\end{lem}
\begin{proof}
$Y_1\toto M$ is Morita equivalent to
$N\toto N$, and under this Morita equivalence,
the module $Z(\liekernel) \to M$ over $Y_1\toto M$ 
becomes the trivial one $Z (\liekernel)_N\to N$.
Thus the conclusion  follows from a  general fact
regarding  groupoid cohomology  of Morita equivalent
groupoids \cite{Crainic}.
\end{proof}

Now we are ready to  state the main result of this section.

\begin{thm} 
Let $X_1\xrightarrow{\phi}Y_1\toto M$ be a $G$-extension with
 kernel  $\kernel\xrightarrow{\phi}M$, and 
 $\alpha \in \Omega^1 (X_1, \bt^* \liekernel  )$ a right connection 1-form
  with Ehresmann curvature $\omega \in \Omega^2 (X_1, \bt^* \liekernel  )$.
Assume that $\omega^\kernel=0$, i.e. the group bundle $\kernel\xrightarrow{\phi}M$ is flat. Then
  \begin{enumerate}
\item $\omega $ is  valued in the center  ${\bt}^*{\liecenter} $,
i.e. $\omega\in \Omega^2 (X_1, {\bt}^*{\liecenter} )$;
\item  $\covderPB \omega = 0$;
\item there exists a $\eta \in \Omega^2 \big( Y_1,  \mathfrak
  {t}^* \liecenter\big) $ such that  $  \phi^* \eta =\omega $,
and  satisfies
$$  {\pt} \eta =0, \ \ \ \covder \eta =0;$$
\item if $Y_1 \toto M $ is the \v{C}ech groupoid
associated to an open  cover of a manifold $N$, then
there exists a curving $B \in \Omega^2(M,\liecenter)$,
 and the $3$-curvature $\Omega $ descends to
a  $3$-form $ \Omega_N$ in $ \Omega^3 \big(N,\liecenter_N\big)$
  such that $\pi^* \Omega_N =\Omega$, where $\pi: M\to N$ is the
projection.
 \end{enumerate}
\end{thm}
\begin{proof}
\begin{itemize} 
\item[1.] According to Eq.~\eqref{eq:delalpha2},
if ${\omegakernel}=0 $, then $\delta_\kernel \omega=0$.
Hence  $\omega $ must take its values in the center of the Lie algebra
 $\liekernel $.
\item[2.] It follows from 1) and the Bianchi identity \eqref{eq:Bianchi1}.
\item[3.] One needs to  prove that $\omega $ is basic with respect to the right
  action of $\kernel$.
Note that $\omega $ is an horizontal 2-form according
 to Proposition \ref{ex:hor1form}.
For any $x \in X_1$ and $k  \in \kernel_{{\bt}(x)} $, 
let $\sigma $ be  any local section
of $\kernel \to M$  through $k$.
For  any $v_i \in T_x X_1, \ i=1, 2$,  the elements
 $  u_i = \sigma_*\smalcirc  {\bs}_* (v_i) $ and $ v_i $, $i=1,2 $ 
are composable  in the tangent groupoid $TX_1 \toto TM$, and
  $ u_i \cdot v_i = (R_\sigma)_* u_i$, $i=1,2 $.
Then the relation
$\ptr  \omega=0 $  implies  that
\begin{multline*} \omega_{|_{k \cdot x}} \big( (R_\sigma)_* v_1,(R_\sigma)_*  v_2 \big)
= \omega_{|_{k \cdot x}} \big( u_1 \cdot v_1, u_2 \cdot  v_2 \big) \\
= \omega_{|_{k}}(u_1,u_2) + \Ad_k\inv  \cdot \omega_{|_{ x}}(v_1,v_2)
= \omega_{|_{ x}}(v_1,v_2), \end{multline*} 
where, in the second equality, we used the assumption that $\omega^\kernel=0$.
Thus it follows that  $\omega $ is basic. Hence there exists a
$\eta\in\Omega^2\big(Y_1,\bt^*\liecenter\big)$
such that $\omega=\phi^*\eta$.

Since the groupoid morphism  $ \phi: X_{\simplicial} \to Y_{\simplicial}$
commutes with their  actions  on
$\liecenter\xrightarrow{\phi}M$, it follows that
$ \phi^* {\pt} = \ptr   \phi^*$.
Thus we have  ${\pt}\eta =0$ by Eq.~\eqref{eq:Bianchi2}.
On the other hand, we have
$0=\covderPB \omega =\covderPB  \phi^*\eta =\phi^* D\eta$.
It thus follows that $D\eta=0$.
\item[4.] It follows from 3) together with Lemma \ref{lem:par}.
\end{itemize}
\end{proof}

\begin{rmk}
It would be interesting to compare our definition
of curving and $3$-curvature with the one defined by Breen-Messing
\cite{BreenMessing}. It is clear that 
our definition reduces to the standard one 
in \cite{Murray, Hitchin} for
bundle gerbes.
\end{rmk}

\subsection{Connections on $G$-central extensions}
\label{sec:6.6}

We now study connections and curvings  on Lie groupoid
$G$-central extensions.
By Theorem \ref{thm:3.9}, for any central $G$-extension
$X_1\xrightarrow{\phi}Y_1 \toto M$,
there corresponds a $Z(G)$-central extension 
$\tilde{X}_1\to Y_1\toto M$ so that
 \begin{equation}\label{eq:rappel}
 X_1  \isomorphism \frac{ \tilde{X}_1 \times G}{Z(G)},
 \end{equation}
where  $Z(G)$
 acts  on $\tilde{X}_1\times G$ diagonally:
$(\tilde{x},g)\cdot h=(\tilde{x}h_{\bt(x)},h\inv g)$, 
$\forall h\in Z(G)$.
Recall that $X_1$ is a $G$-$G$ principal bibundle, where
the actions are given, respectively, by
$$g\cdot[x,g']=[x,gg'], \qquad [x,g']g=[x,g'g].$$

Denote by $\pi:\tilde{X}_1\times G\to X_1$ the quotient map,
by $\pr_1:\tilde{X}_1\times G\to\tilde{X}_1$ 
and $\pr_2:\tilde{X}_1\times G\to G $
 the projections.
Also, denote by $\tau:\tilde{X_1}\to X_1$, the 
embedding defined by $x\to\pi(x,1)$.

The following result describes the precise  relation between
connections on  these extensions. 
Note that in this case  $\Omega^m(X_n,\bt^*\liekernel)$
(resp. $\Omega^m(X_n,\bt^*\liecenter)$) can be  
naturally identified with
the space of $\mfg$ (resp. $Z(\mfg)$)-valued $m$-forms
on $X_n$. Denote by $\theta$ the left Maurer-Cartan form on $G$.

\begin{thm}
 \label{prop:central11corr}
For a groupoid central $G$-extension
 $X_1\xrightarrow{\phi}Y_1\toto M$,
there is a one-one correspondence between right  
connection 1-forms $\alpha\in\Omega^1(X_1,\mfg)$
  which
are also  right principal  $G$-bundle connections,
 and connections $\tilde{\alpha}\in\Omega^1(\tilde{X}_1,Z(\mfg))$
  of the groupoid central $Z(G)$-extension
 $\tilde{X}_1\xrightarrow{\phi}Y_1\toto M$.
\end{thm}

\begin{proof}
Assume that $\alpha\in\Omega^1(X_1,\mfg)$
is a connection 1-form for the right principal
$G$-bundle $X_1\to Y_1$.
 It is simple to see that
$\tau^*\alpha$
must be $Z(\mfg)$-valued and $\tilde{\alpha}= \tau^*\alpha
\in \Omega^1 (\tilde{X}_1, Z(\mfg)) $  is 
a connection 1-form for the $Z(G)$-principal
bundle $\tilde{X}_1\to Y_1$. Conversely, given
a connection 1-form  $\tilde{\alpha} \in \Omega^1
 (\tilde{X}_1, Z(\mfg))$ for the $Z(G)$-principal
bundle $\tilde{X}_1\to Y_1$, $\pr_1^* \tilde{\alpha}+ \pr_2^* \theta
\in  \Omega^1   (\tilde{X}_1\times G,  \mfg)$ is basic
 with respect   to the diagonal $Z(G)$-action. Hence
it defines a   1-form $\alpha\in \Omega^1 (X_1,  \mfg)$
such that  $\pi^* \alpha = \pr_1^* \tilde{\alpha}+ \pr_2^* \theta$.
It is simple to check directly that $\alpha$
is  a connection 1-form for the right principal
$G$-bundle $X_1\to Y_1$.

Now we have
$$\pi^* \ptl \alpha
= \ptl ( \pr_1^* \tilde{\alpha}+ \pr_2^* \theta)
=\pr_1^* \ptl \tilde{\alpha}+
\pr_2^* \ptl \theta=\pr_1^* \ptl \tilde{\alpha}, $$
since $\theta$ is a left Maurer-Cartan form. Hence it follows
that $\ptl \alpha=0$ if and only if
 $\ptl \tilde{\alpha}=0$. 
Therefore the conclusion follows.
\end{proof}

We end this section with the following important

\begin{cor}
\label{cor:centralrestriction2}
\begin{enumerate}
\item There is a one-one correspondence between
connections $\alpha$ on the  central $G$-extension
 $X_1\to Y_1\toto M$ such that $M\times \{g\}$,
for all $g\in G$, are horizontal, and
connections $\tilde{\alpha}$ on the central ${Z(G)}$-extension
$\tilde{X_1}\to Y_1\toto M$.
\item The form $B\in\Omega^2(M,Z(\mfg))$ is a curving for
  $\alpha$ if and only if
it is a curving for $\tilde{\alpha}$. In this
case, their 3-curvatures coincide.
\item There is a one-one correspondence between flat central $G$-extensions
$X_1\to Y_1\toto M$ such that $M\times \{g\}$,
for all $g\in G$,  are  horizontal, and flat  central $Z(G)$-extensions
$\tilde{X_1}\to Y_1\toto M$.
\end{enumerate}
\end{cor}

\end{document}